\documentclass[12pt]{article}
\usepackage{ct}

\def\+{\includegraphics[scale=0.4]{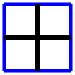}}
\def\elbow{\includegraphics[scale=0.4]{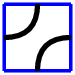}}

\DeclareMathOperator{\Inv}{Inv}

\usepackage[dvipsnames]{xcolor}
\usepackage{amsfonts,enumitem}
\usepackage{nicematrix}
\usepackage{tikz-cd}
\usepackage{url}   
\usepackage{stackengine}
\usepackage{tikz}

\usepackage{environ}
\numberwithin{equation}{section}

\newcommand\sbullet[1][.5]{\mathbin{\vcenter{\hbox{\scalebox{#1}{$\bullet$}}}}}

\newcommand{\rect}{{\textnormal{Rect}}}
\newcommand{\prtnshape}{\textnormal{Shape}}
\newcommand{\maxrow}{{\textnormal{max\_bump\_row}}}
\newcommand{\mincol}{{\textnormal{min\_bump\_col}}}
\newcommand{\lpath}{{\textnormal{Path}}}
\newcommand{\rpd}{{\textnormal{RPD}}}
\newcommand{\B}{{B_x(w)}}

\newcommand{\ladder}{{\mathcal{L}}}
\newcommand{\move}{{\mathcal{M}}}
\newcommand{\upperbump}{h}
\newcommand{\rightbump}{k}
\newcommand{\T}{\mathbf{T}}

\DeclareFontFamily{U}{mathx}{\hyphenchar\font45}
\DeclareFontShape{U}{mathx}{m}{n}{%
<-6> mathx5
<6-7> mathx6
<7-8> mathx7
<8-9> mathx8
<9-10> mathx9
<10-12> mathx10
<12-> mathx12
}{}

\DeclareSymbolFont{mathx}{U}{mathx}{m}{n}
\DeclareFontSubstitution{U}{mathx}{m}{n}

\DeclareMathSymbol{\bigovoid}{\mathop}{mathx}{"EC}

\newcommand{\BigO}{\mathop{\stackinset{c}{}{c}{}{ \scalebox{1.1}{$\bigovoid$}}{ \scalebox{1.15}{$\bigovoid$}}}}

\newcommand{\bigO}{\mathop{\stackinset{c}{-2pt}{c}{}{ \scalebox{0.8}{$\bigovoid$}}{\scalebox{0.85}{$\bigovoid$}}}}

\newcommand{\given}{\,\, | \,\,}
\newcommand{\Dbot}{D_{\textrm{bot}}}
\newcommand{\Dtop}{D_{\textrm{top}}}

\specs{n}{vol (issue)}{year}{}

\dateline{Oct 14, 2025}{Aug 26, 2026}{TBD}



\title{A Proof of Rubey's Lattice Conjecture}

\author[1]{Sara Billey}
\author[2]{Connor McCausland}
\author[3]{Clare Minnerath}

\affil[1,2,3]{%
Department of Mathematics, University of Washington, Seattle

\email{billey@uw.edu}, \email{connortm@uw.edu}, \email{cminnera@uw.edu}%
}

\begin{document}
\maketitle
\begin{abstract}
In 2011, Rubey generalized chute and ladder moves on the set of
reduced pipe dreams for a permutation $w$ and conjectured that
the induced poset on reduced pipe dreams is a lattice. In this paper,
we prove this conjecture. Our key tool is a new type of move operation
$\mathcal{M}_{ij}$, defined as a composite of certain generalized
ladder moves and then simplified in terms of swaps
on a partition shape. We show that joins and meets exist in Rubey's
poset by proving simple recursive formulas in terms of
$\mathcal{M}_{ij}$ operations, so we refer to the poset as the Rubey lattice henceforth. In addition, we give an explicit
criterion to determine if two elements of the Rubey lattice are comparable using an injective map from reduced pipe dreams to tableaux on the
diagram of a permutation. The pipe dream tableaux construction also gives an exact formula for the maximum length of any chain in the Rubey lattice and bounds on the number of reduced pipe dreams for $w$, or equivalently the sum of the coefficients in the corresponding Schubert polynomial. Several open problems for further development of the Rubey lattice related to Markov processes, Schubert polynomials, and Grothendieck polynomials are given.
\end{abstract}

\section{Introduction}

\begin{figure}[h]
\begin{center}
\includegraphics[scale=0.15]{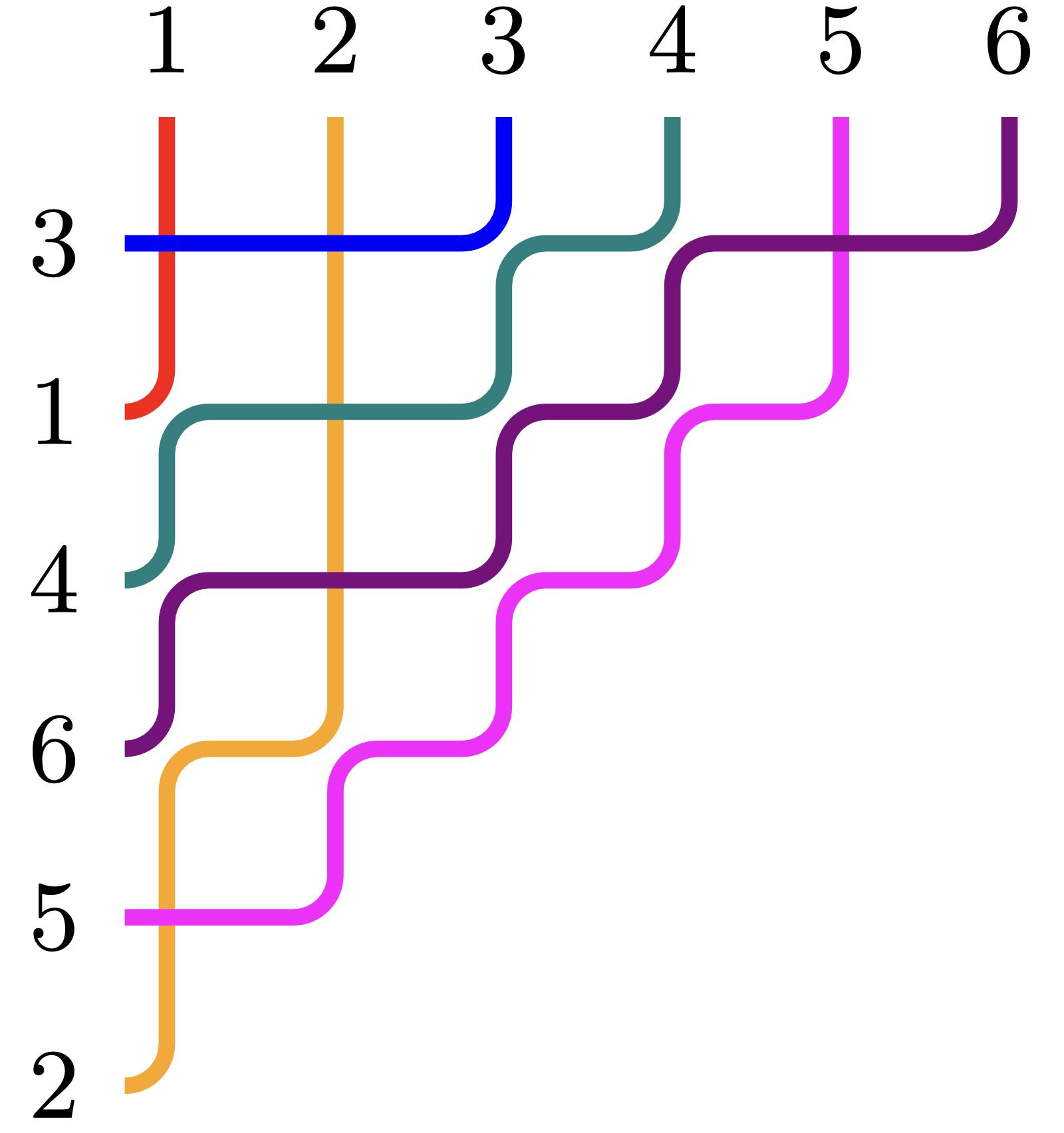}
\end{center}
\caption{A reduced pipe dream $D$ for $w=[3,1,4,6,5,2]$.}
\label{fig:rcgraphs}
\end{figure}

Reduced pipe dreams are combinatorial objects that encode some of the
algebraic, enumerative, geometric, and probabilistic properties
related to Schubert  and Grothendieck polynomials
\cite{billeygaopow2025introductioncohomologyflagvariety,knutson-miller-2005,kogan.phd,MPP.2019,MPPY.2025,anderson2026computationsamplingschubertspecializations}.
See \Cref{fig:rcgraphs} for an example, and see \Cref{sec:background}
for the definition as a union of cross and bump tiles forming the
pipes. They first appeared in the work of Billey-Bergeron
\cite{billey-bergeron} following the work of several authors \cite{b-j-s,FS,FK,LS1,M2}. They showed that there are natural raising and lowering operators, called chutes and ladders,  connecting the
set of all reduced pipe dreams for a fixed permutation.

In \cite{rubey2012maximal}, Rubey generalized the chute and
ladder moves. He conjectured that the poset on all reduced pipe
dreams for any fixed permutation with covering relations given by the
generalized chute moves, or equivalently generalized ladder moves, is a lattice. See Section 2 for the relevant definitions and notation. The main theorem of this paper
is a proof of Rubey's conjecture.

\begin{theorem}\label{thm:rubey.lattice}
For any permutation $w \in S_{n}$, the generalized chute poset on its reduced
pipe dreams is a lattice.
\end{theorem}

Henceforth, we will refer to this poset as the \textit{Rubey lattice} for $w$, denoted $(\rpd(w),<)$. We note that another independent proof
of Rubey's conjecture has recently been found \cite{axelrodfreed2025chuteposetslattices}. An abstract written by both teams and summarizing both proofs has been accepted to FPSAC 2026 \cite{ABDMMMNT26}. Additionally, a special case of the conjecture was previously known by a connection to $\nu$-Tamari lattices \cite{ceballos2019nutamarilatticenutreesnubracket}. Another approach to finding lattices on pipe dreams was given in \cite{BCCP.2025}.

The key tool in our proof of \Cref{thm:rubey.lattice} is both an
algorithmic and explicit realization of a new type of move operation
$\move_{ij}$, defined on any $D \in \rpd(w)$ whenever there is a cross at $(i,j)$ in $D$ and a
bump tile above it. See \Cref{nonrecursive_m} and \Cref{fig:oneshape}
for an explicit description of $\move_{ij}(D)$. These move operators
allow us to identify the unique minimal element above $D$ in $(\rpd(w),<)$ with no cross at $(i,j)$ that can be reached by any sequence of ladder moves
changing positions northeast of $(i,j)$. The operator $\move_{ij}$
swaps cross tiles and bump tiles along the boundary of an easily
identified partition shape with lower left corner at
$(i,j)$. Successive applications of the move operators are used in the
algorithm to compute the join and meet of two pipe dreams. Therefore,
the move operators also give a method to recursively verify if $D_{1}$
and $D_{2}$ are comparable in the Rubey lattice. We show that
$\move_{ij}\circ \move_{kl}(D)=\move_{kl}\circ \move_{ij}(D)$ whenever
$(i,j)$ and $(k,l)$ are southwest incomparable. Therefore, to
state the test, we will use the \textit{big composition} notation, meaning

\[
\BigO_{(i,j) \in S}\move_{ij}(D) = \move_{i_{1},j_{1}}\circ
\move_{i_{2},j_{2}}\circ \cdots \circ \move_{i_{m},j_{m}}(D)
\]
whenever $S=\{(i_{1},j_{1}),\dots, (i_{m},j_{m}) \}$ is a set of
southwest incomparable coordinates so the corresponding move operators
all commute. 

\begin{corollary}\label{cor:comprability.test}
Let $D_{1}$ and $D_{2}$ be reduced pipe dreams for $w \in S_{n}$. Let $S$ be the set of southwest extremal coordinates where their cross/bump tiles
differ. Then $D_{1} \leq D_{2}$ if and only if 
\begin{enumerate}
\item $D_{1}$ contains a cross at each element of $S$, and
\item $D'=\BigO_{(i,j) \in S}\move_{ij}(D_{1}) \leq D_{2}$, which can be
further tested recursively.
\end{enumerate}
\end{corollary}

In addition to the recursive method above, there is another simple method to
verify if $D_{1}\le D_{2}$ for any two pipe dreams for the same
permutation in the Rubey lattice. In \Cref{sec:comparison}, we will
define a tableau $\T(D)$ associated with the reduced pipe dream $D$ to
be a certain filling of the diagram of the corresponding permutation
by positive integers. Such \textit{pipe dream tableaux} naturally sit inside the lattice
on $\mathbb{Z}_{+}^{\ell}$ under coordinate-wise order, denoted
$\unlhd$. We note that a similar notion called inversions tableaux or flagged inversion fillings has been used in \cite{axelrodfreed2025inversionstableaux,Kelly.thesis}.

\begin{theorem}\label{thm:comprable}
Let $w \in S_{n}$. If $D_1$ and $D_{2}$ are reduced pipe dreams for
$w$, then $$D_1\le D_2 \iff \T(D_1)\unlhd \T(D_2).$$
\end{theorem}

Several corollaries follow from \Cref{thm:comprable}. In particular, the number of reduced pipe dreams for $w$ can be bounded above and below by considering the minimal and maximal pipe dream tableaux for $w$, see \Cref{cor:bounds}. Recall that the Schubert polynomial $\mathfrak{S}_w
(x_1,\ldots, x_n)$ is a generating function over the reduced pipe dreams of $w$, so these bounds also have implications for finding the sum of the coefficients in the Schubert polynomial associated to $w$, which have been studied before in \cite{stanley2017schubertshenanigans,Weigandt.2018,Gao.2021,MPPY.2025,anderson2026computationsamplingschubertspecializations}. Furthermore, an exact formula for the maximal length of any chain in the Rubey lattice $(\rpd(w),<)$ is given in \Cref{cor:max.chain.length}.

The outline of this paper is as follows. In \Cref{sec:background}, we
review the necessary concepts for pipe dreams from the literature and
Rubey's definition of the poset on reduced pipe dreams. In
\Cref{sec:move}, we define the move operator algorithmically in terms
of generalized ladder moves and explicitly in terms of swapping bump
and cross tiles along the boundary of a corresponding partition shape. We then
analyze the commutativity properties of these operators. In
\Cref{sec:proof}, we state and prove an algorithm for the join of any
two pipe dreams for a fixed permutation. The meet operation is given
by transposition symmetry on pipe dreams. With these join and meet
operations, we complete the proof of \Cref{thm:rubey.lattice} at the
end of \Cref{sec:proof}. In \Cref{sec:comparison}, we define pipe dream tableaux, which are used to prove \Cref{thm:comprable} and describe certain properties of the joins and meets of pipe dreams. We conclude with some further observations and open problems on the Rubey lattice in \Cref{sec:future}, including characterizing pipe dream tableaux explicitly, characterizing the join algorithm for the Rubey lattice in terms of pipe dream tableaux, identifying a Markov process on $\rpd(w)$ that converges to the uniform distribution, and connections to Stanley's conjecture on the maximal number of reduced pipe dreams for any $w \in S_n$.

\section{Background}\label{sec:background}

Let $S_{n}$ be the symmetric group of permutations on $[n]=\{1,2,\dots
,n \}$ for any positive integer $n$. We will denote a permutation $w
\in S_{n}$ by its one-line notation $w=w_{1}w_{2}\dots w_{n}$,
listing the elements in $[n]$ in some order. Each permutation encodes
a bijection from $[n]$ to itself given by $w(i)=w_{i}$.
Multiplication in $S_{n}$ is defined by composition of bijections,
$vw=v_{w_1}\dots v_{w_n}$. The set of \textit{inversions} for $w$
is
\[
\Inv (w) =\{(i,j) \in [n]^{2} \given i<j, \ w_{i}>w_{j} \}.
\]
The group $S_{n}$ naturally embeds into $S_{n+1}$ by adding $n+1$ as a
fixed point to each bijection. This embedding preserves inversion
sets. Let $S_{\infty}=\bigcup_{n\geq 1} S_{n}$ be the set of
permutations on the positive integers $\mathbb{Z}_+$ that
fix all but finitely many values, where we identify two permutations
if they have the same inversion set.

The \textit{diagram of a permutation} $w \in S_{n}$ is obtained from
the $n\times n$ array with a dot in position $(i,w_{i})$ for $1\leq i\leq n$ by removing all cells in the array which are in the hook shape weakly to the
right of or below a dot, including the dots. The remaining cells
form the (Rothe) diagram $\textbf{D}(w)$, see \Cref{fig:diagram.w}. Thus,
\begin{equation}\label{eq:diagram}
 \textbf{D}(w)= \{(i,j)\in [n]^{2} \given j< w(i) \text{ and } i< w^{-1}(j) \}.
\end{equation}
Equivalently, the cells of $\textbf{D}(w)$ are in bijection with the inversions
of $w$, so
\begin{equation}\label{eq:diagram.2}
 \textbf{D}(w)= \{(i,w(j))\in [n]^{2} \given i<j \text{ and } w(i)> w(j) \}.
\end{equation}

\begin{figure}[h!]
\begin{center}
\begin{minipage}{0.45\linewidth}
\centering
\begin{tikzpicture}[scale=0.6]
\draw[step=1.0,green,thin] (0,0) grid (5,5);
\draw[very thick] (1,1)--(2,1)--(2,2)--(1,2)--(1,1);
\draw[very thick] (0,3)--(2,3)--(2,4)--(3,4)--(3,5)--(0,5)--(0,3);
\node at (0.5,2.5) {$\bullet$};
\node at (1.5,0.5) {$\bullet$};
\node at (2.5,3.5) {$\bullet$};
\node at (3.5,4.5) {$\bullet$};
\node at (4.5,1.5) {$\bullet$};
\draw(0.5,-0.5)--(0.5,2.5)--(5.5,2.5);
\draw(1.5,-0.5)--(1.5,0.5)--(5.5,0.5);
\draw(2.5,-0.5)--(2.5,3.5)--(5.5,3.5);
\draw(3.5,-0.5)--(3.5,4.5)--(5.5,4.5);
\draw(4.5,-0.5)--(4.5,1.5)--(5.5,1.5);
\node at (1.5,1.5) {};
\node at (0.5,3.5) {};
\node at (0.5,4.5) {};
\node at (1.5,3.5) {};
\node at (1.5,4.5) {};
\node at (2.5,4.5) {};
\end{tikzpicture}
\end{minipage}
\end{center}
\caption{The diagram of $w=43152$ is the set of outlined cells above,
so 
$\textbf{D}(w)=\{(1,1), (1,2), (1,3), (2,1), (2,2), (4,2) \}$ using matrix coordinates. }
\label{fig:diagram.w}
\end{figure}

A \emph{pipe dream} $D$ is a finite subset of $\mathbb{Z}_+\times
\mathbb{Z}_+$. One can visualize a pipe dream by placing a \+-tile, called
a ``cross", at every $(i,j)\in D$ using matrix coordinates, and a
\elbow-tile, called a ``bump" or an ``elbow", at every $(i,j)\in
\mathbb{Z}_+\times\mathbb{Z}_+\setminus D$, creating pipes (or wires)
connecting the top boundary to the left boundary. If the pipes are
numbered $1, 2, 3, \ldots$ across the top, then the corresponding
pipes read along the left side of the diagram from top to bottom form
a permutation $w \in S_{\infty}$. Call $w$ the \emph{permutation of $D$} following the
literature.

A pipe dream $D$ is \emph{reduced} if $w$ is the permutation of $D$ and
$|D|$ equals the number of inversions of $w$. Write $\rpd(w)$ for
the set of all reduced pipe dreams for $w$. Observe that two pipes
labeled $x<y$ cross somewhere in a pipe dream $D \in \rpd(w)$ if and
only if $(x,y) \in \Inv(w^{-1})$. The larger labeled pipe necessarily
enters the crossing horizontally in a reduced pipe dream as it
proceeds down and left.

One only needs to draw a finite number of pipes in a triangular array to
represent a pipe dream since it necessarily contains only a finite
number of crossings. Furthermore, we often simplify the pictures by
drawing a matrix with $+$ representing a cross tile and $\sbullet$
representing a bump tile. The pipes can easily be reconstructed from
this data. We abuse notation and write $D(i,j)=+$ to mean $(i,j) \in
D$, and $D(i,j)=\sbullet$ to mean $(i,j) \not\in D$ to better connect the
proofs to the pictures. Also, we often show a subarray of a pipe
dream when discussing local properties.

\begin{figure}[h!]
\centering
\begin{tikzpicture}[scale=0.4]
\node at (0,0) {$\cdot$};
\node at (1,0) {$+$};
\node at (2,0) {$+$};
\node at (3,0) {$+$};
\node at (4,0) {$\cdot$};
\node at (0,-1) {$+$};
\node at (1,-1) {$+$};
\node at (2,-1) {$+$};
\node at (3,-1) {$+$};
\node at (4,-1) {$\cdot$};
\def\a{-8};
\node at (\a+0,0) {$\cdot$};
\node at (\a+1,0) {$+$};
\node at (\a+2,0) {$+$};
\node at (\a+3,0) {$+$};
\node at (\a+4,0) {$+$};
\node at (\a+0,-1) {$\cdot$};
\node at (\a+1,-1) {$+$};
\node at (\a+2,-1) {$+$};
\node at (\a+3,-1) {$+$};
\node at (\a+4,-1) {$\cdot$};
\draw[very thick,->] (-3,-0.5)--(-1,-0.5);
\end{tikzpicture}
\qquad
\begin{tikzpicture}[scale=0.2]
\def\a{0};
\def\b{0};
\draw[step=2, color=green](\a+2,-6+\b)grid(\a+12,-2+\b);
\draw[very thick](\a+3,-2+\b)--(\a+3,-2.5+\b);
\draw[very thick](\a+3,-3.5+\b)--(\a+3,-4+\b);
\draw[very thick](\a+2,-3+\b)--(\a+2.5,-3+\b);
\draw[very thick](\a+3.5,-3+\b)--(\a+4,-3+\b);
\draw[very thick](\a+2.5,-3+\b)arc(270:360:0.5);
\draw[very thick](\a+3.5,-3+\b)arc(90:180:0.5);
\draw[very thick](\a+4,-3+\b)--(\a+6,-3+\b);
\draw[very thick](\a+5,-2+\b)--(\a+5,-4+\b);
\draw[very thick](\a+6,-3+\b)--(\a+8,-3+\b);
\draw[very thick](\a+7,-2+\b)--(\a+7,-4+\b);
\draw[very thick](\a+8,-3+\b)--(\a+10,-3+\b);
\draw[very thick](\a+9,-2+\b)--(\a+9,-4+\b);
\draw[very thick](\a+11,-2+\b)--(\a+11,-2.5+\b);
\draw[very thick](\a+11,-3.5+\b)--(\a+11,-4+\b);
\draw[very thick](\a+10,-3+\b)--(\a+10.5,-3+\b);
\draw[very thick](\a+11.5,-3+\b)--(\a+12,-3+\b);
\draw[very thick](\a+10.5,-3+\b)arc(270:360:0.5);
\draw[very thick](\a+11.5,-3+\b)arc(90:180:0.5);
\draw[very thick](\a+2,-5+\b)--(\a+4,-5+\b);
\draw[very thick](\a+3,-4+\b)--(\a+3,-6+\b);
\draw[very thick](\a+4,-5+\b)--(\a+6,-5+\b);
\draw[very thick](\a+5,-4+\b)--(\a+5,-6+\b);
\draw[very thick](\a+6,-5+\b)--(\a+8,-5+\b);
\draw[very thick](\a+7,-4+\b)--(\a+7,-6+\b);
\draw[very thick](\a+8,-5+\b)--(\a+10,-5+\b);
\draw[very thick](\a+9,-4+\b)--(\a+9,-6+\b);
\draw[very thick](\a+11,-4+\b)--(\a+11,-4.5+\b);
\draw[very thick](\a+11,-5.5+\b)--(\a+11,-6+\b);
\draw[very thick](\a+10,-5+\b)--(\a+10.5,-5+\b);
\draw[very thick](\a+11.5,-5+\b)--(\a+12,-5+\b);
\draw[very thick](\a+10.5,-5+\b)arc(270:360:0.5);
\draw[very thick](\a+11.5,-5+\b)arc(90:180:0.5);

\def\a{-16};
\def\b{0};
\draw[step=2, color=green](\a+2,-6+\b)grid(\a+12,-2+\b);
\draw[very thick](\a+3,-2+\b)--(\a+3,-2.5+\b);
\draw[very thick](\a+3,-3.5+\b)--(\a+3,-4+\b);
\draw[very thick](\a+2,-3+\b)--(\a+2.5,-3+\b);
\draw[very thick](\a+3.5,-3+\b)--(\a+4,-3+\b);
\draw[very thick](\a+2.5,-3+\b)arc(270:360:0.5);
\draw[very thick](\a+3.5,-3+\b)arc(90:180:0.5);
\draw[very thick](\a+4,-3+\b)--(\a+6,-3+\b);
\draw[very thick](\a+5,-2+\b)--(\a+5,-4+\b);
\draw[very thick](\a+6,-3+\b)--(\a+8,-3+\b);
\draw[very thick](\a+7,-2+\b)--(\a+7,-4+\b);
\draw[very thick](\a+8,-3+\b)--(\a+10,-3+\b);
\draw[very thick](\a+9,-2+\b)--(\a+9,-4+\b);
\draw[very thick](\a+10,-3+\b)--(\a+12,-3+\b);
\draw[very thick](\a+11,-2+\b)--(\a+11,-4+\b);
\draw[very thick](\a+3,-4+\b)--(\a+3,-4.5+\b);
\draw[very thick](\a+3,-5.5+\b)--(\a+3,-6+\b);
\draw[very thick](\a+2,-5+\b)--(\a+2.5,-5+\b);
\draw[very thick](\a+3.5,-5+\b)--(\a+4,-5+\b);
\draw[very thick](\a+2.5,-5+\b)arc(270:360:0.5);
\draw[very thick](\a+3.5,-5+\b)arc(90:180:0.5);
\draw[very thick](\a+4,-5+\b)--(\a+6,-5+\b);
\draw[very thick](\a+5,-4+\b)--(\a+5,-6+\b);
\draw[very thick](\a+6,-5+\b)--(\a+8,-5+\b);
\draw[very thick](\a+7,-4+\b)--(\a+7,-6+\b);
\draw[very thick](\a+8,-5+\b)--(\a+10,-5+\b);
\draw[very thick](\a+9,-4+\b)--(\a+9,-6+\b);
\draw[very thick](\a+11,-4+\b)--(\a+11,-4.5+\b);
\draw[very thick](\a+11,-5.5+\b)--(\a+11,-6+\b);
\draw[very thick](\a+10,-5+\b)--(\a+10.5,-5+\b);
\draw[very thick](\a+11.5,-5+\b)--(\a+12,-5+\b);
\draw[very thick](\a+10.5,-5+\b)arc(270:360:0.5);
\draw[very thick](\a+11.5,-5+\b)arc(90:180:0.5);
\draw[very thick,->] (-3,-4)--(1,-4);
\end{tikzpicture}
\caption{Chute moves on pipe dreams move one cross down and to the left, 
preserving the permutation.}
\label{fig:chute-move-PD}
\end{figure}

\begin{figure}[h]
\centering
\begin{tikzpicture}[scale=0.4]
\node at (0,0) {$\cdot$};
\node at (1,0) {$\cdot$};
\node at (0,-1) {$+$};
\node at (1,-1) {$+$};
\node at (0,-2) {$+$};
\node at (1,-2) {$+$};
\node at (0,-3) {$+$};
\node at (1,-3) {$\cdot$};
\def\a{5};
\node at (\a+0,0) {$\cdot$};
\node at (\a+1,0) {$+$};
\node at (\a+0,-1) {$+$};
\node at (\a+1,-1) {$+$};
\node at (\a+0,-2) {$+$};
\node at (\a+1,-2) {$+$};
\node at (\a+0,-3) {$\cdot$};
\node at (\a+1,-3) {$\cdot$};
\draw[very thick,->] (2,-1.5)--(4,-1.5);
\end{tikzpicture}
\qquad\qquad
\begin{tikzpicture}[scale=0.2]
\def\a{0};
\def\b{0};
\draw[step=2, color=green](\a+2,-10+\b)grid(\a+6,-2+\b);
\draw[very thick](\a+3,-2+\b)--(\a+3,-2.5+\b);
\draw[very thick](\a+3,-3.5+\b)--(\a+3,-4+\b);
\draw[very thick](\a+2,-3+\b)--(\a+2.5,-3+\b);
\draw[very thick](\a+3.5,-3+\b)--(\a+4,-3+\b);
\draw[very thick](\a+2.5,-3+\b)arc(270:360:0.5);
\draw[very thick](\a+3.5,-3+\b)arc(90:180:0.5);
\draw[very thick](\a+5,-2+\b)--(\a+5,-2.5+\b);
\draw[very thick](\a+5,-3.5+\b)--(\a+5,-4+\b);
\draw[very thick](\a+4,-3+\b)--(\a+4.5,-3+\b);
\draw[very thick](\a+5.5,-3+\b)--(\a+6,-3+\b);
\draw[very thick](\a+4.5,-3+\b)arc(270:360:0.5);
\draw[very thick](\a+5.5,-3+\b)arc(90:180:0.5);
\draw[very thick](\a+2,-5+\b)--(\a+4,-5+\b);
\draw[very thick](\a+3,-4+\b)--(\a+3,-6+\b);
\draw[very thick](\a+4,-5+\b)--(\a+6,-5+\b);
\draw[very thick](\a+5,-4+\b)--(\a+5,-6+\b);
\draw[very thick](\a+2,-7+\b)--(\a+4,-7+\b);
\draw[very thick](\a+3,-6+\b)--(\a+3,-8+\b);
\draw[very thick](\a+4,-7+\b)--(\a+6,-7+\b);
\draw[very thick](\a+5,-6+\b)--(\a+5,-8+\b);
\draw[very thick](\a+2,-9+\b)--(\a+4,-9+\b);
\draw[very thick](\a+3,-8+\b)--(\a+3,-10+\b);
\draw[very thick](\a+5,-8+\b)--(\a+5,-8.5+\b);
\draw[very thick](\a+5,-9.5+\b)--(\a+5,-10+\b);
\draw[very thick](\a+4,-9+\b)--(\a+4.5,-9+\b);
\draw[very thick](\a+5.5,-9+\b)--(\a+6,-9+\b);
\draw[very thick](\a+4.5,-9+\b)arc(270:360:0.5);
\draw[very thick](\a+5.5,-9+\b)arc(90:180:0.5);

\def\a{10};
\def\b{0};
\draw[step=2, color=green](\a+2,-10+\b)grid(\a+6,-2+\b);
\draw[very thick](\a+3,-2+\b)--(\a+3,-2.5+\b);
\draw[very thick](\a+3,-3.5+\b)--(\a+3,-4+\b);
\draw[very thick](\a+2,-3+\b)--(\a+2.5,-3+\b);
\draw[very thick](\a+3.5,-3+\b)--(\a+4,-3+\b);
\draw[very thick](\a+2.5,-3+\b)arc(270:360:0.5);
\draw[very thick](\a+3.5,-3+\b)arc(90:180:0.5);
\draw[very thick](\a+4,-3+\b)--(\a+6,-3+\b);
\draw[very thick](\a+5,-2+\b)--(\a+5,-4+\b);
\draw[very thick](\a+2,-5+\b)--(\a+4,-5+\b);
\draw[very thick](\a+3,-4+\b)--(\a+3,-6+\b);
\draw[very thick](\a+4,-5+\b)--(\a+6,-5+\b);
\draw[very thick](\a+5,-4+\b)--(\a+5,-6+\b);
\draw[very thick](\a+2,-7+\b)--(\a+4,-7+\b);
\draw[very thick](\a+3,-6+\b)--(\a+3,-8+\b);
\draw[very thick](\a+4,-7+\b)--(\a+6,-7+\b);
\draw[very thick](\a+5,-6+\b)--(\a+5,-8+\b);
\draw[very thick](\a+3,-8+\b)--(\a+3,-8.5+\b);
\draw[very thick](\a+3,-9.5+\b)--(\a+3,-10+\b);
\draw[very thick](\a+2,-9+\b)--(\a+2.5,-9+\b);
\draw[very thick](\a+3.5,-9+\b)--(\a+4,-9+\b);
\draw[very thick](\a+2.5,-9+\b)arc(270:360:0.5);
\draw[very thick](\a+3.5,-9+\b)arc(90:180:0.5);
\draw[very thick](\a+5,-8+\b)--(\a+5,-8.5+\b);
\draw[very thick](\a+5,-9.5+\b)--(\a+5,-10+\b);
\draw[very thick](\a+4,-9+\b)--(\a+4.5,-9+\b);
\draw[very thick](\a+5.5,-9+\b)--(\a+6,-9+\b);
\draw[very thick](\a+4.5,-9+\b)arc(270:360:0.5);
\draw[very thick](\a+5.5,-9+\b)arc(90:180:0.5);

\draw[very thick,->] (7,-6)--(11,-6);
\end{tikzpicture}
\caption{Ladder moves on pipe dreams move one cross up and to the right, preserving 
the permutation.}
\label{fig:ladder-move-PD}
\end{figure}

There are two types of simple local moves that preserve the set of
reduced pipe dreams for $w$: \textit{chute moves}
(\Cref{fig:chute-move-PD}) and \textit{ladder moves}
(\Cref{fig:ladder-move-PD}). A chute move can be thought of as
swapping a \+-tile with the first \elbow-tile on the row below and to
the left in such a way as to preserve the permutation. Similarly, a
ladder move swaps a \+-tile with the first \elbow-tile above and in the
column to its right. The \textit{generalized chute} and
\textit{generalized ladder} moves are defined by stretching out the
corresponding rectangles and requiring all entries other than the
corners to be crosses. The generalized chutes and ladders are formally defined below.

\begin{definition}\cite[Def. 2.6]{rubey2012maximal}\label{def:ladder}
For $w \in S_{n}$ and $D\in\textnormal{RPD}(w)$, assume there exists
$h<i$ and $j<k$ such that
\begin{itemize}
\item $D(i,j)=+$,
\item $D(h,j)=D(h,k)=D(i,k)=\sbullet$, and
\item every other tile in the rectangle with corners
$(i,j),(h,j),(h,k),(i,k)$ in $D$ is a $+$. 
\end{itemize}
Then the \textit{generalized ladder move} $\ladder_{ij}$ applied to $D$ is given by swapping the tiles at $(i,j)$ and $(h,k)$. The result is a new pipe dream $\ladder_{ij}(D)\in\rpd(w)$, see \Cref{fig:gen_ladder_example}. In this case, we say the cross tile at $(i,j)$ is \textit{ladder movable} in $D$. Furthermore, the inverse operation to a generalized ladder move $\ladder_{ij}^{-1}$ is a \textit{generalized chute move}. 
\end{definition}

Observe that if the cross tile at $(i,j)$ is ladder movable in
$D$, then there is a unique bump tile northeast of $(i,j)$ that
can be swapped with $(i,j)$ by a generalized ladder move. Therefore, if $(i,j)$ is ladder movable in $D$, then there is a unique $D'$ such that
$D'=\ladder_{ij}(D)$.

\begin{figure}[h]
    \centering
    \begin{tikzpicture}[scale=0.4]
\node at (0,0) {$\cdot$};
\node at (1,0) {$+$};
\node at (2,0) {$+$};
\node at (3,0) {$+$};
\node at (4,0) {$+$};
\node at (0,-1) {$+$};
\node at (1,-1) {$+$};
\node at (2,-1) {$+$};
\node at (3,-1) {$+$};
\node at (4,-1) {$+$};
\node at (0,-2) {$\cdot$};
\node at (1,-2) {$+$};
\node at (2,-2) {$+$};
\node at (3,-2) {$+$};
\node at (4,-2) {$\cdot$};
\def\a{-8};
\node at (\a+0,0) {$\cdot$};
\node at (\a+1,0) {$+$};
\node at (\a+2,0) {$+$};
\node at (\a+3,0) {$+$};
\node at (\a+4,0) {$\cdot$};
\node at (\a+0,-1) {$+$};
\node at (\a+1,-1) {$+$};
\node at (\a+2,-1) {$+$};
\node at (\a+3,-1) {$+$};
\node at (\a+4,-1) {$+$};
\node at (\a+0,-2) {$+$};
\node at (\a+1,-2) {$+$};
\node at (\a+2,-2) {$+$};
\node at (\a+3,-2) {$+$};
\node at (\a+4,-2) {$\cdot$};

\node at (-2.25,-0.25) {$\ladder_{31}$};
\draw[very thick,->] (-3,-1)--(-1,-1);
\end{tikzpicture}
\qquad
\begin{tikzpicture}[scale=0.2]
\def\a{-16};
\def\b{0};
\draw[step=2, color=green](\a+2,-8+\b)grid(\a+12,-2+\b);
\draw[very thick](\a+3,-2+\b)--(\a+3,-2.5+\b);
\draw[very thick](\a+3,-3.5+\b)--(\a+3,-4+\b);
\draw[very thick](\a+2,-3+\b)--(\a+2.5,-3+\b);
\draw[very thick](\a+3.5,-3+\b)--(\a+4,-3+\b);
\draw[very thick](\a+2.5,-3+\b)arc(270:360:0.5);
\draw[very thick](\a+3.5,-3+\b)arc(90:180:0.5);
\draw[very thick](\a+4,-3+\b)--(\a+6,-3+\b);
\draw[very thick](\a+5,-2+\b)--(\a+5,-4+\b);
\draw[very thick](\a+6,-3+\b)--(\a+8,-3+\b);
\draw[very thick](\a+7,-2+\b)--(\a+7,-4+\b);
\draw[very thick](\a+8,-3+\b)--(\a+10,-3+\b);
\draw[very thick](\a+9,-2+\b)--(\a+9,-4+\b);
\draw[very thick](\a+11,-2+\b)--(\a+11,-2.5+\b);
\draw[very thick](\a+11,-3.5+\b)--(\a+11,-4+\b);
\draw[very thick](\a+10,-3+\b)--(\a+10.5,-3+\b);
\draw[very thick](\a+11.5,-3+\b)--(\a+12,-3+\b);
\draw[very thick](\a+10.5,-3+\b)arc(270:360:0.5);
\draw[very thick](\a+11.5,-3+\b)arc(90:180:0.5);

\draw[very thick](\a+2,-5+\b)--(\a+4,-5+\b);
\draw[very thick](\a+3,-4+\b)--(\a+3,-6+\b);
\draw[very thick](\a+4,-5+\b)--(\a+6,-5+\b);
\draw[very thick](\a+5,-4+\b)--(\a+5,-6+\b);
\draw[very thick](\a+6,-5+\b)--(\a+8,-5+\b);
\draw[very thick](\a+7,-4+\b)--(\a+7,-6+\b);
\draw[very thick](\a+8,-5+\b)--(\a+10,-5+\b);
\draw[very thick](\a+9,-4+\b)--(\a+9,-6+\b);
\draw[very thick](\a+10,-5+\b)--(\a+12,-5+\b);
\draw[very thick](\a+11,-4+\b)--(\a+11,-6+\b);

\draw[very thick](\a+2,-7+\b)--(\a+4,-7+\b);
\draw[very thick](\a+3,-6+\b)--(\a+3,-8+\b);
\draw[very thick](\a+4,-7+\b)--(\a+6,-7+\b);
\draw[very thick](\a+5,-6+\b)--(\a+5,-8+\b);
\draw[very thick](\a+6,-7+\b)--(\a+8,-7+\b);
\draw[very thick](\a+7,-6+\b)--(\a+7,-8+\b);
\draw[very thick](\a+8,-7+\b)--(\a+10,-7+\b);
\draw[very thick](\a+9,-6+\b)--(\a+9,-8+\b);
\draw[very thick](\a+11,-6+\b)--(\a+11,-6.5+\b);
\draw[very thick](\a+11,-7.5+\b)--(\a+11,-8+\b);
\draw[very thick](\a+10,-7+\b)--(\a+10.5,-7+\b);
\draw[very thick](\a+11.5,-7+\b)--(\a+12,-7+\b);
\draw[very thick](\a+10.5,-7+\b)arc(270:360:0.5);
\draw[very thick](\a+11.5,-7+\b)arc(90:180:0.5);

\node at (-1.5,-3.5) {$\ladder_{31}$};
\draw[very thick,->] (-3,-5)--(1,-5);

\def\a{0};
\def\b{0};
\draw[step=2, color=green](\a+2,-8+\b)grid(\a+12,-2+\b);
\draw[very thick](\a+3,-4+\b)--(\a+3,-4.5+\b);
\draw[very thick](\a+3,-2+\b)--(\a+3,-2.5+\b);
\draw[very thick](\a+3,-3.5+\b)--(\a+3,-4+\b);
\draw[very thick](\a+2,-3+\b)--(\a+2.5,-3+\b);
\draw[very thick](\a+3.5,-3+\b)--(\a+4,-3+\b);
\draw[very thick](\a+2.5,-3+\b)arc(270:360:0.5);
\draw[very thick](\a+3.5,-3+\b)arc(90:180:0.5);
\draw[very thick](\a+4,-3+\b)--(\a+6,-3+\b);
\draw[very thick](\a+5,-2+\b)--(\a+5,-4+\b);
\draw[very thick](\a+6,-3+\b)--(\a+8,-3+\b);
\draw[very thick](\a+7,-2+\b)--(\a+7,-4+\b);
\draw[very thick](\a+8,-3+\b)--(\a+10,-3+\b);
\draw[very thick](\a+9,-2+\b)--(\a+9,-4+\b);
\draw[very thick](\a+10,-3+\b)--(\a+12,-3+\b);
\draw[very thick](\a+11,-2+\b)--(\a+11,-4+\b);

\draw[very thick](\a+2,-5+\b)--(\a+4,-5+\b);
\draw[very thick](\a+3,-4+\b)--(\a+3,-6+\b);
\draw[very thick](\a+4,-5+\b)--(\a+6,-5+\b);
\draw[very thick](\a+5,-4+\b)--(\a+5,-6+\b);
\draw[very thick](\a+6,-5+\b)--(\a+8,-5+\b);
\draw[very thick](\a+7,-4+\b)--(\a+7,-6+\b);
\draw[very thick](\a+8,-5+\b)--(\a+10,-5+\b);
\draw[very thick](\a+9,-4+\b)--(\a+9,-6+\b);
\draw[very thick](\a+10,-5+\b)--(\a+12,-5+\b);
\draw[very thick](\a+11,-4+\b)--(\a+11,-6+\b);

\draw[very thick](\a+3,-6+\b)--(\a+3,-6.5+\b);
\draw[very thick](\a+3,-7.5+\b)--(\a+3,-8+\b);
\draw[very thick](\a+2,-7+\b)--(\a+2.5,-7+\b);
\draw[very thick](\a+3.5,-7+\b)--(\a+4,-7+\b);
\draw[very thick](\a+2.5,-7+\b)arc(270:360:0.5);
\draw[very thick](\a+3.5,-7+\b)arc(90:180:0.5);
\draw[very thick](\a+4,-7+\b)--(\a+6,-7+\b);
\draw[very thick](\a+5,-6+\b)--(\a+5,-8+\b);
\draw[very thick](\a+6,-7+\b)--(\a+8,-7+\b);
\draw[very thick](\a+7,-6+\b)--(\a+7,-8+\b);
\draw[very thick](\a+8,-7+\b)--(\a+10,-7+\b);
\draw[very thick](\a+9,-6+\b)--(\a+9,-8+\b);
\draw[very thick](\a+11,-6+\b)--(\a+11,-6.5+\b);
\draw[very thick](\a+11,-7.5+\b)--(\a+11,-8+\b);
\draw[very thick](\a+10,-7+\b)--(\a+10.5,-7+\b);
\draw[very thick](\a+11.5,-7+\b)--(\a+12,-7+\b);
\draw[very thick](\a+10.5,-7+\b)arc(270:360:0.5);
\draw[very thick](\a+11.5,-7+\b)arc(90:180:0.5);
\end{tikzpicture}
    \caption{A pipe dream in which the cross tile at $(3,1)$ is ladder movable by swapping it with the bump tile $(1,5)$ is shown in the simplified form and with the wire tiles drawn. }
    \label{fig:gen_ladder_example}
\end{figure}

\begin{definition}\cite[Conj. 2.8]{rubey2012maximal}\label{def:rubey.lattice}
Let $w \in S_{n}$. The \textit{Rubey lattice} $(\rpd(w),<)$ is the
partial order on $\rpd(w)$ with covering relations given by $D\lessdot
\ladder_{ij}(D)$ whenever $(i,j)$ is ladder movable in $D$. We abuse
notation and also use $\rpd(w)$ to refer to the Rubey lattice.
\end{definition}

\begin{figure}[h]
\centering
\begin{tikzpicture}[scale=0.25]
\def\a{0};
\def\b{0};
\draw[color=green](\a+2,-2+\b)--(\a+10,-2+\b);
\draw[color=green](\a+2,-2+\b)--(\a+2,-10+\b);
\draw[color=green](\a+2,-4+\b)--(\a+10,-4+\b);
\draw[color=green](\a+4,-2+\b)--(\a+4,-10+\b);
\draw[color=green](\a+2,-6+\b)--(\a+8,-6+\b);
\draw[color=green](\a+6,-2+\b)--(\a+6,-8+\b);
\draw[color=green](\a+2,-8+\b)--(\a+6,-8+\b);
\draw[color=green](\a+8,-2+\b)--(\a+8,-6+\b);
\draw[color=green](\a+2,-10+\b)--(\a+4,-10+\b);
\draw[color=green](\a+10,-2+\b)--(\a+10,-4+\b);
\draw[very thick](\a+3,-2+\b)--(\a+3,-2.50000000000000+\b);
\draw[very thick](\a+3,-3.50000000000000+\b)--(\a+3,-4+\b);
\draw[very thick](\a+2,-3+\b)--(\a+2.50000000000000,-3+\b);
\draw[very thick](\a+3.50000000000000,-3+\b)--(\a+4,-3+\b);
\draw[very thick](\a+2.50000000000000,-3+\b)arc(270:360:0.500000000000000);
\draw[very thick](\a+3.50000000000000,-3+\b)arc(90:180:0.500000000000000);
\draw[very thick](\a+5,-2+\b)--(\a+5,-2.50000000000000+\b);
\draw[very thick](\a+5,-3.50000000000000+\b)--(\a+5,-4+\b);
\draw[very thick](\a+4,-3+\b)--(\a+4.50000000000000,-3+\b);
\draw[very thick](\a+5.50000000000000,-3+\b)--(\a+6,-3+\b);
\draw[very thick](\a+4.50000000000000,-3+\b)arc(270:360:0.500000000000000);
\draw[very thick](\a+5.50000000000000,-3+\b)arc(90:180:0.500000000000000);
\draw[very thick](\a+7,-2+\b)--(\a+7,-2.50000000000000+\b);
\draw[very thick](\a+7,-3.50000000000000+\b)--(\a+7,-4+\b);
\draw[very thick](\a+6,-3+\b)--(\a+6.50000000000000,-3+\b);
\draw[very thick](\a+7.50000000000000,-3+\b)--(\a+8,-3+\b);
\draw[very thick](\a+6.50000000000000,-3+\b)arc(270:360:0.500000000000000);
\draw[very thick](\a+7.50000000000000,-3+\b)arc(90:180:0.500000000000000);
\draw[very thick](\a+2,-5+\b)--(\a+4,-5+\b);
\draw[very thick](\a+3,-4+\b)--(\a+3,-6+\b);
\draw[very thick](\a+4,-5+\b)--(\a+6,-5+\b);
\draw[very thick](\a+5,-4+\b)--(\a+5,-6+\b);
\draw[very thick](\a+2,-7+\b)--(\a+4,-7+\b);
\draw[very thick](\a+3,-6+\b)--(\a+3,-8+\b);
\draw[very thick](\a+8,-3+\b)arc(270:360:1);
\draw[very thick](\a+6,-5+\b)arc(270:360:1);
\draw[very thick](\a+4,-7+\b)arc(270:360:1);
\draw[very thick](\a+2,-9+\b)arc(270:360:1);

\def\a{12};
\def\b{5};
\draw[color=green](\a+2,-2+\b)--(\a+10,-2+\b);
\draw[color=green](\a+2,-2+\b)--(\a+2,-10+\b);
\draw[color=green](\a+2,-4+\b)--(\a+10,-4+\b);
\draw[color=green](\a+4,-2+\b)--(\a+4,-10+\b);
\draw[color=green](\a+2,-6+\b)--(\a+8,-6+\b);
\draw[color=green](\a+6,-2+\b)--(\a+6,-8+\b);
\draw[color=green](\a+2,-8+\b)--(\a+6,-8+\b);
\draw[color=green](\a+8,-2+\b)--(\a+8,-6+\b);
\draw[color=green](\a+2,-10+\b)--(\a+4,-10+\b);
\draw[color=green](\a+10,-2+\b)--(\a+10,-4+\b);
\draw[very thick](\a+3,-2+\b)--(\a+3,-2.50000000000000+\b);
\draw[very thick](\a+3,-3.50000000000000+\b)--(\a+3,-4+\b);
\draw[very thick](\a+2,-3+\b)--(\a+2.50000000000000,-3+\b);
\draw[very thick](\a+3.50000000000000,-3+\b)--(\a+4,-3+\b);
\draw[very thick](\a+2.50000000000000,-3+\b)arc(270:360:0.500000000000000);
\draw[very thick](\a+3.50000000000000,-3+\b)arc(90:180:0.500000000000000);
\draw[very thick](\a+4,-3+\b)--(\a+6,-3+\b);
\draw[very thick](\a+5,-2+\b)--(\a+5,-4+\b);
\draw[very thick](\a+7,-2+\b)--(\a+7,-2.50000000000000+\b);
\draw[very thick](\a+7,-3.50000000000000+\b)--(\a+7,-4+\b);
\draw[very thick](\a+6,-3+\b)--(\a+6.50000000000000,-3+\b);
\draw[very thick](\a+7.50000000000000,-3+\b)--(\a+8,-3+\b);
\draw[very thick](\a+6.50000000000000,-3+\b)arc(270:360:0.500000000000000);
\draw[very thick](\a+7.50000000000000,-3+\b)arc(90:180:0.500000000000000);
\draw[very thick](\a+2,-5+\b)--(\a+4,-5+\b);
\draw[very thick](\a+3,-4+\b)--(\a+3,-6+\b);
\draw[very thick](\a+4,-5+\b)--(\a+6,-5+\b);
\draw[very thick](\a+5,-4+\b)--(\a+5,-6+\b);
\draw[very thick](\a+3,-6+\b)--(\a+3,-6.50000000000000+\b);
\draw[very thick](\a+3,-7.50000000000000+\b)--(\a+3,-8+\b);
\draw[very thick](\a+2,-7+\b)--(\a+2.50000000000000,-7+\b);
\draw[very thick](\a+3.50000000000000,-7+\b)--(\a+4,-7+\b);
\draw[very thick](\a+2.50000000000000,-7+\b)arc(270:360:0.500000000000000);
\draw[very thick](\a+3.50000000000000,-7+\b)arc(90:180:0.500000000000000);
\draw[very thick](\a+8,-3+\b)arc(270:360:1);
\draw[very thick](\a+6,-5+\b)arc(270:360:1);
\draw[very thick](\a+4,-7+\b)arc(270:360:1);
\draw[very thick](\a+2,-9+\b)arc(270:360:1);

\def\a{12};
\def\b{-5};
\draw[color=green](\a+2,-2+\b)--(\a+10,-2+\b);
\draw[color=green](\a+2,-2+\b)--(\a+2,-10+\b);
\draw[color=green](\a+2,-4+\b)--(\a+10,-4+\b);
\draw[color=green](\a+4,-2+\b)--(\a+4,-10+\b);
\draw[color=green](\a+2,-6+\b)--(\a+8,-6+\b);
\draw[color=green](\a+6,-2+\b)--(\a+6,-8+\b);
\draw[color=green](\a+2,-8+\b)--(\a+6,-8+\b);
\draw[color=green](\a+8,-2+\b)--(\a+8,-6+\b);
\draw[color=green](\a+2,-10+\b)--(\a+4,-10+\b);
\draw[color=green](\a+10,-2+\b)--(\a+10,-4+\b);
\draw[very thick](\a+3,-2+\b)--(\a+3,-2.50000000000000+\b);
\draw[very thick](\a+3,-3.50000000000000+\b)--(\a+3,-4+\b);
\draw[very thick](\a+2,-3+\b)--(\a+2.50000000000000,-3+\b);
\draw[very thick](\a+3.50000000000000,-3+\b)--(\a+4,-3+\b);
\draw[very thick](\a+2.50000000000000,-3+\b)arc(270:360:0.500000000000000);
\draw[very thick](\a+3.50000000000000,-3+\b)arc(90:180:0.500000000000000);
\draw[very thick](\a+5,-2+\b)--(\a+5,-2.50000000000000+\b);
\draw[very thick](\a+5,-3.50000000000000+\b)--(\a+5,-4+\b);
\draw[very thick](\a+4,-3+\b)--(\a+4.50000000000000,-3+\b);
\draw[very thick](\a+5.50000000000000,-3+\b)--(\a+6,-3+\b);
\draw[very thick](\a+4.50000000000000,-3+\b)arc(270:360:0.500000000000000);
\draw[very thick](\a+5.50000000000000,-3+\b)arc(90:180:0.500000000000000);
\draw[very thick](\a+6,-3+\b)--(\a+8,-3+\b);
\draw[very thick](\a+7,-2+\b)--(\a+7,-4+\b);
\draw[very thick](\a+2,-5+\b)--(\a+4,-5+\b);
\draw[very thick](\a+3,-4+\b)--(\a+3,-6+\b);
\draw[very thick](\a+5,-4+\b)--(\a+5,-4.50000000000000+\b);
\draw[very thick](\a+5,-5.50000000000000+\b)--(\a+5,-6+\b);
\draw[very thick](\a+4,-5+\b)--(\a+4.50000000000000,-5+\b);
\draw[very thick](\a+5.50000000000000,-5+\b)--(\a+6,-5+\b);
\draw[very thick](\a+4.50000000000000,-5+\b)arc(270:360:0.500000000000000);
\draw[very thick](\a+5.50000000000000,-5+\b)arc(90:180:0.500000000000000);
\draw[very thick](\a+2,-7+\b)--(\a+4,-7+\b);
\draw[very thick](\a+3,-6+\b)--(\a+3,-8+\b);
\draw[very thick](\a+8,-3+\b)arc(270:360:1);
\draw[very thick](\a+6,-5+\b)arc(270:360:1);
\draw[very thick](\a+4,-7+\b)arc(270:360:1);
\draw[very thick](\a+2,-9+\b)arc(270:360:1);

\def\a{24};
\def\b{-5};
\draw[color=green](\a+2,-2+\b)--(\a+10,-2+\b);
\draw[color=green](\a+2,-2+\b)--(\a+2,-10+\b);
\draw[color=green](\a+2,-4+\b)--(\a+10,-4+\b);
\draw[color=green](\a+4,-2+\b)--(\a+4,-10+\b);
\draw[color=green](\a+2,-6+\b)--(\a+8,-6+\b);
\draw[color=green](\a+6,-2+\b)--(\a+6,-8+\b);
\draw[color=green](\a+2,-8+\b)--(\a+6,-8+\b);
\draw[color=green](\a+8,-2+\b)--(\a+8,-6+\b);
\draw[color=green](\a+2,-10+\b)--(\a+4,-10+\b);
\draw[color=green](\a+10,-2+\b)--(\a+10,-4+\b);
\draw[very thick](\a+3,-2+\b)--(\a+3,-2.50000000000000+\b);
\draw[very thick](\a+3,-3.50000000000000+\b)--(\a+3,-4+\b);
\draw[very thick](\a+2,-3+\b)--(\a+2.50000000000000,-3+\b);
\draw[very thick](\a+3.50000000000000,-3+\b)--(\a+4,-3+\b);
\draw[very thick](\a+2.50000000000000,-3+\b)arc(270:360:0.500000000000000);
\draw[very thick](\a+3.50000000000000,-3+\b)arc(90:180:0.500000000000000);
\draw[very thick](\a+4,-3+\b)--(\a+6,-3+\b);
\draw[very thick](\a+5,-2+\b)--(\a+5,-4+\b);
\draw[very thick](\a+6,-3+\b)--(\a+8,-3+\b);
\draw[very thick](\a+7,-2+\b)--(\a+7,-4+\b);
\draw[very thick](\a+3,-4+\b)--(\a+3,-4.50000000000000+\b);
\draw[very thick](\a+3,-5.50000000000000+\b)--(\a+3,-6+\b);
\draw[very thick](\a+2,-5+\b)--(\a+2.50000000000000,-5+\b);
\draw[very thick](\a+3.50000000000000,-5+\b)--(\a+4,-5+\b);
\draw[very thick](\a+2.50000000000000,-5+\b)arc(270:360:0.500000000000000);
\draw[very thick](\a+3.50000000000000,-5+\b)arc(90:180:0.500000000000000);
\draw[very thick](\a+5,-4+\b)--(\a+5,-4.50000000000000+\b);
\draw[very thick](\a+5,-5.50000000000000+\b)--(\a+5,-6+\b);
\draw[very thick](\a+4,-5+\b)--(\a+4.50000000000000,-5+\b);
\draw[very thick](\a+5.50000000000000,-5+\b)--(\a+6,-5+\b);
\draw[very thick](\a+4.50000000000000,-5+\b)arc(270:360:0.500000000000000);
\draw[very thick](\a+5.50000000000000,-5+\b)arc(90:180:0.500000000000000);
\draw[very thick](\a+2,-7+\b)--(\a+4,-7+\b);
\draw[very thick](\a+3,-6+\b)--(\a+3,-8+\b);
\draw[very thick](\a+8,-3+\b)arc(270:360:1);
\draw[very thick](\a+6,-5+\b)arc(270:360:1);
\draw[very thick](\a+4,-7+\b)arc(270:360:1);
\draw[very thick](\a+2,-9+\b)arc(270:360:1);

\def\a{36};
\def\b{-5};
\draw[color=green](\a+2,-2+\b)--(\a+10,-2+\b);
\draw[color=green](\a+2,-2+\b)--(\a+2,-10+\b);
\draw[color=green](\a+2,-4+\b)--(\a+10,-4+\b);
\draw[color=green](\a+4,-2+\b)--(\a+4,-10+\b);
\draw[color=green](\a+2,-6+\b)--(\a+8,-6+\b);
\draw[color=green](\a+6,-2+\b)--(\a+6,-8+\b);
\draw[color=green](\a+2,-8+\b)--(\a+6,-8+\b);
\draw[color=green](\a+8,-2+\b)--(\a+8,-6+\b);
\draw[color=green](\a+2,-10+\b)--(\a+4,-10+\b);
\draw[color=green](\a+10,-2+\b)--(\a+10,-4+\b);
\draw[very thick](\a+3,-2+\b)--(\a+3,-2.50000000000000+\b);
\draw[very thick](\a+3,-3.50000000000000+\b)--(\a+3,-4+\b);
\draw[very thick](\a+2,-3+\b)--(\a+2.50000000000000,-3+\b);
\draw[very thick](\a+3.50000000000000,-3+\b)--(\a+4,-3+\b);
\draw[very thick](\a+2.50000000000000,-3+\b)arc(270:360:0.500000000000000);
\draw[very thick](\a+3.50000000000000,-3+\b)arc(90:180:0.500000000000000);
\draw[very thick](\a+4,-3+\b)--(\a+6,-3+\b);
\draw[very thick](\a+5,-2+\b)--(\a+5,-4+\b);
\draw[very thick](\a+6,-3+\b)--(\a+8,-3+\b);
\draw[very thick](\a+7,-2+\b)--(\a+7,-4+\b);
\draw[very thick](\a+3,-4+\b)--(\a+3,-4.50000000000000+\b);
\draw[very thick](\a+3,-5.50000000000000+\b)--(\a+3,-6+\b);
\draw[very thick](\a+2,-5+\b)--(\a+2.50000000000000,-5+\b);
\draw[very thick](\a+3.50000000000000,-5+\b)--(\a+4,-5+\b);
\draw[very thick](\a+2.50000000000000,-5+\b)arc(270:360:0.500000000000000);
\draw[very thick](\a+3.50000000000000,-5+\b)arc(90:180:0.500000000000000);
\draw[very thick](\a+4,-5+\b)--(\a+6,-5+\b);
\draw[very thick](\a+5,-4+\b)--(\a+5,-6+\b);
\draw[very thick](\a+3,-6+\b)--(\a+3,-6.50000000000000+\b);
\draw[very thick](\a+3,-7.50000000000000+\b)--(\a+3,-8+\b);
\draw[very thick](\a+2,-7+\b)--(\a+2.50000000000000,-7+\b);
\draw[very thick](\a+3.50000000000000,-7+\b)--(\a+4,-7+\b);
\draw[very thick](\a+2.50000000000000,-7+\b)arc(270:360:0.500000000000000);
\draw[very thick](\a+3.50000000000000,-7+\b)arc(90:180:0.500000000000000);
\draw[very thick](\a+8,-3+\b)arc(270:360:1);
\draw[very thick](\a+6,-5+\b)arc(270:360:1);
\draw[very thick](\a+4,-7+\b)arc(270:360:1);
\draw[very thick](\a+2,-9+\b)arc(270:360:1);

\draw[very thick, ->](10,-5)--(12,-2);
\draw[very thick, ->](10,-7)--(12,-10);
\draw[very thick, ->](21,-11)--(24,-11);
\draw[very thick, ->](33,-11)--(36,-11);
\draw[very thick, ->](24,0)--(36,-6);   
\end{tikzpicture}
\caption{Construction of the Rubey lattice on $\rpd(1432)$ by
generalized ladder moves. $\Dbot(1432)$ is on the left and
$\Dtop(1432)$ is on the right.}
\label{fig:1432.ladders}
\end{figure}

It is known that the Rubey lattice is connected and has unique minimal
and maximal elements by the following theorem. Furthermore, $\rpd(w)$
and the dual of $\rpd(w^{-1})$ are isomorphic as posets. The
isomorphism is given by transposing the tilings of pipe dreams and interchanging generalized ladder moves $\ladder_{ij}$ in $\rpd(w)$ with generalized chute moves $\ladder^{-1}_{ji}$ in the dual of $\rpd(w^{-1})$. The Rubey lattice has many more edges than necessary to simply connect the reduced pipe dreams for $w$.

\begin{theorem}\cite[Thm 3.7]{billey-bergeron}\label{thm:chutes.and.ladders}
For any $w \in S_{n}$, every reduced pipe dream for $w$ can be
obtained by a sequence of ladder moves from $\Dbot(w)$, and every
reduced pipe dream for $w$ can be obtained by a sequence of chute
moves from $\Dtop(w)$. The minimal pipe dream $\Dbot(w)$ is obtained from the
diagram of the permutation $\textbf{D}(w)$ by replacing each cell by a \+
and left justifying the crosses in each row. The maximal pipe dream
$\Dtop(w)$ is obtained from the diagram of the permutation $\textbf{D}(w)$ by
replacing each cell by a \+ and upper justifying the crosses in
each column.
\end{theorem}

\begin{example}\label{ex:1432.1}
The reduced pipe dreams for $w=1432$, along with all possible generalized ladder moves, are shown in
\Cref{fig:1432.ladders}.
\end{example}

\section{The Move Operation}\label{sec:move}

Throughout the rest of this paper, assume $n$ is a fixed positive integer and $w \in S_{n}$ is a fixed
permutation. All pipe dreams discussed below are reduced pipe dreams for
$w$. In this section, we introduce the move operation $\move_{ij}$
and examine its properties. The pipe dream $\move_{ij}(D)$ will be
defined for $D\in\rpd(w)$ such that $D(i,j)=+$ and there exists a bump
tile in column $j$ of $D$ above $(i,j)$. The effect of
$\move_{ij}$ will be to turn the tile $(i,j)$ in $D$ into a
bump tile by a particular sequence of generalized ladder moves depending on $D$.

\subsection{Definitions and basic properties of move operators}

\begin{definition}\label{movable_defn}
     Let $D\in\textnormal{RPD}(w)$. If $D(i,j)=+$ and $\{h\in[1,i-1]\mid D(h,j)=\sbullet\}$ is nonempty, we say that $(i,j)$ is \textit{movable} in $D$, and we define
     \[
     \upperbump_{ij}(D)=\max\{h\in[1,i-1]\mid D(h,j)=\sbullet\}.\]
     Similarly, define
     \[\rightbump_{ij}(D)=\min\{k\in[j+1,n]\mid D(i,k)=\sbullet\}.\]
    Thus, $\upperbump_{ij}(D)$ is the row of the lowest bump tile above $(i,j)$ in $D$, and $\rightbump_{ij}(D)$ is the column of the leftmost bump tile to the right of $(i,j)$. For any cross tile $(i,j)$ in $D$, $\rightbump_{ij}(D)$ is well-defined since $D(i,n)=\sbullet$ for all reduced pipe dreams $D$ representing $w\in S_n$.
\end{definition}

\begin{definition}
    Let $D\in\textnormal{RPD}(w)$ such that $(i,j)$ is a movable cross tile in $D$, and let $h=h_{ij}(D)$ and $k=k_{ij}(D)$. Define $\rect_{ij}(D)=[h,i]\times[j,k]$. Let 
    \[
    \maxrow_{ij}(D)=\max\big\{p\in[h,i)\mid\textnormal{row $p$ of $D$ contains a $\sbullet$ tile in $\rect_{ij}(D)$}\big\},\]
    and let
    \[\mincol_{ij}(D)=\min\big\{q\in (j,k]\mid\textnormal{column $q$ of $D$ contains a $\sbullet$ tile in $\rect_{ij}(D)$}\big\}.\]
\end{definition}

\def\RPlus{\mathbin{\rule[0.3em]{.7em}{.16em}\hspace{-.455em}\rule[0em]{.16em}{.7em}}}


\begin{example}
    Consider the pipe dream $D\in\rpd(126543)$ shown in \Cref{fig:boxed.rectangle}, where $(3,2)$ is
    movable. 
    We have that $\upperbump_{32}(D)=1$ and $\rightbump_{32}(D)=4$, so $\rect_{32}(D)=[1,3]\times[2,4]$. Furthermore, $\maxrow_{32}(D)=2$, and $\mincol_{32}(D)=3$.
\end{example}
\begin{figure}[h]
\[
\begin{NiceMatrix}[columns-width=auto]
    \sbullet & \textcolor{cyan}{\sbullet} & \textcolor{cyan}{\sbullet} & \textcolor{cyan}{\sbullet} & \sbullet & \sbullet \\
    \sbullet & \textcolor{cyan}{+}  & \textcolor{cyan}{\sbullet} & \textcolor{cyan}{\sbullet} & \sbullet & \\
    +        & \textcolor{cyan}{+}  & \textcolor{cyan}{+}  & \textcolor{cyan}{+}  & & \\
    +        & +                         & \sbullet                  & & & \\
    \sbullet & \sbullet                  &                           & & & \\
    \sbullet &                           &                           & & &
    \CodeAfter
    \tikz \draw[green, thin]
        ([xshift=-7pt,yshift=7pt]1-2.north west)
        rectangle
        ([xshift=3pt,yshift=-3pt]3-4.south east);
    \end{NiceMatrix}
\]
    \caption{Pipe dream with $\rect_{32}$ boxed.}
    \label{fig:boxed.rectangle}
\end{figure}

\begin{lemma}\label{lmove_conditions}
    Let $D\in\rpd(w)$ such that $(i,j)$ is a movable cross tile in $D$. Then $(i,j)$ is ladder movable if and only if $\upperbump_{ij}(D)=\maxrow_{ij}(D)$ and $\rightbump_{ij}(D)=\mincol_{ij}(D)$.
\end{lemma}
\begin{proof}
    Let $h=\upperbump_{ij}(D)$ and $k=\rightbump_{ij}(D)$. We have that $(i,j)$ is ladder movable in $D$ if and only if the set of bump tiles in $\rect_{ij}(D)$ consists exactly of the corner tiles $(h,j)$, $(h,k)$, and $(i,k)$. Thus, if $(i,j)$ is ladder movable, then $h=\maxrow_{ij}(D)$ and $k=\mincol_{ij}(D)$. 
    
    On the other hand, if $h=\maxrow_{ij}(D)$ and $k=\mincol_{ij}(D)$, then the only tiles in $\rect_{ij}(D)$ which could be bump tiles are the corners, and $D(h,j)=D(i,k)=\sbullet$ by assumption. If $D(h,k)=+$, then $D$ is not reduced contrary to the hypothesis $D \in \rpd(w)$, see \Cref{fig:nonreduced_pipe_dream}. Therefore, $D(h,k)=\sbullet$, so $(i,j)$ is ladder movable in $D$.
\end{proof}
\begin{figure}[h]
    \centering
    \begin{tikzpicture}[scale=0.3]
\def\a{-16};
\def\b{0};
\draw[step=2, color=green,thin](\a+2,-6+\b)grid(\a+12,\b);

\draw[very thick](\a+3,\b)--(\a+3,-0.5+\b);
\draw[color=blue][very thick](\a+3,-1.5+\b)--(\a+3,-2+\b);
\draw[very thick](\a+2,-1+\b)--(\a+2.5,-1+\b);
\draw[color=blue][very thick](\a+3.5,-1+\b)--(\a+4,-1+\b);
\draw[very thick](\a+2.5,-1+\b)arc(270:360:0.5);
\draw[color=blue][very thick](\a+3.5,-1+\b)arc(90:180:0.5);
\draw[color=blue][very thick](\a+4,-1+\b)--(\a+6,-1+\b);
\draw[very thick](\a+5,\b)--(\a+5,-2+\b);
\draw[color=blue][very thick](\a+6,-1+\b)--(\a+8,-1+\b);
\draw[very thick](\a+7,\b)--(\a+7,-2+\b);
\draw[color=blue][very thick](\a+8,-1+\b)--(\a+10,-1+\b);
\draw[very thick](\a+9,\b)--(\a+9,-2+\b);
\draw[color=blue][very thick](\a+10,-1+\b)--(\a+12,-1+\b);
\draw[color=orange][very thick](\a+11,\b)--(\a+11,-2+\b);

\draw[very thick](\a+2,-3+\b)--(\a+4,-3+\b);
\draw[color=blue][very thick](\a+3,-2+\b)--(\a+3,-4+\b);
\draw[very thick](\a+4,-3+\b)--(\a+6,-3+\b);
\draw[very thick](\a+5,-2+\b)--(\a+5,-4+\b);
\draw[very thick](\a+6,-3+\b)--(\a+8,-3+\b);
\draw[very thick](\a+7,-2+\b)--(\a+7,-4+\b);
\draw[very thick](\a+8,-3+\b)--(\a+10,-3+\b);
\draw[very thick](\a+9,-2+\b)--(\a+9,-4+\b);
\draw[very thick](\a+10,-3+\b)--(\a+12,-3+\b);
\draw[color=orange][very thick](\a+11,-2+\b)--(\a+11,-4+\b);

\draw[color=orange][very thick](\a+2,-5+\b)--(\a+4,-5+\b);
\draw[color=blue][very thick](\a+3,-4+\b)--(\a+3,-6+\b);
\draw[color=orange][very thick](\a+4,-5+\b)--(\a+6,-5+\b);
\draw[very thick](\a+5,-4+\b)--(\a+5,-6+\b);
\draw[color=orange][very thick](\a+6,-5+\b)--(\a+8,-5+\b);
\draw[very thick](\a+7,-4+\b)--(\a+7,-6+\b);
\draw[color=orange][very thick](\a+8,-5+\b)--(\a+10,-5+\b);
\draw[very thick](\a+9,-4+\b)--(\a+9,-6+\b);
\draw[color=orange][very thick](\a+11,-4+\b)--(\a+11,-4.5+\b);
\draw[very thick](\a+11,-5.5+\b)--(\a+11,-6+\b);
\draw[color=orange][very thick](\a+10,-5+\b)--(\a+10.5,-5+\b);
\draw[very thick](\a+11.5,-5+\b)--(\a+12,-5+\b);
\draw[color=orange][very thick](\a+10.5,-5+\b)arc(270:360:0.5);
\draw[very thick](\a+11.5,-5+\b)arc(90:180:0.5);
\end{tikzpicture}
    \caption{The submatrix above cannot belong to any reduced pipe dream, since the blue and orange wires cross at both (3,1) and $(h_{31}(D),k_{31}(D))=(1,5)$.}
    \label{fig:nonreduced_pipe_dream}
\end{figure}

We are now ready to define the $\move_{ij}$ operation. We first give a recursive algorithm for $\move_{ij}(D)$ in terms of generalized ladder moves, and we then give an explicit characterization of $\move_{ij}(D)$.

\begin{definition}\label{M_recursion}
    Let $D\in\rpd(w)$ such that $(i,j)$ is a movable cross tile in
$D$. Let $a=\maxrow_{ij}(D)$, $b=\mincol_{ij}(D)$,
$h=\upperbump_{ij}(D)$, and $k=\rightbump_{ij}(D)$. Define the \textit{move operation}
$\move_{ij}(D)$ recursively as
follows.
\begin{enumerate}[label=(\roman*)]
\item If $a=h$ and $b=k$, then define $\move_{ij}(D)=\ladder_{ij}(D)$,
where $(i,j)$ is ladder movable in $D$ by \Cref{lmove_conditions}.

    \item If $a>h$, define $\move_{ij}(D)=\move_{ij}\circ \move_{aj}(D)$.

    \item Otherwise, define $\move_{ij}(D)=\move_{ij}\circ \move_{ib}(D)$.
\end{enumerate}
\end{definition}

Observe that the tiles $(a,j)$ and $(i,b)$ in cases (ii) and (iii),
respectively, are movable in $D$. Therefore, the algorithm
above is well-defined and terminates for any cross tile $(i,j)$ that
is movable in $D$. Furthermore, if $D \in \rpd(w)$, then
$\move_{ij}(D) \in \rpd(w)$ and $\move_{ij}(D)>D$. We emphasize that $\move_{ij}(D)=\ladder_{ij}(D)$ if and only if $(i,j)$ is ladder movable in $D$. Before giving the explicit characterization of $\move_{ij}(D)$, we need another definition and a lemma.

\begin{definition}\label{path_defn}
     Let $D\in\rpd(w)$ such that $(i,j)$ is a movable cross tile in
     $D$. Define $\lpath_{ij}(D)$ to be the lattice path generated by the procedure below.
     \begin{enumerate}[label=(\roman*)]
         \item Let $(p,q)=(i,\rightbump_{ij}(D))$ and initialize $\lpath_{ij}(D)$ as $\{(p,q)\}$.
         \item Let $p'=\maxrow_{pj}(D)$ and trace up from $(p,q)$ to $(p',q)$ --- include $\{(p',q),(p'+1,q),\dots,(p-1,q)\}$ in $\lpath_{ij}(D)$.
         \item Let $q'=\min\{t\ge j\mid D(p',t)=\sbullet\}$ and trace left from $(p',q)$ to $(p',q')$ --- include $\{(p',q'),(p',q'+1),\dots,(p',q-1)\}$ in $\lpath_{ij}(D)$.
         \item If $q'=j$, return $\lpath_{ij}(D)$. Otherwise, set $(p,q)=(p',q')$ and return to step (ii).
     \end{enumerate}
     A tile $(p,q)\in\lpath_{ij}(D)$ is a \textit{corner} of $\lpath_{ij}(D)$ if either $\{(p-1,q),(p,q+1)\}\subseteq\lpath_{ij}(D)$ or $\{(p,q-1),(p+1,q)\}\subseteq\lpath_{ij}(D)$. The tiles $(i,\rightbump_{ij}(D))$ and $(\upperbump_{ij}(D),j)$ are the \textit{endpoints} of $\lpath_{ij}(D)$.
\end{definition}
\begin{definition}\label{def:shape}
    Let $D\in\rpd(w)$ such that $(i,j)$ is a movable cross tile in
    $D$. Define $\prtnshape_{ij}(D)$ to be the partition shape given by the tiles $(p,q)\in\rect_{ij}(D)$ such that there exists $p'\le p$ with $(p',q)\in\lpath_{ij}(D)$. In particular, if $(i,j)$ is ladder movable in $D$, then $\prtnshape_{ij}(D)=\rect_{ij}(D)$. 
\end{definition}
\begin{example}\label{M_example}
    In the submatrix of a pipe dream shown in \Cref{fig:move.bullets}, $\lpath_{61}(D)$ consists
    of the tiles in cyan, and $\prtnshape_{61}(D)$ consists of the
    tiles weakly below $\lpath_{61}(D)$. The cross tile $(6,1)$ is
    movable in $D$, and $\move_{61}(D)$ can be computed by
    \Cref{M_recursion}. We note that $\upperbump_{61}(D)=1$, $\rightbump_{61}(D)=8$, $\maxrow_{61}(D)=4$, and $\mincol_{61}(D)=2$.
\end{example}

\begin{figure}
\[
    \begin{NiceMatrix}[columns-width=auto]
        \textcolor{cyan}{\bullet} & \textcolor{cyan}{\bullet} & \textcolor{cyan}{\bullet} & + & \sbullet & + & \sbullet & + \\
        \textcolor{black}{+} & \textcolor{black}{+} & \textcolor{cyan}{+} & + & + & + & \sbullet & + \\
        \textcolor{black}{+} & \textcolor{black}{+} & \textcolor{cyan}{\bullet} & \textcolor{cyan}{+} & \textcolor{cyan}{\bullet} & \textcolor{cyan}{\bullet} & + & + \\
        \textcolor{black}{+} & \textcolor{black}{+} & \textcolor{black}{+} & \textcolor{black}{+} & \textcolor{black}{+} & \textcolor{cyan}{\bullet} & \textcolor{cyan}{+} & \textcolor{cyan}{\bullet} \\
        \textcolor{black}{+} & \textcolor{black}{+} & \textcolor{black}{+} & \textcolor{black}{+} & \textcolor{black}{+} & \textcolor{black}{+} & \textcolor{black}{+} & \textcolor{cyan}{+} \\
        \textcolor{black}{+} & \textcolor{black}{+} & \textcolor{black}{+} & \textcolor{black}{+} & \textcolor{black}{+} & \textcolor{black}{+} & \textcolor{black}{+} & \textcolor{cyan}{\bullet}
    \end{NiceMatrix}
\xrightarrow{\move_{61}}
\begin{matrix}
        \textcolor{cyan}{\sbullet} & + & + & + & \sbullet & + &  \sbullet & + \\
        \textcolor{cyan}{+} & + & + & + & + & + & \sbullet & + \\
        \textcolor{cyan}{\bullet} & + & + & + & + & + & + & + \\
        \textcolor{cyan}{\bullet} & + & + & + & + & + & + & + \\
        \textcolor{cyan}{+} & + & + & + & + & + & + & + \\
        \textcolor{cyan}{\bullet} & \textcolor{cyan}{\bullet} & \textcolor{cyan}{\bullet} & \textcolor{cyan}{+} & \textcolor{cyan}{\bullet} & 
        \textcolor{cyan}{\bullet} & \textcolor{cyan}{+} & \textcolor{cyan}{\sbullet}
    \end{matrix}
    \]
    \caption{A move operation is shown above. The diagram on the left is $\rect_{61}(D)$ for a pipe dream $D$. The path $\lpath_{61}(D)$ is colored cyan, and all bump tiles along the path are bold. The diagram on the right is the result of $\move_{61}(D)$, and new bump tiles along the southwest boundary of the shape are bold to highlight the changes. Note also the new cross tiles along $\lpath_{61}(D)$ in $\move_{61}(D)$.}
    \label{fig:move.bullets}
\end{figure}

\begin{lemma}\label{path_reduced_condition}
    Let $D\in\rpd(w)$ such that $(i,j)$ is a movable cross tile in $D$. If $(p,q)$ is a corner of $\lpath_{ij}(D)$, then $D(p,q)=\sbullet$.
\end{lemma}
\begin{proof}
    If $(p-1,q)$ and $(p,q+1)$ are both in $\lpath_{ij}(D)$, then $(p,q)$ corresponds to one of the tiles $(p',q')$ from step (iii) in the algorithm of \Cref{path_defn}, so $D(p,q)=\sbullet$ by construction.
    
     Now suppose that $(p,q-1),(p+1,q)\in\lpath_{ij}(D)$. Then $(p,q)$ must be a $(p',q)$ added in step (ii). Furthermore, by step (iii), there must be a column number $t<q$ such that $(p,t)\in\lpath_{ij}(D)$ and $D(p,t)=\sbullet$; let $t_{\max}$ be the largest such $t<q$. Similarly, by construction, there must be a row number $s>p$ such that $(s,q)\in\lpath_{ij}(D)$ and $D(s,q)=\sbullet$; let $s_{\min}$ be the smallest such $s>p$. By the choice of $t_{\max}$ and $s_{\min}$, $\rect_{s_{\min} t_{\max}}(D)=[p,s_{\min}]\times[t_{\max},q]$, and it follows from \Cref{path_defn} that the only bump tiles in this rectangle are $(p,t_{\max})$, $(s_{\min},q)$, and possibly $(p,q)$. However, if $D(p,q)=+$, then two pipes cross twice in $\rect_{s_{\min} t_{\max}}(D)$, contradicting the assumption that $D$ is reduced, see \Cref{fig:nonreduced_pipe_dream}. It follows that $D(p,q)=\sbullet$.
\end{proof}
\begin{proposition}\label{nonrecursive_m}
    Let $D\in\rpd(w)$ such that $(i,j)$ is a movable cross tile in $D$. Let $h=\upperbump_{ij}(D)$ and $k=\rightbump_{ij}(D)$. Then $\move_{ij}(D)$ is equal to the output of the procedure below.
\end{proposition}
    \begin{enumerate}[label=(\roman*)]
        \item For each $r\in(h,i]$, if there is a bump tile in row $r$ of $\lpath_{ij}(D)$, replace $(r,j)$ with a bump tile.
        \item For each $c\in[j,k)$, if there is a bump tile in column $c$ of $\lpath_{ij}(D)$, replace $(i,c)$ with a bump tile.
        \item Replace each bump tile in $\lpath_{ij}(D)$ with a cross tile, except for the endpoints $(h,j)$ and $(i,k)$.
    \end{enumerate}
    
\begin{figure}[h]
    \centering
\begin{tikzpicture}[scale=0.3]
  \foreach \x/\y in {1/0, 2/0, 3/0, 4/0, 5/0, 6/0, 7/0, 8/0, 9/0, 10/0, 11/0, 12/0, 13/0, 14/0, 15/0, 16/0, 17/0, 0/1, 1/1, 2/1, 3/1, 4/1, 5/1, 6/1, 7/1, 8/1, 9/1, 10/1, 11/1, 12/1, 13/1, 14/1, 15/1, 16/1, 17/1, 0/2, 1/2, 2/2, 3/2, 4/2, 5/2, 6/2, 7/2, 8/2, 9/2, 10/2, 11/2, 12/2, 13/2, 14/2, 15/2, 16/2, 0/3, 1/3, 2/3, 3/3, 4/3, 5/3, 6/3, 7/3, 8/3, 9/3, 10/3, 11/3, 12/3, 13/3, 14/3, 15/3, 16/3, 0/4, 1/4, 2/4, 3/4, 4/4, 5/4, 6/4, 7/4, 8/4, 0/5, 1/5, 2/5, 3/5, 4/5, 5/5, 6/5, 7/5, 8/5, 0/6, 1/6, 2/6, 3/6, 4/6, 0/7, 1/7, 2/7, 3/7, 4/7, 0/8, 1/8, 2/8, 3/8, 4/8, 0/9, 1/9, 0/0
  }{\fill[yellow!50] (\x,\y) rectangle ++(1,1); }
\draw[blue,   thick]   (0,0) -- (0,10) -- (2,10) -- (2,9) -- (5,9) -- (5,8) -- (5,8) -- (5,7) -- (5,7) -- (5,6) -- (9,6) -- (9,5) -- (9,5) -- (9,4) -- (17,4) -- (17,3) -- (17,3) -- (17,2) -- (18,2) -- (18,1) -- (18,1) -- (18,0) -- (0,0) -- (0,10) -- cycle;

  \tikzset{
    dot/.style={circle,fill=black,inner sep=1.5pt},
    ocdot/.style={circle,draw=black,inner sep=1.5pt},
    plus/.style={font=\large},
    arr/.style={->,>=latex,thick}
  }

  \foreach \x/\y in {
2/9,5/6,9/4,17/2,2/10,5/9,9/6,17/4,18/2,5/8,10/4,12/4
  }{
    \node[dot] at (\x,\y) {};
  }
  \foreach \x/\y in {
    0/0,0/9,0/6,0/4,0/2,2/0,5/0,9/0,17/0,0/8,10/0,12/0
  }{
    \node[plus] at (\x,\y) {+};
  }

  \foreach \x/\y in {
0/10,18/0
  }{
    \node[ocdot] at (\x,\y) {};
  }

   \draw[arr] (1.6,9)  -- (0.4,9);
   \draw[arr] (2,8.6)  -- (2,0.4);
   \draw[arr] (4.6,6)  -- (0.4,6);
   \draw[arr] (5,5.6)  -- (5,0.4);
   \draw[arr] (4.6,8)  -- (0.4,8); 
   \draw[arr] (8.6,4)  -- (0.4,4);
   \draw[arr] (9,3.6)  -- (9,0.4);
   \draw[arr] (10,3.6)  -- (10,0.4);
   \draw[arr] (12,3.6)  -- (12,0.4);
   \draw[arr] (16.6,2)  -- (0.4,2);
   \draw[arr] (17,1.6)  -- (17,0.4);

   \node at (-1,0) {\textit{i}};
   \node at (0,-1) {\textit{j}};
   \node at (-1,10) {\textit{h}};
   \node at (18,-1) {\textit{k}};
\end{tikzpicture}
    \caption{A depiction of the move operation $D \to \move_{ij}(D)$, where $(i,j)$ is the tile in the bottom left corner of the partition shape. Each tile that is pointed to by an arrow, along with $(i,j)$ itself, will become a $\sbullet$ tile after $\move_{ij}$. Each tile on $\lpath_{ij}(D)$, with the exception of the endpoints $(h,j)$ and $(i,k)$, will become a + tile.}
    \label{fig:oneshape}
\end{figure}

In \Cref{fig:oneshape}, we give a visualization of the operations described in \Cref{nonrecursive_m}. Note that one could have equivalently defined the move operators using this algorithm and then proved $\move_{ij}(D) \in \rpd(w)$ and $D<\move_{ij}(D)$. 

\begin{proof}
    Let $R_{ij}(D)$ denote the set of rows $r\in[h,i]$ such
that there is a bump tile in row $r$ of $\lpath_{ij}(D)$, and let
$C_{ij}(D)$ denote the set of columns $c\in[j,k]$ such that
there is a bump tile in column $c$ of $\lpath_{ij}(D)$. Let $\mathcal{N}_{ij}(D)$ denote the pipe dream obtained from the algorithm above. To show that
$\move_{ij}(D)=\mathcal{N}_{ij}(D)$, we induct on the
difference $\maxrow_{ij}(D)-h$ for $i>h$. See \Cref{M_computation} for a visualization of the induction process.

    In the base case $\maxrow_{ij}(D)=h$, let
$j=c_{1}<c_2<\dots < c_{m-1}<c_m=k$ be the column numbers in $C_{ij}(D)$. Then $D(h,k)=\sbullet$ by \Cref{path_reduced_condition}, so the bump tiles
in $\rect_{ij}(D)$ are exactly
$\{(h,j),(h,c_2),\dots,(h,c_{m-1}),(h,k)\}\cup\{(i,k)\}$. One can induct on
$m$ to show that \[\move_{ij}(D)=\ladder_{ij}\circ
\ladder_{ic_2}\circ\dots\circ \ladder_{ic_{m-1}}(D)\] by
\Cref{M_recursion}. Observe that the only bump tile in subrow
$\{h\}\times[j,k]$ of $\move_{ij}(D)$ is $(h,j)$, all tiles in the
rectangle $[h+1,i-1]\times[j,k]$ of $\move_{ij}(D)$ are crosses, and
the bump tiles in subrow $\{i\}\times[j,k]$ of $\move_{ij}(D)$ are
exactly the tiles $(i,c)$ for $c\in C_{ij}(D)$. Therefore, comparing
$\move_{ij}(D)$ with $\mathcal{N}_{ij}(D)$ tile-by-tile for
$(p,q)\in[h,i]\times[j,k]$ yields that $\move_{ij}(D)=\mathcal{N}_{ij}(D)$.
    
    Now consider the induction step, in which $a=\maxrow_{ij}(D)>h$ and $\move_{aj}(D)=\mathcal{N}_{aj}(D)$ by the induction hypothesis. Let $c_{\min}=\min\{c>j\mid D(a,c)=\sbullet\}$, so that $\rect_{aj}(D)=[h,a]\times[j,c_{\min}]$. Then $\lpath_{aj}(D)$ consists exactly of the elements in $\lpath_{ij}(D)$ with a column number between $j$ and $c_{\min}$. By the induction hypothesis, $\move_{aj}(D)$ has a bump tile at each $(p,j)$ for $p$ in $R_{aj}(D)=R_{ij}(D)\setminus\{i\}$, and at each $(a,q)$ with $q\in C_{aj}(D)$. Since $\move_{aj}(D)$ does not change any tiles in columns $q>c_{\min}$, $\move_{aj}(D)$ also has a bump tile at each $(a,q)$ with $q\in C_{ij}(D)$. Furthermore, all tiles in $\lpath_{aj}(D)$ above row $a$ are cross tiles in $\move_{aj}(D)$ except for $(h,j)$. Therefore, $\move_{aj}(D)$ agrees with $\mathcal{N}_{ij}(D)$ on all tiles in all rows $r<a$, and $C_{ij}(\move_{aj}(D))=C_{ij}(D)$, see \Cref{M_computation}. Since $\maxrow_{ij}(\move_{aj}(D))=\upperbump_{ij}(\move_{aj}(D))=a$, it follows from the work in the base case and \Cref{M_recursion} that $\move_{ij}\circ\move_{aj}(D)$ agrees with $\mathcal{N}_{ij}(D)$ on all tiles in all rows $r\in[a,i]$ as well. Since $\move_{ij}(D)$ and $\mathcal{N}_{ij}(D)$ agree with $D$ and therefore each other on all rows $r>i$, it follows that $\move_{ij}(D)=\move_{ij}\circ \move_{aj}(D)=\mathcal{N}_{ij}(D)$.
\end{proof}

\begin{example}\label{M_computation}
    Let $D$ be the pipe dream from \Cref{M_example}. By definition, $\move_{61}(D)=\move_{61}\circ \move_{41}(D)$. By \Cref{nonrecursive_m}, $\move_{41}(D)$ is given as below. The tiles of $\lpath_{41}(D)$ are shown in cyan on the left. The bottom and left boundary tiles of $\prtnshape_{41}(D)$, which have been altered by $\move_{41}$, are shown in cyan on the right.
    \[
    \begin{matrix}
        \textcolor{cyan}{\bullet} & \textcolor{cyan}{+} & \textcolor{cyan}{\bullet} & + & \sbullet & + & \sbullet & + \\
        + & + & \textcolor{cyan}{+} & + & + & + & \sbullet & + \\
        + & + & \textcolor{cyan}{\bullet} & \textcolor{cyan}{+} & \textcolor{cyan}{\bullet} & \textcolor{cyan}{\bullet} & + & + \\
        + & + & + & + & + & \textcolor{cyan}{\bullet} & + & \sbullet \\
        + & + & + & + & + & + & + & + \\
        + & + & + & + & + & + & + & \sbullet
    \end{matrix}\xrightarrow{\move_{41}}
    \begin{matrix}
        \textcolor{cyan}{\bullet} & + & + & + & \sbullet & + & \sbullet & + \\
        \textcolor{cyan}{+} & + & + & + & + & + & \sbullet & + \\
        \textcolor{cyan}{\bullet} & + & + & + & + & + & + & + \\
        \textcolor{cyan}{\bullet} & \textcolor{cyan}{+} & \textcolor{cyan}{\bullet} & \textcolor{cyan}{+} & \textcolor{cyan}{\bullet} & \textcolor{cyan}{\bullet} & + & \sbullet \\
        + & + & + & + & + & + & + & + \\
        + & + & + & + & + & + & + & \sbullet
    \end{matrix}
    \]
    We now compute $\move_{61}(\move_{41}(D))=\move_{61}\circ
    \move_{63}(\move_{41}(D))=\ladder_{61}\circ \ladder_{63}\circ
    \ladder_{65}\circ \ladder_{66}(\move_{41}(D))$ by \Cref{M_recursion}. The tiles of $\lpath_{61}(\move_{41}(D))$ are shown in cyan on the left, and the bottom and left boundary tiles of $\prtnshape_{61}(\move_{41}(D))$ are shown in cyan on the right.
    \[
    \begin{matrix}
        \sbullet & + & + & + & \sbullet & + &  \sbullet & + \\
        + & + & + & + & + & + & \sbullet & + \\
        \sbullet & + & + & + & + & + & + & + \\
        \textcolor{cyan}{\bullet} & \textcolor{cyan}{+} & \textcolor{cyan}{\bullet} & \textcolor{cyan}{+} & \textcolor{cyan}{\sbullet} & \textcolor{cyan}{\bullet} & \textcolor{cyan}{+} & \textcolor{cyan}{\bullet} \\
        + & + & + & + & + & + & + & \textcolor{cyan}{+} \\
        + & + & + & + & + & + & + & \textcolor{cyan}{\bullet}
    \end{matrix}\xrightarrow{\move_{61}\circ \move_{63}}
    \begin{matrix}
        \sbullet & + & + & + & \sbullet & + &  \sbullet & + \\
        + & + & + & + & + & + & \sbullet & + \\
        \sbullet & + & + & + & + & + & + & + \\
        \textcolor{cyan}{\bullet} & + & + & + & + & + & + & + \\
        \textcolor{cyan}{+} & + & + & + & + & + & + & + \\
        \textcolor{cyan}{\bullet} & \textcolor{cyan}{+} & \textcolor{cyan}{\bullet} & \textcolor{cyan}{+} & \textcolor{cyan}{\bullet} & 
        \textcolor{cyan}{\bullet} & \textcolor{cyan}{+} & \textcolor{cyan}{\bullet}
    \end{matrix}
    \]
\end{example}
We can derive some easy consequences of the move operations $\move_{ij}(D)$ directly from \Cref{nonrecursive_m}. The first corollary gives a special case in which the move operators commute, and the second corollary shows that there is a symmetry in the recursive definition of $\move_{ij}(D)$ between rows and columns.
\begin{corollary}\label{disjoint_commutativity}
    Let $D\in\rpd(w)$, and suppose that $(i,j),(p,q)$ are movable cross tiles in $D$ such that $\prtnshape_{ij}(D)$ and $\prtnshape_{pq}(D)$ are disjoint. Then $\move_{ij}\circ \move_{pq}(D)=\move_{pq}\circ \move_{ij}(D)$.
\end{corollary}
\begin{proof}
    This follows immediately from \Cref{nonrecursive_m}.
\end{proof}
\begin{corollary}\label{M'}
    Let $D\in\rpd(w)$ such that $(i,j)$ is a movable cross tile in $D$. Let $a=\maxrow_{ij}(D)$, $b=\mincol_{ij}(D)$, $h=\upperbump_{ij}(D)$, and $k=\rightbump_{ij}(D)$. Define $\move'_{ij}(D)$ recursively as follows:
\begin{enumerate}[label=\textnormal{(\roman*')}]
    \item If $a=h$ and $b=k$, define $\move'_{ij}(D)=\ladder_{ij}(D)$.

    \item If $b<k$, define $\move'_{ij}(D)=\move'_{ij}\circ \move'_{ib}(D)$.

    \item Otherwise, define $\move'_{ij}(D)=\move'_{ij}\circ \move'_{aj}(D)$.
\end{enumerate}
    Then $\move'_{ij}(D)=\move_{ij}(D)$.
\end{corollary}
\begin{proof}
    One can use a similar argument to the proof of \Cref{nonrecursive_m}, inducting on the difference $\rightbump_{ij}(D)-b$, to show that $\move'_{ij}(D)$ follows the same rule as $\move_{ij}(D)$ given in \Cref{nonrecursive_m}.
\end{proof}

\subsection{Minimality property}

Let $D\in\rpd(w)$, and let $\langle D\rangle$ be the upper order ideal generated by $D$, where the poset $(\rpd(w),<)$ has covering relations given by generalized ladder moves as defined in \Cref{def:ladder}. We will show that $\move_{ij}(D)$ is the unique minimal element of the set $V_{ij}(D)$, defined below.

\begin{definition}
    Let $D\in\rpd(w)$ and suppose that $D(i,j)=+$. Define
    \[
    V_{ij}(D)=\left\{
    Q\in \langle D\rangle\,\middle\vert
    \begin{array}{c}
    Q(i,j)=\sbullet,\, Q=\ladder_{p_mq_m}\circ \ladder_{p_{m-1}q_{m-1}}\circ\dots\circ \ladder_{p_1q_1}(D) \\
    \textnormal{ for some }(p_1,q_1),\dots,(p_m,q_m)\textnormal{ all northeast of }(i,j)
    \end{array}
    \right\}
    \]
where \textit{northeast} means $p_a\leq i$ and $q_a \geq j$ for all indices $a$.
\end{definition}
\begin{example}
    Let $D\in\rpd(126543)$ be the pipe dream below.
    \[
    D=\begin{NiceMatrix}
        \sbullet & \sbullet & \sbullet & \sbullet & \sbullet & \sbullet \\
        \sbullet & + & \sbullet & \sbullet & \sbullet \\
        + & + & + & \sbullet \\
        + & + & \sbullet \\
        \sbullet & \sbullet \\
        \sbullet
    \end{NiceMatrix}
    \]
    By definition, the pipe dream 
    \[
    \ladder_{32}\circ \ladder_{22}\circ \ladder_{33}(D)=\begin{NiceMatrix}
        \sbullet & \sbullet & + & \sbullet & \sbullet & \sbullet \\
        \sbullet & \sbullet & + & + & \sbullet \\
        + & \sbullet & \sbullet & \sbullet \\
        + & + & \sbullet \\
        \sbullet & \sbullet \\
        \sbullet
    \end{NiceMatrix}
    \]
    is in $V_{32}(D)$.
    On the other hand, one can verify that the pipe dream 
    \[
    \ladder_{32}\circ \ladder_{33}\circ \ladder_{22}\circ \ladder_{23}\circ \ladder_{41}(D)=\begin{NiceMatrix}
        \sbullet & \sbullet & + & + & \sbullet & \sbullet \\
        \sbullet & \sbullet & + & + & \sbullet \\
        + & \sbullet & \sbullet & \sbullet \\
        \sbullet & + & \sbullet \\
        \sbullet & \sbullet \\
        \sbullet
    \end{NiceMatrix}
    \]
    is greater than $D$ in $(\rpd(w),<)$ but is not in $V_{32}(D)$. Therefore,
    $V_{32}(D)$ is properly contained in the set of all $Q>D$ such that $Q(3,2)=\sbullet$.
\end{example}

\begin{lemma}\label{movable_conditions}
     Let $D\in\textnormal{RPD}(w)$ such that $D(i,j)=+$. Then tile $(i,j)$ is movable in $D$ if and only if $V_{ij}(D)$ is nonempty.
\end{lemma}
\begin{proof}
If $(i,j)$ is movable in $D$, then $\move_{ij}(D)\in V_{ij}(D)$ by \Cref{M_recursion}, and so $V_{ij}(D)$ is nonempty. On the other hand, if $(i,j)$ is not movable in $D$, then there is no $h<i$ such that $D(h,j)=\sbullet$, so there is no generalized ladder move on $D$ which could change any tile $(h,j)$ into a bump tile for $h\le i$. Therefore,
there is no sequence of generalized ladder moves which could turn $(i,j)$ into a bump tile, and in particular $V_{ij}(D)$ is empty.
\end{proof}

To show that $\move_{ij}(D)$ is the unique minimal element of $V_{ij}(D)$, we will use a double-induction argument on the upper order ideal $\langle D\rangle$ and on the southwest partial order of coordinates, where we consider $(i,j)\le_{sw} (i',j')$ if $(i,j)$ is southwest of $(i',j')$. In particular, we will often reference the induction hypotheses below.

\begin{definition}
    Suppose that $(i,j)$ is a movable cross tile in $D$. Then $D$ satisfies the $\move_{ij}$ \textit{induction hypotheses} if the following conditions hold:
    \begin{enumerate}[label=(\roman*)]
        \item For all $Q\in\langle D\rangle\setminus\{D\}$ and all $(p,q)$ such that $(p,q)$ is movable in $Q$, $\move_{pq}(Q)$ is the unique minimal element of $V_{pq}(Q)$.
        \item For all $(p,q)\neq(i,j)$ northeast of $(i,j)$ such that $(p,q)$ is movable in $D$, $\move_{pq}(D)$ is the unique minimal element of $V_{pq}(D)$.
    \end{enumerate}
\end{definition}

\begin{lemma}\label{laddermovablecompatibility}
    Let $D\in\textnormal{RPD}(w)$ such that $(i,j)$ is a ladder movable cross tile in $D$. Assume that $D$ satisfies the $\move_{ij}$ induction hypotheses. Then $\move_{ij}(D)=\ladder_{ij}(D)$ is the unique minimal element of $V_{ij}(D)$.
\end{lemma}
\begin{proof}
    We wish to show that, for all $Q\in V_{ij}(D)$,
$\ladder_{ij}(D)\le Q$, and hence $\ladder_{ij}(D)=\move_{ij}(D)$ is the unique minimal element of $V_{ij}(D)$. Given any $Q\in V_{ij}(D)$, there exists a sequence $\ladder_{p_1q_1},\ladder_{p_2q_2},\dots,\ladder_{p_mq_m}$ such that $Q=\ladder_{p_mq_m}\circ\dots \ladder_{p_2q_2}\circ\ladder_{p_1q_1}(D)$, where $(p_\ell,q_\ell)$ is northeast of $(i,j)$ for all $\ell$. Let $(p,q)=(p_1,q_1)$, and let $E=\ladder_{pq}(D)$. If $(p,q)=(i,j)$, then  $Q\ge\ladder_{ij}(D)$ immediately, so assume for the rest of the proof that $(p,q)\neq(i,j)$. 

If $(p,q)\notin\rect_{ij}(D)$, then since $(p,q)$ is northeast of $(i,j)$ and all tiles in $\prtnshape_{pq}(D)$ are northeast of $(p,q)$, $\rect_{ij}(D)$ and $\prtnshape_{pq}(D)$ are disjoint. Therefore, we have $\ladder_{ij}\circ \ladder_{pq}(D)=\ladder_{pq}\circ \ladder_{ij}(D)$ by \Cref{disjoint_commutativity}, where $(i,j)$ remains ladder movable in $E=\ladder_{pq}(D)$ since $\rect_{ij}(D)$ and $\rect_{pq}(D)$ are disjoint. Putting these observations together, we have
\[
\move_{ij}(E)=\ladder_{ij}(E)=\ladder_{ij}\circ\ladder_{pq}(D)=\ladder_{pq}\circ\ladder_{ij}(D)\ge \ladder_{ij}(D).
\]
Since $E>D$, we have from the first $\move_{ij}$ induction hypothesis that $\move_{ij}(E)$ is the unique minimal element of $V_{ij}(E)$. Therefore, $Q\ge  \move_{ij}(E) \ge \ladder_{ij}(D)$, as required.

For the remainder of the proof, suppose that $(p,q)\in\rect_{ij}(D)$. Let $h=\upperbump_{ij}(D)$ and $k=\rightbump_{ij}(D)$, so  we know $(h,k)$ is a bump tile in $D$ by \Cref{path_reduced_condition} and the hypothesis that $(i,j)$ is ladder movable. Since $(p,q)$ is also ladder movable by assumption, we claim that $q=k$ or $p=h$. Otherwise, if $j<q<k$, then every entry in column $q$ of $\rect_{ij}(D)$ is a cross tile by the ladder movability assumption, so $\upperbump_{pq}(D)<h$. Furthermore, the only bump tiles in $\rect_{ij}(D)\setminus\{(h,j)\}$ are in column $k$, so $\rightbump_{pq}(D)\ge k$ and hence $\rect_{pq}(D)$ contains the bump tile at $(h,k)$. Therefore, $\maxrow_{pq}(D)\ge h>\upperbump_{pq}(D)$, so $(p,q)$ is not ladder movable in $D$ by \Cref{lmove_conditions}. If $q=j$ and $h<p<i$, then one can show in a similar manner that $\upperbump_{pq}(D)=h$ and $\mincol_{pj}(D)=k<\rightbump_{pj}(D)$, so $(p,j)$ is not ladder movable in $D$, completing the proof of the claim.
    
We will now consider the case in which $q=k$ and $(p,k)$ is ladder movable in $D$ for $h<p<i$. By direct computation,
    \[
    \move_{ij}(E)=\ladder_{ij}\circ \ladder_{pj}(E)=\ladder_{ij}\circ \ladder_{pj}\circ \ladder_{pk}(D)=\ladder_{pj}\circ \ladder_{ij}(D)\ge \ladder_{ij}(D),\]
see \Cref{fig:ladder_move_combination}. Since $E=\ladder_{pk}(D)>D$, it follows from the first $\move_{ij}$ induction hypothesis that $\move_{ij}(E)$ is the unique minimal element of $V_{ij}(E)$. Therefore, $Q\ge\move_{ij}(E)\ge\ladder_{ij}(D)$, as required. The case in which $p=h$ and $(h,q)$ is ladder movable for $j<q<k$ is similar.
\end{proof}
\begin{figure}[h]
    \centering
    \begin{tikzcd}[ampersand replacement=\&]
        \begin{matrix}
            \sbullet & + & \sbullet & + & \sbullet \\
            + & + & + & + & +\\
            + & + & + & + & \sbullet \\
            + & + & \sbullet &
        \end{matrix}\arrow{r}{\ladder_{33}}\arrow{rd}{\ladder_{41}} \& \begin{matrix}
            \sbullet & + & \sbullet & + & +\\
            + & + & + & + & +\\
            + & + & \sbullet & + & \sbullet \\
            + & + & \sbullet &
        \end{matrix}\arrow{r}{\ladder_{31}} \&\begin{matrix}
            \sbullet & + & + & + & + \\
            + & + & + & + & + \\
            \sbullet & + & \sbullet & + & \sbullet \\
            + & + & \sbullet &
        \end{matrix}\arrow{d}{\ladder_{41}} \\
        \&\begin{matrix}
            \sbullet & + & + & + & \sbullet \\
            + & + & + & + & + \\
            + & + & + & + & \sbullet \\
            \sbullet & + & \sbullet &
        \end{matrix}\arrow{r}{\ladder_{31}}\&\begin{matrix}
            \sbullet & + & + & + & + \\
            + & + & + & + & + \\
            \sbullet & + & + & + & \sbullet \\
            \sbullet & + & \sbullet &
        \end{matrix}
    \end{tikzcd}
    \caption{Example computation in which $\ladder_{ij}\circ \ladder_{pj}\circ \ladder_{pk}(D)=\ladder_{pj}\circ \ladder_{ij}(D)$ for $i=4$, $j=1$, $p=3$, and $k=3$.}
    \label{fig:ladder_move_combination}
\end{figure}

\begin{lemma}\label{maxrowcompatibility}
    Let $D\in\textnormal{RPD}(w)$ such that $(i,j)$ is a movable cross tile in $D$. Assume that $D$ satisfies the $\move_{ij}$ induction hypotheses and $\maxrow_{ij}(D)>\upperbump_{ij}(D)$. Then $\move_{ij}(D)$ is the unique minimal element of $V_{ij}(D)$.
\end{lemma}
\begin{proof}
    Let $a=\maxrow_{ij}(D)$. Since $a>\upperbump_{ij}(D)$, we know $(i,j)$ is not ladder movable and $\move_{ij}(D)=\move_{ij}\circ\move_{aj}(D)$ by \Cref{M_recursion}. Every $Q\in V_{ij}(D)$ is obtained from $D$ by a sequence of ladder moves $Q=\ladder_{p_mq_m}\circ\dots\circ \ladder_{p_1q_1}(D)$, where each $(p_\ell,q_\ell)$ is northeast of $(i,j)$, and in particular $Q\in V_{ij}(\ladder_{p_1q_1}(D))$. Therefore, $Q\ge \move_{ij}\circ \ladder_{p_1q_1}(D)$ by the first $\move_{ij}$ induction hypothesis. If we can show that $ \move_{ij}\circ \ladder_{p_1q_1}(D) \geq \move_{aj}(D)$ for all ladder movable $(p_1,q_1)$ northeast of $(i,j)$, then we will have that $V_{ij}(D)=V_{ij}(\move_{aj}(D))$. Since $\move_{aj}(D)>D$, we have from the first $\move_{ij}$ induction hypothesis that $\move_{ij}\circ \move_{aj}(D)=\move_{ij}(D)$ is the unique minimal element of $V_{ij}(\move_{aj}(D))=V_{ij}(D)$. Therefore, if we assume $(p,q)$ is northeast of $(i,j)$ and ladder movable in $D$, then to prove the lemma it suffices to show 
    \begin{equation}\label{eq:mijleqmijlpq}
      \move_{ij}\circ \ladder_{pq}(D)\geq   \move_{aj}(D).
    \end{equation}
    If $(p,q)\notin\prtnshape_{ij}(D)$, then by \Cref{disjoint_commutativity},
    \[
    \move_{ij}\circ \ladder_{pq}(D)=\ladder_{pq}\circ \move_{ij}(D)=\ladder_{pq}\circ \move_{ij}\circ \move_{aj}(D)\ge \move_{aj}(D).
    \]
    Hence, one may further assume that $(p,q) \in \prtnshape_{ij}(D)$, $h=\upperbump_{ij}(D)$, and $k=\rightbump_{ij}(D)$. We consider several cases based on the relative order of $a,i,p$ given $a<i$ and $p\leq i$. 

    \begin{figure}[ht]
    \centering
    \[
    \begin{NiceMatrix}[columns-width=auto]
        \color{cyan}{\bullet} & \color{cyan}{\bullet} & \color{cyan}{\bullet} & \sbullet & \sbullet & \sbullet & \sbullet \\
        + & + & \color{cyan}{+} & \sbullet & \sbullet & \sbullet & \sbullet \\
        + & \textcolor{orange}{+} &  \color{cyan}{\bullet} & \color{cyan}{\bullet} & \color{cyan}{+} & \color{cyan}{\bullet} & \sbullet \\
        + & + & + & + & + & \color{cyan}{+} & + \\
        + & + & + & + & + & \textcolor{orange}{+} & \sbullet \\
        + & + & + & \textcolor{orange}{+} & + &
\color{cyan}{\bullet} & \sbullet \end{NiceMatrix} \]
\caption{A
section of a pipe dream $D$ in which the tile at $(6,1)$ is movable and
$a=\maxrow_{61}(D)=3>h_{61}(D)=1$, so $\move_{61}(D)=\move_{61}\circ
\move_{31}(D)$. The ladder movable tiles $(p,q)$ in
$\prtnshape_{61}(D)$ are $(3,2)$, $(5,6)$, and $(6,4)$, shown in
orange. The other tiles in $\lpath_{61}(D)$ are shown in cyan with bold bump tiles. }
\label{fig:maxrow_example}
\end{figure}

   \medskip
    \noindent
    \textbf{Case 1}: Assume $p=i$. The tile $(p,q)=(i,q)\in \prtnshape_{ij}(D)$ is a cross different from $(i,j)$ so $j< q<k$ and $k_{iq}(D)=k$, as is the case with  $(6,4)$ in \Cref{fig:maxrow_example}. By definition, $\rightbump_{aj}(D)$ is the minimal column number $t\in[j+1,k]$ such that $D(a,t)=\sbullet$. Therefore, if $q<\rightbump_{aj}(D)$, then $\upperbump_{iq}(D)<a$ and the rectangle $\rect_{iq}(D)=[h_{iq}(D),i]\times[q,k]$ must contain the bump tile $(a,\rightbump_{aj}(D))$ somewhere other than one of its corners,  contradicting the assumption that $(i,q)$ is ladder movable. Hence,  $q\ge\rightbump_{aj}(D)$ and $(a,\rightbump_{aj}(D))$ remains a bump tile in $\ladder_{iq}(D)$. It follows that $\maxrow_{ij}(\ladder_{iq}(D))=a>h$ by hypothesis. Therefore, by the recursive definition of $\move_{ij}$, we have $\move_{ij}(\ladder_{iq}(D))=\move_{ij}\circ \move_{aj}(\ladder_{iq}(D))$. Since $q\ge \rightbump_{aj}(D)$, $\prtnshape_{aj}(\ladder_{iq}(D))$ and $\rect_{iq}(D)$ can only overlap at the tile $(a,\rightbump_{aj}(D))$, which would necessarily be an endpoint along the corresponding paths. So by \Cref{nonrecursive_m}, $\move_{aj}\circ \ladder_{iq}(D)=\ladder_{iq}\circ \move_{aj}(D)$. Thus, \eqref{eq:mijleqmijlpq} holds whenever $p=i$ since 
    \[
    \move_{ij}\circ \ladder_{pq}(D)=\move_{ij}\circ \move_{aj}\circ \ladder_{iq}(D) = \move_{ij}\circ \ladder_{iq}\circ \move_{aj}(D) \ge \move_{aj}(D).
    \]

\medskip
    \noindent
    \textbf{Case 2}: Assume $a<p<i$. If $q=k$, then $(p,q)=(p,k)$ and $\maxrow_{ij}(\ladder_{pk}(D))=p$, as is the case with  $(5,6)$ in \Cref{fig:maxrow_example}. Hence, $\move_{ij}(\ladder_{pk}(D))=\move_{ij}\circ \move_{pj}(\ladder_{pk}(D))$. Observe $\upperbump_{pj}(\ladder_{pk}(D))=h<a$ and $\rightbump_{pj}(\ladder_{pk}(D))=k$. Therefore,  $\rect_{pj}(\ladder_{pk}(D))$ and $\rect_{ij}(D)$ agree on 
    rows $h,h+1,\dots,p-1$, so $\maxrow_{pj}(\ladder_{pk}(D))=a > h.$  It follows that $\move_{pj}(\ladder_{pk}(D))=\move_{pj}\circ\move_{aj}(\ladder_{pk}(D))$, so
    \[\move_{ij}\circ \ladder_{pk}(D)=\move_{ij}\circ\move_{pj}\circ\ladder_{pk}(D)=\move_{ij}\circ \move_{pj}\circ \move_{aj}\circ \ladder_{pk}(D).
    \]
     Furthermore, $\prtnshape_{aj}(D)$ and $\rect_{pk}(D)$ can only overlap at the tile $(a,k)$, which would necessarily be an endpoint along the corresponding paths, so by \Cref{nonrecursive_m}, we know $(a,k)$ is a bump tile in $\move_{aj}(D)$ and $\ladder_{pk}(D)$. Therefore, $\move_{aj}$ changes no tiles in $\rect_{pk}(D)$, and $\ladder_{pk}$ changes no tiles in $\prtnshape_{aj}(D)$, so $\move_{aj}\circ \ladder_{pk}(D)=\ladder_{pk}\circ \move_{aj}(D)$. Thus, \eqref{eq:mijleqmijlpq} holds since 
    \[
    \move_{ij}\circ \ladder_{pq}(D)=\move_{ij}\circ \move_{pj}\circ \move_{aj}\circ \ladder_{pk}(D)=\move_{ij}\circ \move_{pj}\circ(\ladder_{pk}\circ \move_{aj})(D)\ge \move_{aj}(D).
    \]
On the other hand, if $q<k$, then $\upperbump_{pq}(D)\le a$ and $\rightbump_{pq}(D)>k$ since 
$(p,q)\in\prtnshape_{ij}(D)$. Therefore, $\rect_{pq}(D)$ contains $[a,p]\times[q,k+1]$. However, we have  $D(a,k)=\sbullet$ by \Cref{path_reduced_condition}, which contradicts the assumption that  $(p,q)$ is ladder movable in $D$.

\medskip
    \noindent
    \textbf{Case 3}: Assume $p\leq a$. If $(p,q)=(a,j)$, then $\move_{aj}(D)=\ladder_{pq}(D)$ by definition, so 
 \[
\move_{ij}\circ \ladder_{pq}(D) = \move_{ij}\circ\move_{aj}(D) \geq \move_{aj}(D).
   \]
   If $(p,q)\neq(a,j)$, then  $\maxrow_{ij}(\ladder_{pq}(D))=a$. It follows that
   \[
   \move_{ij}\circ \ladder_{pq}(D)=\move_{ij}\circ \move_{aj}\circ \ladder_{pq}(D)\ge \move_{aj}\circ \ladder_{pq}(D).\]
   Since $(p,q)$ is northeast of and not equal to $(a,j)$, we have $\move_{aj}\circ \ladder_{pq}(D)\in V_{aj}(D)$. By the second $\move_{ij}$ induction hypothesis, $\move_{aj}(D)$ is the unique minimal element of $V_{aj}(D)$, so $\move_{aj}(D)\le \move_{aj}\circ \ladder_{pq}(D)\le \move_{ij}\circ \ladder_{pq}(D)$. Therefore, \eqref{eq:mijleqmijlpq} holds in all possible cases.
\end{proof}

\begin{lemma}\label{mincolcompatibility}
    Let $D\in\textnormal{RPD}(w)$ such that $(i,j)$ is a movable cross tile in $D$. Assume that $D$ satisfies the $\move_{ij}$ induction hypotheses, $\maxrow_{ij}(D)=\upperbump_{ij}(D)$, and $\mincol_{ij}(D)<\rightbump_{ij}(D)$. Then $\move_{ij}(D)$ is the unique minimal element of $V_{ij}(D)$.
\end{lemma}

\begin{proof}
    A parallel argument to the proof of \Cref{maxrowcompatibility} shows that $\move'_{ij}(D)$, which was shown to be equal to $\move_{ij}(D)$ in \Cref{M'}, is the unique minimal element of $V_{ij}(D)$.
\end{proof}

\begin{proposition}\label{M_minimal_characterization}
    Suppose that $D\in\textnormal{RPD}(w)$ and $(i,j)$ is movable in $D$. Then $\move_{ij}(D)$ is the unique minimal element of $V_{ij}(D)$.
\end{proposition}
\begin{proof}
    We induct on the southwest order of coordinate pairs and the size of the upper order ideal $\langle D\rangle$. Our induction hypotheses are precisely the $\move_{ij}$ induction hypotheses. In the base case where $(i,j)$ is movable in $D$ and there are no movable tiles $(p,q)\neq(i,j)$ northeast of $(i,j)$, we must have that $(i,j)$ is ladder movable in $D$ and the only ladder movable tile northeast of $(i,j)$. Hence $\ladder_{ij}(D)$ is the unique minimal element of $V_{ij}(D)$. The base case for the order ideal induction, in which $D$ is the maximal element $\Dtop(w)$ of $\rpd(w)$, is trivial, since no tiles are ladder movable in $\Dtop(w)$. The induction step follows directly from Lemmas \ref{laddermovablecompatibility}, \ref{maxrowcompatibility}, and \ref{mincolcompatibility}, where we showed that the $\move_{ij}$ induction hypotheses imply that $\move_{ij}(D)$ is the unique minimal element of $V_{ij}(D)$. 
\end{proof}

\subsection{Commutation rules for move operators}

In order to prove our join algorithm in the next section, we need to show that for any movable and southwest incomparable pairs $(i,j),(p,q)$, we have $\move_{ij}\circ \move_{pq}(D)=\move_{pq}\circ \move_{ij}(D)$. We begin by showing that $\move_{ij}\circ \move_{pq}(D)$ and $\move_{pq}\circ \move_{ij}(D)$ are well-defined whenever $(i,j)$ and $(p,q)$ are movable in $D$.

\begin{lemma}\label{lem:swincomprable.moves}
    Let $D\in\rpd(w)$, and let $(i,j),(p,q)$ be movable cross tiles in $D$. Suppose that $(i,j)$ and $(p,q)$ are southwest incomparable. Then $(i,j)$ is movable in $\move_{pq}(D)$ and $(p,q)$ is movable in $\move_{ij}(D)$.
\end{lemma}
\begin{proof}
    Assume without loss of generality that $(i,j)$ is northwest of $(p,q)$, so that $i<p$ and $j<q$. For any positive integer $t$, it follows from \Cref{nonrecursive_m} that subcolumn $[1,i]\times\{t\}$ contains a bump tile in $\move_{ij}(D)$ if and only if subcolumn $[1,i]\times\{t\}$ contains a bump tile in $D$. Therefore, since $p>i$, $(p,q)$ is movable in $\move_{ij}(D)$ if and only if it is movable in $D$. On the other hand, since $j<q$ and $\move_{pq}(D)$ agrees with $D$ on all columns $t<q$, $(i,j)$ is movable in $\move_{pq}(D)$ if and only if it is movable in $D$.
        
\end{proof}

\begin{lemma}\label{L_and_M_sw_commutativity}
    Let $D\in\rpd(w)$, let $(i,j)$ be a movable cross tile in $D$, and let $(p,q)$ be a ladder movable cross tile in $D$. If $p>i$ and $q\ge j$, then $(p,q)$ is ladder movable in $\move_{ij}(D)$, and $\ladder_{pq}\circ \move_{ij}(D)=\move_{ij}\circ \ladder_{pq}(D)$. Similarly, if $p\le i$ and $q<j$, then $\ladder_{pq}\circ \move_{ij}(D)=\move_{ij}\circ \ladder_{pq}(D)$.
\end{lemma}
\begin{figure}[h]
    \centering
    \[
    \begin{tikzcd}[ampersand replacement=\&]
        \begin{matrix}
            \sbullet & + & \sbullet &  &  &  &  \\
            + & + & \sbullet & + & \sbullet & + & \sbullet \\
            + & + & + & + & + & + & + & \\
            + & + & + & + & + & + & \sbullet \\
             &  & + & + & + & & \\
             &  & + & + & \sbullet &  & 
        \end{matrix}\arrow{r}{\move_{41}}\arrow{d}{\ladder_{63}}
        \& \begin{matrix}
            \sbullet & + & + &  &  &  &  \\
            \sbullet & + & + & + & + & + & + \\
            + & + & + & + & + & + & + & \\
            \sbullet & + & \sbullet & + & \sbullet & + & \sbullet \\
             &  & + & + & + & & \\
             &  & + & + & \sbullet &  & 
        \end{matrix}\arrow{d}{\ladder_{63}} \\
        \begin{matrix}
            \sbullet & + & \sbullet &  &  &  &  \\
            + & + & \sbullet & + & + & + & \sbullet \\
            + & + & + & + & + & + & + & \\
            + & + & + & + & + & + & \sbullet \\
             &  & + & + & + & & \\
             &  & \sbullet & + & \sbullet &  & 
        \end{matrix}\arrow{r}{\move_{41}}
        \& \begin{matrix}
            \sbullet & + & + &  &  &  &  \\
            \sbullet & + & + & + & + & + & + \\
            + & + & + & + & + & + & + & \\
            \sbullet & + & \sbullet & + & + & + & \sbullet \\
             &  & + & + & + & & \\
             &  & \sbullet & + & \sbullet &  & 
        \end{matrix}
    \end{tikzcd}
    \]
    \caption{A subsection of a pipe dream $D$ in which $(4,1)$ is movable, $(6,3)$ is ladder movable, and $\ladder_{63}\circ \move_{41}(D)=\move_{41}\circ \ladder_{63}(D)$.}
    \label{fig:ladder_and_M_example}
\end{figure}
\begin{proof}
    We consider the case in which $p>i$ and $q\ge j$; the other case is similar. Let $h=\upperbump_{pq}(D)$ and $k=\rightbump_{pq}(D)$. Note that $(h,q)$ and $(p,k)$ are bump tiles, since $(p,q)$ is ladder movable in $D$ by hypothesis. If $h>i$, then $\prtnshape_{ij}(D)$ and $\rect_{pq}(D)$ are disjoint, so $\move_{ij}\circ \ladder_{pq}(D)=\ladder_{pq}\circ \move_{ij}(D)$ by \Cref{disjoint_commutativity}. If $h=i$, then since $D(h,q)=\sbullet$, $\prtnshape_{ij}(D)$ and $\rect_{pq}(D)$ can have a non-empty intersection only if $(h,q)=(i,\rightbump_{ij}(D))$, so that the northwest corner of $\rect_{pq}(D)$ is also the southeast corner of $\prtnshape_{ij}(D)$. In this case, the tile $(h,q)$ remains a bump after both $\ladder_{pq}$ and $\move_{ij}$, so $\ladder_{pq}$ changes no tiles in $\prtnshape_{ij}(D)$ and $\move_{ij}$ changes no tiles in $\rect_{pq}(D)$ by \Cref{nonrecursive_m}. Hence $\move_{ij}\circ \ladder_{pq}(D)=\ladder_{pq}\circ \move_{ij}(D)$. Hence, we may suppose for the remainder of the proof that $h<i<p$. 
    
    By hypothesis, $(p,q)$ is ladder movable, so all tiles in $[h,i]\times[q,k]$ other than $(h,q)$ and $(h,k)$  are cross tiles. Since all corners and endpoints of $\lpath_{ij}(D)$ must be bumps by \Cref{path_reduced_condition}, the intersection of $\lpath_{ij}(D)$ and $[h,i]\times[q,k]$ must be contained in a single row if it is nonempty. Therefore, either $(h,q)$ and $(h,k)$ are both in $\lpath_{ij}(D)$, or neither $(h,q)$ nor $(h,k)$ are in $\lpath_{ij}(D)$. 
    
    If neither $(h,q)$ nor $(h,k)$ are in $\lpath_{ij}(D)$, then $\lpath_{ij}(D)$ contains no bumps in $\rect_{pq}(D)$, so by \Cref{nonrecursive_m}, $\move_{ij}$ changes no tiles in $\rect_{pq}(D)$. Similarly, $\ladder_{pq}$ changes no tiles in $\prtnshape_{ij}(D)$, so by \Cref{nonrecursive_m} $\ladder_{pq}\circ \move_{ij}(D)=\move_{ij}\circ \ladder_{pq}(D)$.

    Now assume both $(h,q)$ and $(h,k)$ are in $\lpath_{ij}(D)$. Let $R_{ij}(D)$ be the set of row numbers $r\in[\upperbump_{ij}(D),i]$ such that $\lpath_{ij}(D)$ contains a bump tile in row $r$, and let $C_{ij}(D)$ be the set of column numbers $c\in[j,\rightbump_{ij}(D)]$ such that $\lpath_{ij}(D)$ contains a bump tile in column $c$. By \Cref{nonrecursive_m}, $\move_{ij}(D)$ only changes tiles in $\prtnshape_{ij}(D)$, and it changes the tiles in $\prtnshape_{ij}(D)$ as follows:
    \begin{enumerate}[label=(\roman*)]
        \item for each $r\in R_{ij}(D)$, $\move_{ij}(D)(r,j)=\sbullet$,
        \item for each $c\in C_{ij}(D)$, $\move_{ij}(D)(i,c)=\sbullet$,
        \item for each $(s,t)\in\lpath_{ij}(D)\setminus\{(\upperbump_{ij}(D),j),(i,\rightbump_{ij}(D))\}$, $\move_{ij}(D)(s,t)=+$.
    \end{enumerate}
    Observe that there are no column numbers $c\in C_{ij}(D)$ such that $q<c<k$, since 
    all tiles in $[h,i]\times[q,k]$ other than $(h,q)$ and $(h,k)$  are cross tiles. Therefore, $\rect_{pq}(\move_{ij}(D))=[i,p]\times[q,k]$, where the only bump tiles in the rectangle are the corners $(i,q)$, $(i,k)$, and $(p,k)$. It follows that $(p,q)$ is still ladder movable in $\move_{ij}(D)$ by swapping the cross tile at $(p,q)$ with the bump tile at $(i,k)$, see \Cref{fig:ladder_and_M_example}.

    On the other hand, if $\ladder_{pq}(D)$ is performed before $\move_{ij}$, then the bump tile at $(h,k)$ will be swapped with the cross tile at $(p,q)$, so $C_{ij}(\ladder_{pq}(D))=C_{ij}(D)\setminus\{k\}$. Since there is a bump tile $(h,q)$ in $D$ west of $(h,k)$ and in the same row, and $(h,q),(h,k)$ are both in $\lpath_{ij}(D)$, we have $\lpath_{ij}(\ladder_{pq}(D))=\lpath_{ij}(D)$ as sets. Therefore, $\move_{ij}\circ \ladder_{pq}(D)$ is changed as follows:
    \begin{enumerate}[label=(\roman*)]
        \item for each $r\in R_{ij}(D)$, $\move_{ij}\circ \ladder_{pq}(D)(r,j)=\sbullet$,
        \item for each $c\in C_{ij}(D)\setminus\{k\}$, $\move_{ij}\circ \ladder_{pq}(D)(i,c)=\sbullet$,
        \item for each $(s,t)\in\lpath_{ij}(D)\setminus\{(\upperbump_{ij}(D),j),(i,\rightbump_{ij}(D))\}$, $\move_{ij}\circ \ladder_{pq}(D)(s,t)=+$,
        \item $\move_{ij}\circ \ladder_{pq}(D)(p,q)=\sbullet$.
       
    \end{enumerate}
    Comparing the outputs of $\ladder_{pq}\circ \move_{ij}(D)$ and $\move_{ij}\circ\ladder_{pq}(D)$ yields that $\ladder_{pq}\circ \move_{ij}(D)=\move_{ij}\circ \ladder_{pq}(D)$.
\end{proof}

\begin{proposition}\label{southwest_incomparable_moves}
    Let $D\in\textnormal{RPD}(w)$ and let $(i,j),(p,q)$ be movable cross tiles in $D$. Suppose that $(i,j)$ and $(p,q)$ are southwest incomparable. Then $\move_{ij}\circ \move_{pq}(D)=\move_{pq}\circ \move_{ij}(D)$.
\end{proposition}

\begin{proof}
    We will induct on the size of the upper order ideal $\langle D\rangle$, where the base case is trivial since $\Dtop(w)$ has no movable tiles. In the induction step, assume that for all $Q>D$ and all southwest incomparable pairs of movable tiles $(i',j')$ and $(p',q')$, $\move_{i'j'}\circ \move_{p'q'}(Q)=\move_{p'q'}\circ \move_{i'j'}(Q)$. Now let $(i,j)$ and $(p,q)$ be southwest incomparable, movable cross tiles in $D$, and assume for the rest of the argument that $(i,j)$ is northwest of $(p,q)$, i.e. $p>i$ and $q>j$. We will show that $\move_{pq}\circ \move_{ij}(D)\ge \move_{ij}\circ \move_{pq}(D)$, where both sides of the inequality are well-defined by \Cref{lem:swincomprable.moves}. By the symmetry from \Cref{M'}, the same argument will show that $\move_{ij}\circ \move_{pq}(D)\ge \move_{pq}\circ \move_{ij}(D)$, so $\move_{ij}\circ \move_{pq}(D)=\move_{pq}\circ \move_{ij}(D)$.
    
    In the case where $(i,j)$ is ladder movable, $\move_{ij}\circ \move_{pq}(D) = \ladder_{ij}\circ \move_{pq}(D)=\move_{pq}\circ \ladder_{ij}(D)$ by \Cref{L_and_M_sw_commutativity}, so we are done. Therefore, assume $\move_{ij}(D) = \move_{ij}\circ \ladder_{ab}(D)$ for some $(a,b)$ northeast of $(i,j)$. By the induction hypothesis on $Q>D$,
    \begin{equation}\label{eq:first_commutativity_step}
        \move_{pq}\circ \move_{ij}(D) = \move_{pq}\circ \move_{ij}\circ \ladder_{ab}(D)
        =(\move_{ij}\circ \move_{pq})\circ \ladder_{ab}(D).
    \end{equation}
    If $(a,b)$ is not northeast of $(p,q)$, then $\move_{pq}\circ \ladder_{ab}(D) = \ladder_{ab}\circ \move_{pq}(D)$ by \Cref{L_and_M_sw_commutativity}. By \Cref{M_minimal_characterization}, $\move_{ij}\circ \ladder_{ab}(\move_{pq}(D))\ge \move_{ij}(\move_{pq}(D))$, since $\move_{ij}\circ \ladder_{ab}(\move_{pq}(D))\in V_{ij}(\move_{pq}(D))$. Thus, by \Cref{eq:first_commutativity_step}, $\move_{pq}\circ \move_{ij}(D) = \move_{ij}\circ \ladder_{ab}\circ \move_{pq}(D) \ge \move_{ij}\circ \move_{pq}(D)$ whenever $(a,b)$ is not northeast of $(p,q)$.
    
    If $(a,b)$ is northeast of $(p,q)$, then $\move_{pq}\circ \ladder_{ab}(D)$ is in $V_{pq}(D)$, so $\move_{pq}\circ \ladder_{ab}(D)\ge \move_{pq}(D)$ by \Cref{M_minimal_characterization}. Furthermore, there exists a sequence $\ladder_{s_1t_1},\ladder_{s_2t_2},\dots,\ladder_{s_mt_m}$ such that $\move_{pq}\circ \ladder_{ab}(D)=(\ladder_{s_mt_m}\circ\dots \circ \ladder_{s_1t_1})\circ \move_{pq}(D)$ where all tiles $(s_\ell,t_\ell)$ are northeast of $(p,q)$. Hence $\move_{pq}\circ\ladder_{ab}(D)$ agrees with $D$ on all tiles that are not northeast of $(p,q)$.
    Since every pair $(s_\ell,t_\ell)$ is northeast of $(p,q)$, every $(s_\ell,t_\ell)$ is either northeast of $(i,j)$ or southwest incomparable to $(i,j)$. Let $P_0=\move_{pq}(D)$, and for $1\le \ell\le m$, let
    \[
    P_\ell=\ladder_{s_\ell t_\ell}(P_{\ell-1})=\ladder_{s_\ell t_\ell}\circ \ladder_{s_{\ell-1}t_{\ell-1}}\circ\dots\circ \ladder_{s_1t_1}\circ \move_{pq}(D).
    \]
    If $(s_\ell,t_\ell)$ is northeast of $(i,j)$, then 
    \[
    \move_{ij}(P_\ell)=\move_{ij}\circ \ladder_{s_\ell t_\ell}(P_{\ell-1})\ge \move_{ij}(P_{\ell-1}),
    \]
    since $\move_{ij}\circ \ladder_{s_\ell t_\ell}(P_{\ell-1})\in V_{ij}(P_{\ell-1})$ and $\move_{ij}(P_{\ell-1})$ is the unique minimum of $V_{ij}(P_{\ell-1})$. If $(s_\ell,t_\ell)$ is southwest incomparable to $(i,j)$, then we have from \Cref{L_and_M_sw_commutativity} that
    \begin{equation}\label{eq:final_commutativity_step}
    \move_{ij}(P_\ell)=\move_{ij}\circ \ladder_{s_\ell t_\ell}(P_{\ell-1})=\ladder_{s_\ell t_\ell}\circ \move_{ij}(P_{\ell-1})\ge \move_{ij}(P_{\ell-1}).
    \end{equation}
    Thus, for all $\ell\ge 1$, $\move_{ij}(P_\ell)\ge \move_{ij}(P_{\ell-1})$. It follows from \Cref{eq:first_commutativity_step} and the definition of $P_m$ that $\move_{pq}\circ \move_{ij}(D)=\move_{ij}(P_m)\ge \move_{ij}(P_0)=\move_{ij}\circ \move_{pq}(D)$ whenever $(a,b)$ is northeast of $(p,q)$. In conclusion, $\move_{pq}\circ \move_{ij}(D)\ge \move_{ij}\circ \move_{pq}(D)$ in all cases, as claimed.
\end{proof}

\section{Proof of Rubey's lattice conjecture}\label{sec:proof}

In \Cref{join_algorithm}, we give a recursive algorithm for computing the join $D_1\lor D_2$ for any pipe dreams $D_1,D_2\in\rpd(w)$. It will follow from the correctness of the algorithm that joins exist for all $D_1,D_2\in\rpd(w)$, so $\rpd(w)$ is a join-semilattice in the partial order induced by generalized chute and ladder moves. Recall from the background section that there is an order-reversing bijection between $\rpd(w)$ and $\rpd(w^{-1})$ given by the transpose $D\mapsto D^t$. Therefore, $D_1\land D_2=(D_1^t\lor D_2^t)^t$, so if $\rpd(w)$ is a join-semilattice for all $w$, then it is also a meet-semilattice for all $w$. We conclude as a corollary of the join algorithm that the poset $\rpd(w)$, which we named the Rubey lattice, is in fact a lattice.

In order to state our join algorithm, we choose a specific linear extension of the southwest partial order on $[n]\times[n]$, where we consider $(i,j)\le (i',j')$ if $(i,j)$ is southwest of $(i',j')$. The reader can observe that any linear extension of southwest order will suffice to prove \Cref{join_algorithm}, but it simplifies the exposition to make a concrete choice.

\begin{definition}\label{def:principal}
    Let $\prec$ be the linear extension of the southwest partial order defined by $(i,j)\prec(i',j')$ if $i>i'$ or $i=i'$ and $j<j'$. Let $D_1,D_2\in\rpd(w)$ such that $D_1\neq D_2$, and define the \textit{principal disagreement} between $D_1$ and $D_2$ to be the minimal tile $(i,j)$ with respect to $\prec$ such that $D_1(i,j)\neq D_2(i,j)$.
\end{definition}

\begin{theorem}\label{join_algorithm}
    Let $D_1,D_2\in\rpd(w)$. Then $D_1\lor D_2$ exists and is given by the output of the following recursive algorithm:
\end{theorem}
\begin{enumerate}[label=(\roman*)]

\item If $D_1=D_2$, return $D_1$.

\item If $D_1\neq D_2$, let $(i,j)$ be the principal disagreement between $D_1$ and $D_2$. If $D_1(i,j)=+$ and $D_2(i,j)=\sbullet$, then return $\move_{ij}(D_1)\lor D_2$. Otherwise, return $D_1\lor \move_{ij}(D_2)$.

\end{enumerate}

\begin{example}
    Let $D_1,D_2\in\rpd(126453)$ as below. We compute $D_1\lor D_2$; cf. Figure 4 in \cite{rubey2012maximal}.
    \[
    D_1=\begin{NiceMatrix}[columns-width=auto]
        \sbullet & \sbullet & + & \sbullet & \sbullet & \sbullet \\
        \sbullet & \sbullet & \sbullet & \sbullet & \sbullet \\
        + & + & + & \sbullet \\
        \sbullet & \sbullet & \sbullet \\
        + & \sbullet \\
        \sbullet
    \end{NiceMatrix},\qquad
    D_2=\begin{NiceMatrix}[columns-width=auto]
        \sbullet & \sbullet & \sbullet & + & + & \sbullet \\
        \sbullet & + & \sbullet & \sbullet & \sbullet \\
        \sbullet & + & \sbullet & \sbullet \\
        \sbullet & \sbullet & \sbullet \\
        + & \sbullet \\
        \sbullet
    \end{NiceMatrix}
    \]
    By \Cref{def:principal}, the principal disagreement between $D_1$ and $D_2$ is $(3,1)$, so by \Cref{join_algorithm}, $D_1\lor D_2=\move_{31}(D_1)\lor D_2.$ We compute $\move_{31}(D_1)$ using \Cref{nonrecursive_m}.
    \[
    \move_{31}(D_1)=\begin{NiceMatrix}[columns-width=auto]
        \sbullet & \sbullet & + & \sbullet & \sbullet & \sbullet \\
        \sbullet & + & + & + & \sbullet \\
        \sbullet & \sbullet & \sbullet & \sbullet \\
        \sbullet & \sbullet & \sbullet \\
        + & \sbullet \\
        \sbullet
    \end{NiceMatrix}
    \]
    The principal disagreement between $\move_{31}(D_1)$ and $D_2$ is $(3,2)$, so $D_1\lor D_2=\move_{31}(D_1)\lor D_2=\move_{31}(D_1)\lor \move_{32}(D_2).$ We again compute $\move_{32}(D_2)$ using \Cref{nonrecursive_m}.
    \[
    \move_{32}(D_2)=\begin{NiceMatrix}[columns-width=auto]
        \sbullet & \sbullet & + & + & + & \sbullet \\
        \sbullet & \sbullet & + & \sbullet & \sbullet \\
        \sbullet & \sbullet & \sbullet & \sbullet \\
        \sbullet & \sbullet & \sbullet \\
        + & \sbullet \\
        \sbullet
    \end{NiceMatrix}
    \]
    In the last step of the algorithm, the principal disagreement between $\move_{31}(D_1)$ and $\move_{32}(D_2)$ is $(2,2)$, and one can check that $\move_{22}\circ \move_{31}(D_1)=\move_{32}(D_2)$. Therefore, $D_1\lor D_2=\move_{22}\circ \move_{31}(D_1)=\move_{32}(D_2)$.
\end{example}

Before proving \Cref{join_algorithm}, we introduce notation for the function defined by the algorithm of \Cref{join_algorithm}. Then we can refer to this algorithm in the proofs below without assuming that it is equivalent to computing a join or that any joins exist in $\rpd(w)$. We develop several lemmas and combine them to give the proof toward the end of this section. 

\begin{definition}\label{def:J_function}
    Let $J:\rpd(w)\times\rpd(w)\to\rpd(w)$ be the function defined by the recursion of \Cref{join_algorithm}.
\end{definition}
\begin{lemma}\label{upper_bump_exists}
    Let $D_1,D_2\in\textnormal{RPD}(w)$ such that $D_1\neq D_2$. Let $(i,j)$ be the principal disagreement between $D_1$ and $D_2$, and assume that $D_1(i,j)=+$ while $D_2(i,j)=\sbullet$. Then $(i,j)$ is movable in $D_1$, and hence $J$ is a well-defined function.
\end{lemma}
\begin{proof}
    By the choice of $(i,j)$, $D_1$ and $D_2$ agree on all tiles $(p,q)$ southwest of and not equal to $(i,j)$, so the same pipe enters $(i,j)$ from the bottom in both $D_1$ and $D_2$. Let $p$ be the label of this pipe. Since $D_2(i,j)=\sbullet$, pipe $p$ exits $(i,j)$ to the right in $D_2$, so $p>j$. If $(i,j)$ is not movable in $D_1$, then every tile above $(i,j)$ in column $j$ is a cross tile. Hence $p=j$, contradicting the assumption that $D_1$ and $D_2$ represent the same permutation $w$.
\end{proof}

\begin{definition}\label{def:bar_M}
    For $(i,j)\in[n]^2$, let $\overline{\move}_{ij}:\rpd(w)\to\rpd(w)$ be defined as $\overline{\move}_{ij}(D)=\move_{ij}(D)$ whenever $(i,j)$ is movable in $D$, and $\overline{\move}_{ij}(D)=D$ otherwise.
\end{definition}
\begin{lemma}\label{lemma:J_facts}
    For all $D_1,D_2\in\rpd(w)$, we have $J(D_1,D_2)\ge D_1,D_2$. Furthermore, if $D_1\neq D_2$ and $(i,j)$ is the principal disagreement between $D_1$ and $D_2$, then $J(D_1,D_2)=J(\overline{\move}_{ij}(D_1),\overline{\move}_{ij}(D_2))$.
\end{lemma}
\begin{proof}
    The fact that $J(D_1,D_2)=J(\overline{\move}_{ij}(D_1),\overline{\move}_{ij}(D_2))$ follows immediately from Definitions \ref{def:J_function} and \ref{def:bar_M}. To show that $J(D_1,D_2)\ge D_1,D_2$, observe that the claim holds in the base case $D_1=D_2$ since $J(D_1,D_1)=D_1$. If $D_1\neq D_2$, then inducting on the number of recursive steps to compute $J$, we have $J(\overline{\move}_{ij}(D_1),\overline{\move}_{ij}(D_2))\ge\overline{\move}_{ij}(D_1)\ge D_1$ and $J(\overline{\move}_{ij}(D_1),\overline{\move}_{ij}(D_2))\ge\overline{\move}_{ij}(D_2)\ge D_2$, as desired.
\end{proof}

\begin{proof}[Proof of \Cref{join_algorithm}]
    Let $D_1,D_2\in\rpd(w)$. By definition, the join $D_1\lor D_2$ is the unique minimal element greater than both $D_1$ and $D_2$ when it exists. To show that $D_1\lor D_2$ exists and is equal to $J(D_1,D_2)$, we will induct on the size of the upper order ideal $\langle(D_1,D_2)\rangle$ in $\rpd(w)^2$ equipped with the product order. All cases in which $D_1=D_2$ are trivial, and in particular the base case $\#\langle (\Dtop(w),\Dtop(w))\rangle=1$ is trivial.
    
    For the induction step, assume that $D_1\neq D_2$. Let $(i,j)$ be the principal disagreement between $D_1$ and $D_2$, and assume without loss of generality that $D_1(i,j)=+$ and $D_2(i,j)=\sbullet$. For any $Q\ge D_1,D_2$, we have $Q>D_1$ or $Q>D_2$ since $D_1\neq D_2$. Therefore, either $Q\ge \ladder_{pq}(D_1)$ for some $(p,q)$ and hence $Q\ge J(\ladder_{pq}(D_1),D_2)=\ladder_{pq}(D_1)\lor D_2$ by the induction hypothesis; or $Q\ge \ladder_{st}(D_2)$ for some $(s,t)$ and hence $Q\ge J(D_1,\ladder_{st}(D_2))=D_1\lor \ladder_{st}(D_2)$. Thus, if we show that for any ladder movable tiles $(p,q)$ in $D_1$ and $(s,t)$ in $D_2$, we have both $\ladder_{pq}(D_1)\lor D_2\ge\move_{ij}(D_1)$ and $D_1\lor\ladder_{st}(D_2)\ge \move_{ij}(D_1)$, then we will have that $D_1\lor D_2$ exists and is given by $J(\move_{ij}(D_1),D_2)=\move_{ij}(D_1)\lor D_2$, as desired.

    We reduce the problem of checking $D_1\lor \ladder_{st}(D_2)\ge \move_{ij}(D_1)$ for all ladder movable $(s,t)$ in $D_2$ to the problem of checking $\ladder_{pq}(D_1)\lor D_2\ge \move_{ij}(D_1)$ for all ladder movable $(p,q)$ in $D_1$. Let $(s,t)$ be ladder movable in $D_2$, and suppose that $(i,j)\prec(s,t)$. Then the principal disagreement between $D_1$ and $\ladder_{st}(D_2)$ will still be $(i,j)$. Because $(D_1,\ladder_{st}(D_2))>(D_1,D_2)$, we have by induction on $\#\langle(D_1,D_2)\rangle$ that
    \[D_1\lor \ladder_{st}(D_2)=J(D_1,\ladder_{st}(D_2))=J(\move_{ij}(D_1),\ladder_{st}(D_2))\ge \move_{ij}(D_1).\]
    On the other hand, if $(s,t)\prec(i,j)$, then since $D_1$ and $D_2$ agree on all $(c,d)\prec(i,j)$, the principal disagreement between $D_1$ and $\ladder_{st}(D_2)$ will be $(s,t)$. By induction on $\#\langle (D_1,D_2)\rangle$,
    \[
        D_1\lor \ladder_{st}(D_2)=\move_{st}(D_1)\lor \ladder_{st}(D_2)\ge \move_{st}(D_1)\lor D_2,
    \]
    and from the recursive definition of $\move_{st}$ there is some $(s',t')$ such that $\move_{st}(D_1)\ge \ladder_{s't'}(D_1)$. Hence $D_1\lor \ladder_{st}(D_2)\ge \move_{st}(D_1)\lor D_2\ge \ladder_{s't'}(D_1)\lor D_2$.
    
    To conclude, it suffices to check that for any ladder move $\ladder_{pq}$, $\ladder_{pq}(D_1)\lor D_2\ge \move_{ij}(D_1)$. In \Cref{join_case_checks} below, we show that for any $(p,q)$ that is ladder movable in $D_1$ and not southwest of $(i,j)$, $J(\ladder_{pq}(D_1),D_2)=\ladder_{pq}(D_1)\lor D_2\ge \move_{ij}(D_1)$. If $(p,q)$ is ladder movable in $D_1$ and southwest of $(i,j)$, then we show in \Cref{join_check_sw} that $\ladder_{pq}(D_1)\lor D_2\ge \move_{ij}(D_1)$ by induction on $\#\langle(D_1,D_2)\rangle$, which will complete the proof of \Cref{join_algorithm}.
\end{proof}

\begin{lemma}\label{join_case_checks}
    Let $D_1,D_2\in\rpd(w)$ such that $D_1\neq D_2$, and let $(i,j)$ be the principal disagreement between $D_1$ and $D_2$. Assume that $D_1(i,j)=+$ and $D_2(i,j)=\sbullet$. Then for all tiles $(p,q)$ that are ladder movable in $D_1$ and not southwest of $(i,j)$,
    \[
    J(\ladder_{pq}(D_1),D_2)\ge \move_{ij}(D_1).
    \]
\end{lemma}
\begin{proof}
    We will consider separately the three cases of tile $(p,q)$'s position relative to $(i,j)$:
    \begin{enumerate}[label=\arabic*.]
        \item $p<i$ and $q<j$ (strictly north and strictly west),
        \item $p\le i$ and $q\ge j$ (northeast),
        \item $p>i$ and $q>j$ (strictly south and strictly east).
    \end{enumerate}
    
    \medskip
    \noindent
    \textbf{Case 1}: The principal disagreement between $\ladder_{pq}(D_1)$ and $D_2$ remains $(i,j)$, since $p<i$ and $\ladder_{pq}$ changes no tiles south of row $p$. Therefore,
    \begin{align*}
    J(\ladder_{pq}(D_1), D_2)&=J(\move_{ij}\circ \ladder_{pq}(D_1), D_2) \\
    &=J(\ladder_{pq}\circ \move_{ij}(D_1),D_2) &\textnormal{by \Cref{L_and_M_sw_commutativity}} \\
    &\ge \ladder_{pq}\circ \move_{ij}(D_1) &\textnormal{by \Cref{lemma:J_facts}}\\
    &\ge \move_{ij}(D_1)
    \end{align*}
    as desired.

    \medskip
    \noindent
    \textbf{Case 2}: If $(p,q)=(i,j)$, then $J(\ladder_{ij}(D_1),D_2)\ge\ladder_{ij}(D_1)=\move_{ij}(D_1)$ by \Cref{lemma:J_facts} and the definition of $\move_{ij}$. If $(p,q)\neq(i,j)$, then since $(p,q)$ is northeast of $(i,j)$, the principal disagreement between $\ladder_{pq}(D_1)$ and $D_2$ is again $(i,j)$. Furthermore, $\move_{ij}\circ \ladder_{pq}(D_1)\ge \move_{ij}(D_1)$ by \Cref{M_minimal_characterization}. Therefore,
    \[J(\ladder_{pq}(D_1),D_2)=J(\move_{ij}\circ \ladder_{pq}(D_1),D_2)\ge \move_{ij}\circ \ladder_{pq}(D_1)\ge \move_{ij}(D_1).\]
    
    \medskip
    \noindent
    \textbf{Case 3}: Let $m$ be the number of iterations of $J$ required to compute $J(\ladder_{pq}(D_1),D_2)$. Define a sequence $(D_1^{(\ell)},D_2^{(\ell)})$ as follows: let $(D_1^{(0)},D_2^{(0)})=(\ladder_{pq}(D_1),D_2).$ For $0\le\ell< m$, let $(p_\ell,q_\ell)$ be the principal disagreement between $D_1^{(\ell)}$ and $D_2^{(\ell)}$, and let
    \[
    (D_1^{(\ell+1)},D_2^{(\ell+1)})=\big(\overline{\move}_{p_\ell q_\ell}(D_1^{(\ell)}),\overline{\move}_{p_\ell q_\ell}(D_2^{(\ell)})\big).
    \]
    Observe that $J(D_1^{(\ell)},D_2^{(\ell)})=D_1^{(m)}=D_2^{(m)}$ for all $\ell$. By \Cref{nonrecursive_m}, $\overline{\move}_{p_\ell q_\ell}$ changes no tiles $(c,d)\prec(p_\ell,q_\ell)$, so by the definition of the principal disagreement, $(p_0,q_0)\prec(p_1,q_1)\prec\dots\prec(p_{m-1},q_{m-1})$. Observe that the principal disagreement $(p_0,q_0)$ between $D_1^{(0)}=\ladder_{pq}(D_1)$ and $D_2^{(0)}=D_2$ is $(p,q)$, since $p>i$ and $q>j$ in this case.
    Therefore, $\{\ell\mid p_\ell>i\}$ is non-empty, and hence $L=1+\max\{\ell\mid p_\ell> i\}$ is well-defined.

    We claim by induction on $\ell$ that $D_1^{(\ell)}$ and $D_2^{(\ell)}$ both agree with $D_1$ and $D_2$ on every tile in $S:=\{(c,d)\mid c\ge i,\ d\le j\}\setminus\{(i,j)\}$. This holds for $\ell=0$ because $D_1$ and $D_2$ agree on $S$, and $\ladder_{pq}$ changes no tile in $S$ since $q>j$. If the claim holds at step $\ell$, then the principal disagreement $(p_\ell,q_\ell)$ is not in $S$. By \Cref{nonrecursive_m}, an operation at a tile outside $S$ changes no tile of $S$, so the claim holds at step $\ell+1$.

    Since all tiles in rows strictly below row $i$ precede all tiles in rows weakly above row $i$, the indices with $p_\ell>i$ are exactly $\ell=0,1,\dots,L-1$. Furthermore, 
for each $\ell<L$, the claim above forces $q_\ell>j$. Thus no operation before step $L$ changes the tile $(i,j)$, so
    \[
    D_1^{(L)}(i,j)=D_1(i,j)=+
    \qquad\text{and}\qquad
    D_2^{(L)}(i,j)=D_2(i,j)=\sbullet.
    \]
    Hence $D_1^{(L)}\neq D_2^{(L)}$, so $L<m$ and $(p_L,q_L)$ is well-defined. By the choice of $L$, we have $p_L\le i$. The claim excludes $p_L=i$ and $q_L<j$, while $p_L<i$ or $p_L=i$ and $q_L>j$ would contradict the fact that $(i,j)$ is already a disagreement preceding $(p_L,q_L)$. Therefore $(p_L,q_L)=(i,j)$.

    Let
    \[
    \Phi =
    \begin{cases}
        \mathrm{id}, & L=1,\\
        \overline{\move}_{p_{L-1}q_{L-1}}\circ\dots\circ
        \overline{\move}_{p_1q_1}, & L>1.
    \end{cases}
    \]
    Each nonidentity factor in $\Phi$ is a move at a tile southwest incomparable to $(i,j)$, and identity factors commute trivially. Therefore,
    \begin{align*}
        J(\ladder_{pq}(D_1),D_2)
        &=D_1^{(m)} \ge D_1^{(L+1)} &\textnormal{by construction}\\
        &=\move_{ij}\circ \Phi\circ \ladder_{pq}(D_1) &\textnormal{by the choice of }L \\
        &=\Phi\circ \ladder_{pq}\circ \move_{ij}(D_1)
        &\textnormal{by \Cref{L_and_M_sw_commutativity} and \Cref{southwest_incomparable_moves}} \\
        &\ge \move_{ij}(D_1).
    \end{align*}
   
\end{proof}

\begin{lemma}\label{lem:SW_principal_disagreement}
    Let $D_1,D_2\in\rpd(w)$ such that $D_1\neq D_2$. Let $(i,j)$ be the principal disagreement between $D_1$ and $D_2$, and assume that $D_1(i,j)=+$. Let $(p,q)$ be southwest of $(i,j)$, and suppose that $(p,q)$ is a ladder movable cross tile in both $D_1$ and $D_2$. Then the principal disagreement between $\ladder_{pq}(D_1)$ and $\ladder_{pq}(D_2)$ is $(i,j)$.
\end{lemma}
\begin{proof}
    Suppose that $\ladder_{pq}$ acts on $D_\ell$ by exchanging the cross at $(p,q)$ with a bump tile at $(h_\ell,k_\ell)$ for $\ell=1,2$. If either $(h_1,k_1)\prec(i,j)$ or $(h_2,k_2)\prec(i,j)$, then since $D_1$ and $D_2$ agree on all tiles $(c,d)\prec(i,j)$, $\ladder_{pq}$ must swap the same pair of tiles in $D_2$ as in $D_1$, and hence $(h_2,k_2)=(h_1,k_1)\prec(i,j)$.  Hence, the principal disagreement between $\ladder_{pq}(D_1)$ and $\ladder_{pq}(D_2)$ is $(i,j)$. Observe that neither $(h_1,k_1)$ nor $(h_2,k_2)$ can be equal to $(i,j)$, since then $(h_1,k_1)$ and $(h_2,k_2)$ must both be $(i,j)$, but since $D_1(i,j)=+$, this contradicts the hypothesis that $D_1$ is reduced. Therefore, the only remaining case to check is that in which $(h_1,k_1)\succ(i,j)$ and $(h_2,k_2)\succ(i,j)$. In this case, $(p,q)$ is the only tile $(c,d)\prec(i,j)$ such that $D_\ell(c,d)\neq\ladder_{pq}(D_\ell)(c,d)$, so $\ladder_{pq}(D_1)$ and $\ladder_{pq}(D_2)$ agree on all tiles $(c,d)\prec(i,j)$. Therefore, the principal disagreement between $\ladder_{pq}(D_1)$ and $\ladder_{pq}(D_2)$ is $(i,j)$ in all cases. Furthermore, since $(h_1,k_1)\neq(i,j)$, $\ladder_{pq}(D_1)(i,j)=+$ in all cases.
\end{proof}

\begin{lemma}\label{double_L_move}
    Let $D_1,D_2\in\rpd(w)$ such that $D_1\neq D_2$. Let $(i,j)$ be the principal disagreement between $D_1$ and $D_2$, and assume that $D_1(i,j)=+$. Let $(p,q)$ be southwest of $(i,j)$, and suppose that $(p,q)$ is a ladder movable cross tile in both $D_1$ and $D_2$. Then $J(\ladder_{pq}(D_1),D_2)\ge \move_{ij}(D_1)$.
\end{lemma}
\begin{proof}
    The principal disagreement between $\ladder_{pq}(D_1)$ and $D_2$ is $(p,q)$, so $J(\ladder_{pq}(D_1),D_2)=J(\ladder_{pq}(D_1),\ladder_{pq}(D_2))$. By \Cref{lem:SW_principal_disagreement}, $(i,j)$ is the principal disagreement between $\ladder_{pq}(D_1)$ and $\ladder_{pq}(D_2)$, where $\ladder_{pq}(D_1)(i,j)=+$. Therefore, by \Cref{lemma:J_facts},
    \[
    J(\ladder_{pq}(D_1),\ladder_{pq}(D_2))=J(\move_{ij}\circ \ladder_{pq}(D_1),\ladder_{pq}(D_2))\ge \move_{ij}\circ \ladder_{pq}(D_1).
    \]
    We claim that $\move_{ij}\circ \ladder_{pq}(D_1)=\ladder_{pq}\circ \move_{ij}(D_1)$, so
    \[J(\ladder_{pq}(D_1),D_2)\ge \move_{ij}\circ \ladder_{pq}(D_1)\ge \move_{ij}(D_1),\]
    completing the proof.
    
    To prove the claimed equality, note that if $p=i$ or $q=j$, then $\move_{ij}\circ \ladder_{pq}(D_1)=\ladder_{pq}\circ \move_{ij}(D_1)$ by \Cref{L_and_M_sw_commutativity}, so we are done. Otherwise, we may assume $p> i$ and $q<j$, since $(p,q)$ is southwest of $(i,j)$. We will show that $\prtnshape_{ij}(D_1)$ and $\rect_{pq}(D_1)$ are disjoint, so $\move_{ij}\circ \ladder_{pq}(D_1)=\ladder_{pq}\circ \move_{ij}(D_1)$ by \Cref{disjoint_commutativity}. Assume to the contrary that $\prtnshape_{ij}(D_1)$ and $\rect_{pq}(D_1)$ intersect. Since $(p,q)$ is southwest of $(i,j)$, this is possible if and only if $(i,j)\in\rect_{pq}(D_1)$, so $\upperbump_{pq}(D_1)\le i$ and $\rightbump_{pq}(D_1)\ge j$. Since $D_1$ and $D_2$ agree below row $i$, $\upperbump_{pq}(D_2)\le i$ and $\rightbump_{pq}(D_2)\ge j$ as well. Therefore, $(i,j)\in\rect_{pq}(D_2)$. However, since $D_2(i,j)=\sbullet$, $(p,q)$ is ladder movable in $D_2$, and we are working in the case where $p>i$ and $q<j$, we have from \Cref{lmove_conditions} that $(i,j)$ is the northeast corner of $\rect_{pq}(D_2)$. Hence the corners $(i,q)$ and $(p,j)$ of $\rect_{pq}(D_2)$ are bumps as well. Since $(i,j)$ is the principal disagreement between $D_1$ and $D_2$, and $(i,q),(p,j)\prec(i,j)$, we must have $D_1(i,q)=D_2(i,q)=\sbullet$ and $D_1(p,j)=D_2(p,j)=\sbullet$. Therefore, $\upperbump_{pq}(D_1)=i$ and $\rightbump_{pq}(D_1)=j$. Since $(p,q)$ is ladder movable in $D_1$, this implies that $D_1(i,j)=\sbullet=D_2(i,j)$ by \Cref{path_reduced_condition}, contradicting the choice of $(i,j)$.
\end{proof}

\begin{lemma}\label{join_check_sw}
    Let $D_1,D_2\in\rpd(w)$ such that $D_1\neq D_2$, and assume that for all $(Q_1,Q_2)>(D_1,D_2)$ in the product order, $Q_1\lor Q_2$ exists and is equal to $J(Q_1,Q_2)$. Let $(i,j)$ be the principal disagreement between $D_1$ and $D_2$, and assume that $D_1(i,j)=+$. Then for any ladder movable cross tile $(p,q)$ in $D_1$ that is southwest of $(i,j)$, $\ladder_{pq}(D_1)\lor D_2\ge \move_{ij}(D_1)$.
\end{lemma}
\begin{proof}
    We claim that for all ladder movable $(p,q)$ in $D_1$ that are southwest of $(i,j)$, one of the following conditions must be true:
    \begin{enumerate}[label=(\roman*)]
        \item $\ladder_{pq}(D_1)\lor D_2\ge \move_{ij}(D_1)$,
        \item There exists $(s,t)$ southwest of $(i,j)$ such that $(s,t)$ is ladder movable in $D_1$, $(s-i)+(j-t)<(p-i)+(j-q)$, and $\ladder_{pq}(D_1)\lor D_2\ge \ladder_{st}(D_1)\lor D_2$.
    \end{enumerate}
    Therefore, we can induct on the $L^1$ distance $(p-i)+(j-q)$ between $(i,j)$ and $(p,q)$ to show that for all ladder movable $(p,q)$ southwest of $(i,j)$, $\ladder_{pq}(D_1)\lor D_2\ge \move_{ij}(D_1)$. Whenever $(p-i)+(j-q)$ is minimal among ladder movable tiles $(p,q)$ in $D_1$, $(p,q)$ cannot satisfy condition (ii), so we necessarily have $\ladder_{pq}(D_1)\lor D_2\ge \move_{ij}(D_1)$. It follows that the result holds in all base cases. In the induction step, where $(p-i)+(j-q)$ is not minimal, either $\ladder_{pq}(D_1)\lor D_2\ge \move_{ij}(D_1)$ or there exists $(s,t)$ satisfying condition (ii) above. However, if there exists $(s,t)$ satisfying condition (ii), then by induction on $(p-i)+(j-q)$, $\ladder_{pq}(D_1)\lor D_2\ge \ladder_{st}(D_1)\lor D_2\ge \move_{ij}(D_1).$ Therefore, $\ladder_{pq}(D_1)\lor D_2\ge \move_{ij}(D_1)$ in all cases.

    We now prove that for all ladder movable $(p,q)$ in $D_1$ that are southwest of $(i,j)$, $(p,q)$ satisfies (i) or (ii). Observe that the principal disagreement between $\ladder_{pq}(D_1)$ and $D_2$ is $(p,q)$, so $\ladder_{pq}(D_1)\lor D_2=\ladder_{pq}(D_1)\lor \move_{pq}(D_2)$ by the hypothesis on $(Q_1,Q_2)>(D_1,D_2)$. If $(p,q)$ is ladder movable in $D_2$, then we have $\ladder_{pq}(D_1)\lor D_2=\ladder_{pq}(D_1)\lor \ladder_{pq}(D_2)\ge \move_{ij}(D_1)$ by \Cref{double_L_move}. Therefore, assume that $(p,q)$ is not ladder movable in $D_2$. Then $\move_{pq}(D_2)$ is recursively defined as $\move_{pq}\circ \ladder_{st}(D_2)$ for some ladder movable $(s,t)$ that is northeast of $(p,q)$. Note that $(p,q)\neq(i,j)$ since $D_2(i,j)=\sbullet$ by assumption. We observe that
    \begin{equation}\label{eq:join_inequality}
        \ladder_{pq}(D_1)\lor D_2=\ladder_{pq}(D_1)\lor \move_{pq}(D_2)=\ladder_{pq}(D_1)\lor \move_{pq}\circ \ladder_{st}(D_2)\ge D_1\lor \ladder_{st}(D_2).
    \end{equation}
    If $(s,t)\succ(i,j)$, then the principal disagreement between $D_1$ and $\ladder_{st}(D_2)$ is $(i,j)$. Therefore,
    \begin{align*}
        \ladder_{pq}(D_1)\lor D_2 &\ge D_1\lor \ladder_{st}(D_2) &\textnormal{by \eqref{eq:join_inequality}} \\
        &=J(D_1,\ladder_{st}(D_2)) &\textnormal{by the hypothesis on }(Q_1,Q_2)>(D_1,D_2) \\
        &=J(\move_{ij}(D_1),\ladder_{st}(D_2)) &\textnormal{by the definition of $J$}\\
        &\ge \move_{ij}(D_1) &\textnormal{by \Cref{lemma:J_facts}.}
    \end{align*}

    If instead $(s,t)\prec(i,j)$, then the principal disagreement between $D_1$ and $\ladder_{st}(D_2)$ is $(s,t)$. Let $(s',t')$ be the first ladder movable tile from the recursive definition of $\move_{st}$ so that $\move_{st}(D_1)\ge\ladder_{s't'}(D_1)$. Note that either $(s',t')$ is not southwest of $(i,j)$ or $(s'-i)+(j-t')<(p-i)+(j-q)$. In either case,
    \begin{align*}
        \ladder_{pq}(D_1)\lor D_2 &\ge D_1\lor \ladder_{st}(D_2) &\textnormal{by \eqref{eq:join_inequality}}\\
        &=\move_{st}(D_1)\lor \ladder_{st}(D_2) &\textnormal{by hypothesis and the definition of $J$}\\
        &\ge \ladder_{s't'}(D_1)\lor D_2.
    \end{align*}
    If $(s',t')$ is not southwest of $(i,j)$, then $\ladder_{s't'}(D_1)\lor D_2\ge \move_{ij}(D_1)$ by \Cref{join_case_checks}, so $(p,q)$ satisfies condition (i) above. On the other hand, if $(s',t')$ is southwest of $(i,j)$, then $(p,q)$ satisfies condition (ii) above. Therefore, we have shown that all tiles $(p,q)$ which are ladder movable in $D_1$ and southwest of $(i,j)$ satisfy condition (i) or condition (ii), completing the proof.
\end{proof}

\begin{proof}[Proof of \Cref{thm:rubey.lattice}]
    Let $(\rpd(w), <)$ be the poset on reduced pipe dreams for $w$ with covering relations given by generalized ladder moves as defined in \Cref{def:rubey.lattice}. By definition, $(\rpd(w), <)$ is a lattice whenever there is a well-defined meet and join operation for any two elements in the poset. The algorithm to find the join of any two reduced pipe dreams for $w$ is given by \Cref{join_algorithm}. As mentioned at the beginning of \Cref{sec:proof}, the meet of any two elements in the poset is given by $D_1\land D_2=(D_1^t\lor D_2^t)^t$. 
\end{proof}

In the following corollary, we show that the choice of linear extension $\prec$ of the southwest partial order that we used to define the join computing algorithm does not matter. In order to make this more precise, we define southwest extremal disagreements below.

\begin{definition}\label{def:sw_extremal_disagreement}
    Let $D_1,D_2\in\rpd(w)$. A tile $(p,q)$ is a \textit{southwest extremal disagreement} between $D_1$ and $D_2$ if $D_1(p,q)\neq D_2(p,q)$ and for all $(s,t)\neq(p,q)$ southwest of $(p,q)$, we have $D_1(s,t)=D_2(s,t)$.
\end{definition}

Observe that the principal disagreement of \Cref{def:principal} is a southwest extremal disagreement. Furthermore, if $(i,j)$ and $(p,q)$ are distinct southwest extremal disagreements between $D_1$ and $D_2$, then $(i,j)$ and $(p,q)$ are southwest incomparable.

\begin{corollary}\label{cor:any.sw.dis}
    Let $D_1, D_2 \in \rpd(w)$, and let $(p,q)$ be any southwest extremal disagreement between $D_1$ and $D_2$. Then $D_1 \lor D_2=\overline{\move}_{pq}(D_1) \lor \overline{\move}_{pq}(D_2)$.
\end{corollary}
\begin{proof}
We show that $D_1 \lor D_2=\overline{\move}_{pq}(D_1) \lor \overline{\move}_{pq}(D_2)$ by inducting on the size of upper order ideals in $\rpd(w)^2$. In the base case $\#\langle(D_1,D_2)\rangle=1$, we have $D_1=D_2=\Dtop$, so there are no disagreements between $D_1$ and $D_2$, and the claim of the corollary is vacuously true.

Now suppose that $\#\langle(D_1,D_2)\rangle>1$ and $D_1\neq D_2$, and let $(i,j)$ be the principal disagreement between $D_1$ and $D_2$. If $(p,q)=(i,j)$, then we are done. If $(p,q)\neq(i,j)$, then it follows from the definitions of the principal disagreement and southwest extremal disagreements that $p<i$ and $q<j$, i.e. $(p,q)$ must be strictly north and strictly west of $(i,j)$. Since $\overline{\move}_{pq}$ does not change any tiles in rows $r>p$, the principal disagreement between $\overline{\move}_{pq}(D_1)$ and $\overline{\move}_{pq}(D_2)$ is $(i,j)$. It follows from \Cref{join_algorithm} that
\begin{equation}\label{eq:any.sw.dis}
\overline{\move}_{pq}(D_1)\lor\overline{\move}_{pq}(D_2)=\overline{\move}_{ij}\circ\overline{\move}_{pq}(D_1)\lor\overline{\move}_{ij}\circ\overline{\move}_{pq}(D_2).
\end{equation}
Since $\overline{\move}_{ij}$ changes no tiles in columns $c<j$, $(p,q)$ is a southwest extremal disagreement between $\overline{\move}_{ij}(D_1)$ and $\overline{\move}_{ij}(D_2)$, so
\begin{align*}
    D_1\lor D_2 &=\overline{\move}_{ij}(D_1)\lor\overline{\move}_{ij}(D_2) &\textnormal{by \Cref{join_algorithm}} \\
    &=\overline{\move}_{pq}\circ\overline{\move}_{ij}(D_1)\lor\overline{\move}_{pq}\circ\overline{\move}_{ij}(D_2) &\textnormal{by the induction hypothesis} \\
    &=\overline{\move}_{ij}\circ\overline{\move}_{pq}(D_1)\lor\overline{\move}_{ij}\circ\overline{\move}_{pq}(D_2) &\textnormal{by \Cref{southwest_incomparable_moves}} \\
    &=\overline{\move}_{pq}(D_1) \lor \overline{\move}_{pq}(D_2) &\textnormal{by \eqref{eq:any.sw.dis},}
\end{align*}
as required.

\end{proof}

\begin{figure}[h!]
	\begin{center}
\includegraphics[scale=0.5]{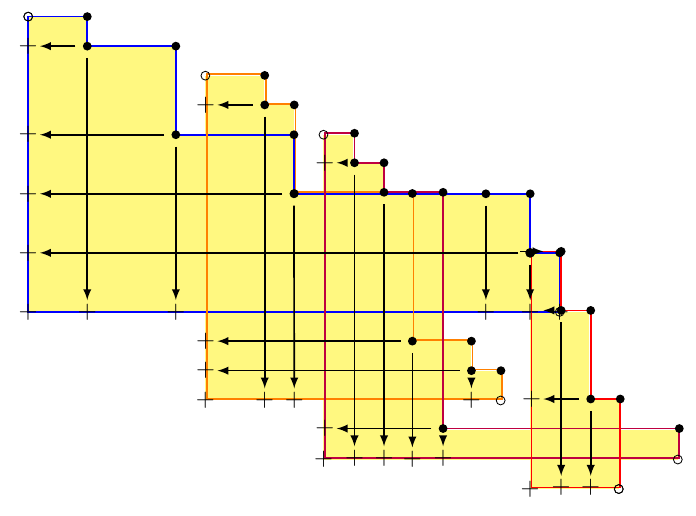}
\end{center}
\caption[]{The shaded region is the union of 4 overlapping partition shapes with incomparable southwest corners. The bump tiles that are changed to a + in the corresponding composition $\BigO_{(i,j) \in S}\move_{ij}(D_{1})$ are represented by a black dot, and the cross tiles that are changed to a bump are represented by a +. The arrows point to the +'s from the dots according to the move operations and \Cref{nonrecursive_m}. The arrows pass over the boundary of a partition, but they cannot pass through a hole in the union of partition shapes.}
\label{fig:polyomino}
\end{figure}

The proof of \Cref{cor:comprability.test} follows straightforwardly from
\Cref{join_algorithm} and \Cref{cor:any.sw.dis}. A depiction of such a composition of
southwest incomparable moves is shown in \Cref{fig:polyomino}.

\begin{proof}[Proof of \Cref{cor:comprability.test}]
Let $D_1,D_2\in\rpd(w)$, and let $S$ be the set of southwest extremal disagreements between $D_1$ and $D_2$. If $D_1(i,j)=+$ for all $(i,j)\in S$ and
\[D'=\BigO_{(i,j) \in S}\move_{ij}(D_{1}) \leq D_{2},\]
then $D_1 \leq D'\le D_2$. For the other direction, suppose that $D_1\le D_2$. If there is any $(i,j)\in S$ such that $D_2(i,j)=+$ and $D_1(i,j)=\sbullet$, then $D_1\lor D_2=D_1\lor\move_{ij}(D_2)$ by \Cref{cor:any.sw.dis}. However, $D_1\lor\move_{ij}(D_2)\ge \move_{ij}(D_2)>D_2$, contradicting the fact that $D_1\lor D_2=D_2$ since $D_1\le D_2$. Therefore, for all $(i,j)\in S$, $D_1(i,j)=+$ while $D_2(i,j)=\sbullet$. Observe that for all distinct $(i,j),(p,q)$ in $S$, $(p,q)$ is a southwest extremal disagreement between $\move_{ij}(D_1)$ and $D_2$. It follows from \Cref{cor:any.sw.dis} that
\[
D_2=D_1\lor D_2=\bigO_{(i,j)\in S}\move_{ij}(D_1)\lor D_2\ge\bigO_{(i,j)\in S}\move_{ij}(D_1),
\]
as desired.

\end{proof}

\section{Comparison Criterion}\label{sec:comparison}

In this section, we give an explicit criterion in terms of a filling of the diagram of $w$ to check if one given pipe dream is less than another given pipe dream in the Rubey lattice. The approach requires the reader to look beyond the cross and bump tiles to how a pipe traverses the tiles from the top to left boundary of the pipe dream. This perspective also yields some enumerative and structural properties of the Rubey lattice.

Recall from Section $2$ that pipes $x<y$ cross in any/every pipe dream for $w$ if and only if $(x,y)\in \Inv(w^{-1})$. Additionally, if pipe $x$ crosses pipe $y$ and $x<y$, then $x$ enters the crossing vertically and $y$ enters horizontally. Such crosses will be referred to as \emph{vertical} or \emph{horizontal crosses} of pipe $x$ or $y$ respectively. Furthermore, a bump tile of pipes $x$ and $y$ where pipe $x$ enters from the north and exits to the west is a \emph{left bump} of pipe $x$ and a \emph{down bump} of pipe $y$. Since each pipe traverses from the top to the left boundary of a pipe dream, each pipe must have at least one left bump. Observe from \Cref{def:ladder} that generalized ladder moves preserve the number of left bumps along a pipe, see \Cref{fig:gen_ladder_example}. Therefore, for any $x\in[n]$, we define \emph{$\B$} to be the number of \emph{left bumps of pipe $x$} in any reduced pipe dream $D\in\rpd(w)$. For $t\in \{1,2,\ldots, \B\}$, the \emph{$t^{th}$ left bump} of pipe $x$ will refer to the left bump of pipe $x$ with $t-1$ left bumps of pipe $x$ in rows to its south. The key object of this section will be a filling of the diagram of a permutation, $\textbf{D}(w)$, defined in Section 2.

\begin{definition}\label{def:pipedreamtab}
    Let $D\in \rpd(w).$ The \emph{pipe dream tableau} $\T(D)$ is a filling of $\textbf{D}(w)$ where $t$ is placed in cell $(w^{-1}(y),x)$ if there are $t$ left bumps of pipe $x$ in any row south of the crossing between pipes $x<y$. Denote the entry in cell $(w^{-1}(y), x)$ corresponding to the crossing between pipes $x < y$ of $\T(D)$ by $\T_D(x,y)$.
\end{definition}

\begin{definition}\label{def:pipedreamtabposet}
    Suppose $D_1, D_2\in \rpd(w)$. Define a partial order $\unlhd$ on pipe dream tableaux by $\T(D_1)\unlhd \T(D_2)$ if $\T_{D_1}(x,y)\le \T_{D_2}(x,y)$ for each $(x,y)\in \Inv(w^{-1})$.
\end{definition}

\begin{observation}\label{obs:bumprows}
    In any reduced pipe dream for $w$, pipe $x$ exits in row $w^{-1}(x)$ and no tiles with pipe $x$ are below row $w^{-1}(x)$. For each row $i< w^{-1}(x)$, there is either a vertical cross of pipe $x$ or a left bump followed by some number of horizontal crosses and a down bump. On row  $w^{-1}(x)$, $x$ must have a left bump.  
\end{observation}

\begin{lemma}
    If $(x,y)\in \Inv(w^{-1})$, then
    \begin{align}
    \T_{\Dtop(w)}(x,y)&=\B  = \#\{i \in [w^{-1}(x)] \mid w(i)\le x\}, \text{ and} \\
    \T_{\Dbot(w)}(x,y)&=\#\{i\in [w^{-1}(y),w^{-1}(x)] \mid w(i)\le x\}.
    \end{align}
\end{lemma}
\begin{proof}
    If pipes $x<y$ cross in $\Dtop(w)$, the crossing must be in column $x$, and every bump tile of pipe $x$ is south of the crossing since all cross tiles are top justified in $\Dtop(w)$. Hence, $\T_{\Dtop(w)}(x,y)=\B$. By \Cref{obs:bumprows}, the number of rows with a left bump of $x$ is
    $$
    \B = w^{-1}(x) - \#\{i\in [w^{-1}(x)] \mid w(i)>x \} = \#\{i\in [w^{-1}(x)] \mid w(i)\le x\}.
    $$

    Since crosses are left justified in $\Dbot(w)$, pipes $x<y$ cross in row $w^{-1}(y)$ where pipe $y$ exits. By \Cref{obs:bumprows}, the crossing between pipes $x<y$ has
    \begin{align*}
    w^{-1}(&x) - w^{-1}(y) - \#\{i\in [w^{-1}(y), w^{-1}(x)] \mid w(i)>x\}\\ &= \#\{i \in [w^{-1}(y), w^{-1}(x)] \mid w(i)\le x\}
    \end{align*}
    left bumps of pipe $x$ to its south and the result follows.
\end{proof}

\begin{figure}
\begin{center}
\includegraphics[scale=0.3]{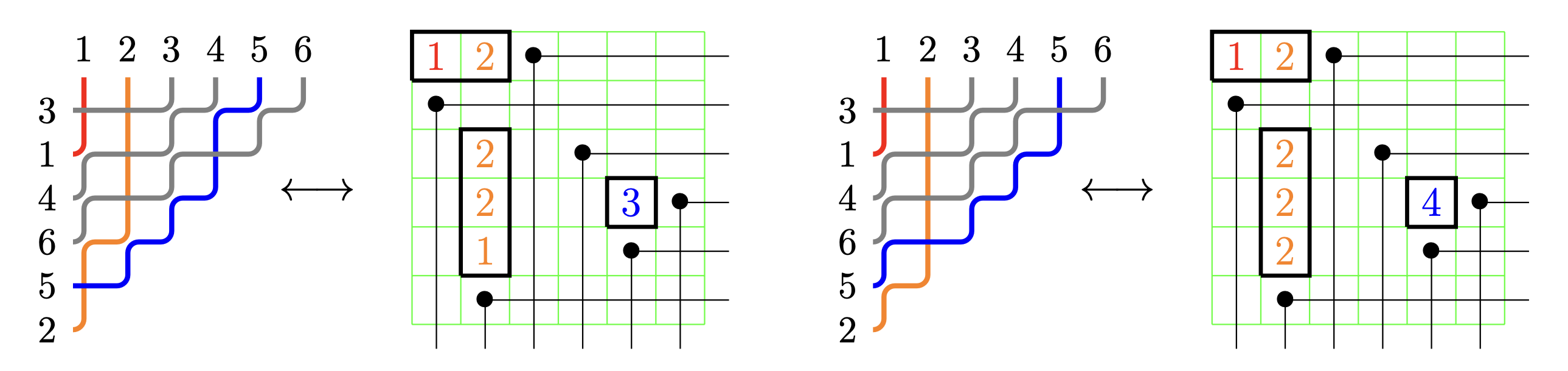}
\end{center}
\caption{Pipe dreams $D, \Dtop\in \rpd(314652)$ with their corresponding pipe dream tableaux on the left and right, respectively. The entries in $\T(D)$ and $\T(\Dtop)$ corresponding to vertical crosses of a given pipe are color coded for visual aid. Pipes with no vertical crosses are displayed in gray.}
\label{fig:pipe.tableau}
\end{figure}

\begin{example}\label{ex:pipe.tableau}
Let $w = 314652$. In \Cref{fig:pipe.tableau}, we show two pipe dreams for $w$. The one on the right is $\Dtop(w)$, so $\T(\Dtop(w))$ has $\B$ in each entry of column $x$. The pipe dream on the left is neither $\Dtop(w)$ nor $\Dbot(w)$. Observe that $\T_D(5,6) = 3$, since there are three left bumps of pipe $5$ south of the cross between pipes $5<6$ in $D$. Therefore, since $w^{-1}(6) = 4$, the value $\T_{D}(5,6)=3$ is entered in cell $(4,5)$. One can similarly check the values for the rest of the entries in $\T(D)$. We observe that both $D\le \Dtop(w)$ and $\T(D)\unlhd\T(\Dtop(w))$ in agreement with \Cref{thm:comprable}.

\end{example}

\begin{lemma}\label{lem:T(D)}
    A pipe dream $D\in \rpd(w)$ is uniquely determined by $\T(D)$.
\end{lemma}
\begin{proof}
       We will reconstruct $D$ uniquely from $\T(D)$ by placing the tiles of each pipe in order. For each pipe $x$, $D$ must necessarily have $$\#\{(x,y)\in \Inv(w^{-1}) \mid \T_D(x,y) =\B\}$$ vertical cross tiles followed by a left bump tile of pipe $x$ at the top of column $x$, so place all of these tiles to start the reconstruction. In particular, all tiles of pipe $1$ have been placed since pipe $1$ reaches the left boundary as it exits its first and only left bump. For $t>1$ and $x>1$, suppose that all tiles of pipes $1,...,x-1$ and down to the $t^{th}$ left bump tile of pipe $x$ have been placed uniquely. Assume the $t^{th}$ left bump tile of $x$ is in position $(i,k)$. A horizontal cross of pipe $x$ is a vertical cross of some pipe in $\{1,2,...,x-1\}$, so the left bump at $(i,k)$ will connect pipe $x$ to some number of previously placed horizontal crosses until one of the following possibilities occurs: pipe $x$ reaches the left boundary, pipe $x$ enters a previously placed down bump in a position $(i,j)$, or pipe $x$ meets an empty position $(i,j)$ where we necessarily place a down bump tile. If pipe $x$ meets a bump tile at $(i,j)$, $D$ must necessarily have $$\#\{(x,y) \in \Inv(w^{-1})\mid \T_D(x,y) =t-1\}$$ vertical cross tiles directly below $(i,j)$ followed by the $(t-1)^{th}$ left bump tile of pipe $x$, so place these tiles in the diagram.  Continuing this way recovers the entire path of each pipe in $D$.
\end{proof}

\begin{lemma}\label{lem:tableau.chute}
    Suppose $(i,j)$ is ladder movable in $D$ and $D(i,j)$ is a vertical cross of pipe $x$. Let $h = h_{ij}(D)$ and $Y$ be the subset of pipes that cross pipe $x$ horizontally in rows $[h+1,i]$ of $D$. For all $(x',y')\in \Inv(w^{-1})$,
    \begin{align*}
        \T_{\ladder_{ij}(D)}(x',y') = 
        \begin{cases}
        \T_D(x',y') + 1 &\text{ if } x' = x \text{ and } y'\in Y\\
        \T_D(x',y') &\text{ otherwise.}
        \end{cases}
    \end{align*}
\end{lemma}
\begin{proof}
    Let $k = k_{ij}(D)$. Applying $\ladder_{ij}$ to $D$ moves the vertical cross of pipe $x$ from $(i,j)$ to $(h, k)$ and as a result pipe $x$ will have a left bump at $(i,k)$ rather than $(h,k)$, see \Cref{fig:gen_ladder_example}. Thus, the vertical crosses of pipe $x$ in rows $[h+1,i]$ of $D$ have one more left bump of pipe $x$ to their south in $\ladder_{ij}(D)$. No other left bump row positions are changed, so all other vertical crosses have the same number of left bumps of pipe $x$ to their south.
\end{proof}

 In the following two corollaries, we show how the framework of pipe dream tableaux yields upper and lower bounds on the number of reduced pipe dreams for a given permutation. To state the bounds, let   
     $$R_{xy}(w) = \T_{\Dtop(w)}(x,y) - \T_{\Dbot(w)}(x,y)=\#\{i< w^{-1}(y) \mid w(i)\le x\}$$
for $(x,y)\in\Inv(w^{-1})$.
\begin{corollary}\label{cor:max.chain.length}
    Let $w$ be a permutation. The maximum cardinality of a chain in the Rubey lattice is 
   \begin{equation}\label{eq:max.chain.length}
    1+\sum_{(x,y)\in\Inv(w^{-1})} R_{xy}(w).
   \end{equation}
\end{corollary}
\begin{proof}
    By \Cref{thm:chutes.and.ladders}, we know that $\Dbot(w)$ is the unique minimal element of $(\rpd(w), <)$ and that 
    for any $D<\Dtop(w)\in \rpd(w)$ there exists a ladder movable cross in $D$. Assume $(i,j)$ is a ladder movable cross in $D$ such that no cross northeast of $(i,j)$ is movable. Then, $h_{ij}(D) =i-1$. Therefore, by \Cref{lem:tableau.chute}, the pipe dream tableaux for $\ladder_{ij}(D)$ and $D$ differ only in the entry corresponding to the crossing $(i,j)$ and the value is increased by exactly $1$. Thus, there is a chain in $\rpd(w)$ of length
    $$1+\sum_{(x,y)\in\Inv(w^{-1})} R_{xy}(w)$$
    from $\Dbot(w)$ to $\Dtop(w)$.
    Since $\sum_{(x,y)\in\Inv(w^{-1})} \left(\T_{\ladder_{ij}(D)}(x,y) - \T_D(x,y)\right) \ge 1$ for any ladder movable $(i,j)$ in $D$ by \Cref{lem:tableau.chute}, we conclude such a chain has the maximum possible length of any chain in $(\rpd(w),<)$.
\end{proof}

Note that $\sum_{(x,y)\in \Inv(w^{-1})} R_{xy}(w)$ is equal to the total number of occurrences of the pattern $132$ in $w$. Thus, the lower bound in \Cref{cor:bounds} below is equivalent to the lower bound from \cite{Weigandt.2018}.

\begin{corollary}\label{cor:bounds}
    Let $w$ be a permutation. Then,
    $$
    1+\sum_{(x,y)\in \Inv(w^{-1})} R_{xy}(w)\, \le \,  \# \rpd(w)\, \le \prod_{(x,y)\in \Inv(w^{-1})} \left( R_{xy}(w)+1 \right).
    $$
\end{corollary}
\begin{proof}
    By \Cref{lem:T(D)}, each pipe dream $D\in\rpd(w)$ corresponds to a unique pipe dream tableau, $\T(D)$, where $\T_{\Dbot(w)}(x,y)\le \T_D(x,y) \le \T_{\Dtop(w)}(x,y)$ for each $(x,y)\in\Inv(w^{-1})$. Therefore, the upper bound holds. The lower bound follows directly from \Cref{cor:max.chain.length}.
\end{proof}

\begin{lemma} \label{lem:shape.tiles}
Suppose $(i,j)$ is movable in $D$. If pipe $x$ has a horizontal cross in $\prtnshape_{ij}(D)\setminus\lpath_{ij}(D)$, then pipe $x$ exits $\prtnshape_{ij}(D)$ through its left boundary. Similarly, if pipe $x$ has a vertical cross in $\prtnshape_{ij}(D)\setminus\lpath_{ij}(D)$, then pipe $x$ exits $\prtnshape_{ij}(D)$ through its bottom boundary.
\end{lemma}
\begin{proof}
    Suppose that there is a horizontal cross of pipe $x$ at $(i',j')$ in $\prtnshape_{ij}(D)\setminus\lpath_{ij}(D)$. It follows from \Cref{def:shape} that each tile in $\prtnshape_{ij}(D)\setminus\lpath_{ij}(D)$ is a cross tile, and in particular, all tiles in $\{i'\} \times [j,j']$ are cross tiles in $D$. Therefore, pipe $x$ exits $\prtnshape_{ij}(D)$ through the left boundary at $(i',j)$. The other case is similar.
\end{proof}

\begin{lemma}\label{lem:tableau.sw}
    Suppose $D_1, D_2\in \rpd(w)$, $\T(D_1)\unlhd \T(D_2)$, and $(i,j)$ is the principal disagreement between $D_1$ and $D_2$ as in \Cref{def:principal}. Then $D_1(i,j) = +$ and $D_2(i,j) = \sbullet$.
\end{lemma}
\begin{proof}
     Suppose for contradiction that $D_1(i,j) = \sbullet$ and $D_2(i,j) =+$. Assume $D_2(i,j)$ is a cross between pipes $x<y$ and that there are $t$ left bumps of pipe $x$ strictly south of row $i$ in $D_2$. 
     We can further assume $D_1(i,j)$ is a bump of pipe $x$ with pipe $y$ and there are $t$ left bumps of pipe $x$ strictly south of row $i$ in $D_1$, since $D_1$ and $D_2$ agree everywhere southwest of $(i,j)$ other than position $(i,j)$ and $D_1, D_2\in \rpd(w)$. Thus, $D_1$ has its $(t+1)^{st}$ left bump of pipe $x$ in row $i$ strictly east of column $j$ and its vertical cross of pipe $x$ with pipe $y$ is strictly north of row $i$. In particular, $\T_{D_1}(x,y)>t$. Therefore, $\T_{D_1}(x, y)>\T_{D_2}(x, y)$ contradicting the hypothesis. Thus, $D_1(i,j) = +$ and $D_2(i,j) = \sbullet$.
\end{proof}

We are now ready to prove  \Cref{thm:comprable} from the introduction, namely  $D_1\leq D_2$ if and only if $\T(D_1) \unlhd \T(D_2)$. In other words, we show that $\T$ is a poset isomorphism between the Rubey lattice and the poset on pipe dream tableaux for $w$ given in \Cref{def:pipedreamtabposet}.

\begin{proof}[Proof of \Cref{thm:comprable}]
    Assume $D_1\le D_2$. Since there is a sequence of ladder moves from $D_1$ to $D_2$, it follows from \Cref{lem:tableau.chute} that $\T(D_1) \unlhd \T(D_2)$.

    Conversely, assume that $\T(D_1) \unlhd \T(D_2)$. We will induct on the size of the upper order ideal $\langle D_1 \rangle$. If $\#\langle D_1\rangle = 1$, then $D_1 = \Dtop(w)$ and all entries in column $x$ of $\T(D_1)$ are $\B$. Since $\T(D_1)\unlhd \T(D_2)$, all entries of $\T(D_2)$ must also be maximal and $\T(D_1) = \T(D_2)$. In this case, and any other cases where $\T(D_1) = \T(D_2)$, \Cref{lem:T(D)} implies $D_1 = D_2$ as desired.
    
    For the induction step, assume $\#\langle D_1 \rangle >1$ and $\T(D_1)\lhd \T(D_2)$. Let $(i,j)$ be the principal disagreement between $D_1$ and $D_2$, as defined in \Cref{def:principal}. By \Cref{lem:tableau.sw}, $D_1(i,j) = +$ and $D_2(i,j) = \sbullet$. By \Cref{upper_bump_exists}, $(i,j)$ is movable in $D_1$. By \Cref{movable_defn} and \Cref{nonrecursive_m}, there exists some ladder movable $(i',j')\in\prtnshape_{ij}(D_1)$ not in $\lpath_{ij}(D_1)$. We will show that $\T(\ladder_{i'j'}(D_1))\unlhd \T(D_2)$. Then, since $D_1 < \ladder_{i'j'}(D_1)$, induction on $\#\langle D_1 \rangle$ implies $D_1< \ladder_{i'j'}(D_1)\le D_2$.

    Let $h' = h_{i'j'}(D_1)$ and $k' = k_{i'j'}(D_1)$, so $D_1(h',k')=\sbullet$. Let $x'<y'$ be the pipes crossing at $D_1(i',j')$. Say there are $t$ left bumps of pipe $x'$ to the south of row $i$ in $D_1$. Then, $\T_{D_1}(x',y') = t$ since $(i',j')\in \prtnshape_{ij}(D_1)\setminus \lpath_{ij}(D_1)$ and all tiles in $\prtnshape_{ij}(D_1)\setminus \lpath_{ij}(D_1)$ are crosses in $D_1$. Let $Y$ be the subset of pipes that cross pipe $x'$ horizontally in rows $[h'+1,i']$ of $D_1$. By \Cref{lem:tableau.chute} and the assumption $\T(D_1)\lhd \T(D_2)$, it suffices to show that $\T_{D_2}(x',y)\ge t+1$ for all $y\in Y$. 
    Since $D_1$ and $D_2$ agree on all rows strictly south of row
    $i$ by the choice of $(i,j)$, $D_2$ must also have $t$ left
    bumps of $x'$ strictly south of row $i$. Furthermore, if
    $D_2(i,j')$ is a bump, then it is a down bump of pipe $x'$. So,
    in $D_2$ the $(t+1)^{th}$ left bump of pipe $x'$ is in row $i$ and
    the vertical crosses of pipe $x'$ with pipes in $Y$ are strictly
    north of row $i$. Thus, the vertical crosses of pipe $x'$ with
    pipes in $Y$ have at least $t+1$ left bumps of $x'$ to their
    south, implying $\T_{\ladder_{i'j'}}(D_1)\unlhd \T(D_2)$. It remains to show that $D_2(i,j')$ must be a bump. 
    
    Suppose for contradiction that $D_2(i,j') = +$. 
    Pipe $x'$ exits tile $(i,j')$ to the south in $D_1$ by \Cref{lem:shape.tiles}. The same is true for $D_2$ since $D_1$ and $D_2$ agree strictly south of row $i$.
    Let $p>x'$ be the pipe crossing horizontally at $(i,j')$ in $D_2$. Then we know $\T_{D_2}(x',p)=t$ since $D_1$ and $D_2$ agree below row $i$. Since $D_2(i,j)$ is a bump, pipe $p$ must exit row $i$ in some column in the range $[j,j')$ in $D_2$, and pipe $p$ necessarily exits row $i$ via the same tile in $D_1$. Since all the tiles in $\prtnshape_{ij}(D_1)\setminus\lpath_{ij}(D_1)$ of $D_1$ are crosses, pipe $p$ enters $\prtnshape_{ij}(D_1)$ northwest of the bump tile at $(h',j')$ in $D_1$. Hence pipes $x'$ and $p$ cross north of the left bump of pipe $x'$ at $(h',k')$ in $D_1$, implying $\T_{D_1}(x',p)> t= \T_{D_2}(x',p)$. This contradicts the assumption that $\T(D_1)\unlhd \T(D_2)$, so $D_2(i,j') = \sbullet$ as required.
\end{proof}

\begin{corollary}\label{cor:join.tableau}
If $\T_{D_1}(x,y) = t_1$ and $\T_{D_2}(x,y) = t_2$, then for all
$(x,y)\in \Inv(w^{-1})$, 
$$\T_{D_1\vee D_2}(x,y) \ge \max\{t_1,t_2\} \text{ and } \T_{D_1\wedge
D_2}(x,y) \le \min\{t_1, t_2\}.$$   In
particular, if $\T_{D_1}(x,y)$ or $\T_{D_2}(x,y)$ equals $\B$, then $\T_{D_1\vee D_2}(x,y) = \B$, and if
$\T_{D_1}(x,y)$ or $\T_{D_2}(x,y)$ equals $\T_{\Dbot(w)}(x,y)$, then
we have $\T_{D_1\wedge D_2}(x,y) = \T_{\Dbot(w)}(x,y)$.
\end{corollary}
\begin{proof}
    The proof is immediate from \Cref{thm:comprable}.
\end{proof}

\begin{corollary}\label{cor:top.row}
    The set of cross tiles in the top row of $D_1\vee D_2$ is exactly the union of crosses in the top rows of $D_1$ and $D_2$, and the set of cross tiles in the first column of $D_1\wedge D_2$ is exactly the union of crosses in the first columns of $D_1$ and $D_2$.
\end{corollary}
\begin{proof}
We only prove the first claim since the second claim follows from the fact that $D_1\wedge D_2$ is given by $(D_1^t\vee D_2^t)^t$.
    If a cross of pipes $x<y$ is in $D_1$ at $(1,x)$, then $\T_{D_1}(x,y) = \B$. By \Cref{cor:join.tableau}, $\T_{D_1\vee D_2}(x,y) = \B$ implying $(D_1\vee D_2)(1,x) = +$ by the reconstruction algorithm in the proof of \Cref{lem:T(D)}. The same argument holds for a cross in the top row of $D_2$, so the set of cross tiles in the top row of $D_1\vee D_2$ must contain the union of crosses in the top rows of $D_1$ and $D_2$. 
    
    It remains to show that every cross tile on the first row of $D_1\vee D_2$ must be on the first row of either $D_1$ or $D_2$. If $D_1=D_2$, then $D_1\vee D_2 = D_1 = D_2$ and the claim follows. Otherwise, let $i_{min}$ be the smallest $i\in[n]$ such that row $i$ of $D_1$ and $D_2$ are not equal. Let $j_{max} = \max\{j\in [n] \mid D_1(i_{min},j)\ne D_2(i_{min},j)\}$. Without loss of generality, assume that $D_1(i_{min},j_{max})$ is a bump between pipes $x$ and $y$ and $D_2(i_{min},j_{max})$ is a cross between pipes $x$ and $y$. Then pipes $x<y$ cross in $D_1$ say at $(i,j)$ where $i>i_{min}$ and $j<j_{max}$. 
    If $\lpath_{ij}(D_1)$ does not intersect the top row of $D_1$, then the top rows of $\move_{ij}(D_1)$ and $D_1$ are the same.
    Otherwise, $h_{ij}(D_1)=1$, and it follows that $(i_{min}, j_{max})=(1,x)$. Since $h_{ij}(D_1)=1$, pipe $x$ crosses vertically through each tile in subcolumn $[2,i]\times\{j\}$ of $D_1$; and since pipe $x$ has a left bump at $(1,x)$, it follows that $x$ crosses horizontally through each tile in subrow $\{1\}\times[j+1,x-1]$ of $D_1$. In particular, each tile in $\{1\}\times[j,x-1]$ is a tile in $\lpath_{ij}(D_1)$ by \Cref{path_reduced_condition}, so $k_{ij}(D_1)\ge x$. However, since pipe $y$ crosses $x$ horizontally at $(i,j)$ in $D_1$ and $y$ has a down bump at $(1,x)$, pipe $y$ has a left bump in subrow $\{i\}\times[j+1,x]$. Hence $k_{ij}(D_1)\le x$, so $k_{ij}(D_1)=x$. Furthermore, since $k_{ij}(D_1)=x$, pipe $y$ must vertically cross each tile in subcolumn $[2,i-1]\times\{x\}$ of $D_1$, so each tile in $\lpath_{ij}(D_1)$ is a cross tile except for $(1,j)$, $(1,x)$, and $(i,x)$. Thus, $\move_{ij}(D_1) = \ladder_{ij}(D_1)$, and $\move_{ij}(D_1)$ has exactly one more cross tile in its top row than $D_1$ at $(1,x)$. Since $D_2(1,x) =+$, the union of crosses in the top rows of $D_1$ and $D_2$ is equal to the union of crosses in the top rows of $\move_{ij}(D_1)$ and $D_2$.
    
    By induction on the size of the upper order ideal $\langle (D_1, D_2)\rangle$, the top row of $\move_{ij}(D_1)\vee D_2$
    is the union of crosses in the top rows of $\move_{ij}(D_1)$ and $D_2$, or equivalently, the union of crosses in the top rows of $D_1$ and $D_2$. Since $D_1\vee D_2\le \move_{ij}(D_1)\vee D_2$, the set of crosses in the top row of $D_1\vee D_2$ must be a subset of the crosses in the top row of $\move_{ij}(D_1)\vee D_2$. Thus, the top row of $D_1\vee D_2$ is also the union of crosses in the top rows of $D_1$ and $D_2$. 
\end{proof}

\begin{corollary}
    There is an entry of column $x$ of $\T(D_1\vee D_2)$ equal to $\B$ if and only if the same is true in column $x$ of either $\T(D_1)$ or $\T(D_2)$. 
\end{corollary}
\begin{proof}
    If $\T_{D_1\vee D_2}(x,y) =\B$ for some $(x,y)\in\Inv(w^{-1})$, then there is a cross in the top row of $D_1\vee D_2$ at $(1,x)$. Hence, by \Cref{cor:top.row}, there is a cross at $(1,x)$ in $D_1$ or $D_2$. Thus, there is an entry of column $x$ equal to $\B$ in $\T(D_1)$ or $\T(D_2)$. Conversely, if $\T_{D_1}(x,y)$ or $\T_{D_2}(x,y)$ is $\B$, for some $(x,y)\in \Inv(w^{-1})$, then $\T_{D_1\vee D_2}(x,y) = \B$ by \Cref{cor:join.tableau}. 
\end{proof}

\begin{example}
The reduced pipe dreams $D_1,D_2, \text{ and } D_1\vee D_2\in \rpd(316542)$ are shown in \Cref{fig:toprowjoin}. The set of cross tiles in the top row of $D_1\vee D_2$ is exactly the union of cross tiles in the top rows of $D_1$ and $D_2$. Notice that the same is not true of the second rows as $(D_1\vee D_2)(2,3) = +$, but $D_1(2,3) = D_2(2,3) = \sbullet$. 
Additionally, this example witnesses that the set of bump tiles in the first column of $D_1 \vee D_2$ contains, but is not always exactly, the union of bump tiles in the first columns of $D_1$ and $D_2$. 
\end{example}

\begin{figure}[h]
    \[
    D_1=\begin{NiceMatrix}[columns-width=auto]
        + & + & \sbullet & + & \sbullet & \sbullet \\
        \sbullet & \sbullet & \sbullet & \sbullet & \sbullet \\
        + & + & + & \sbullet \\
        + & \sbullet & \sbullet \\
        + & \sbullet \\
        \sbullet
    \end{NiceMatrix}
    \hspace{4mm}
    D_2 = \begin{NiceMatrix}[columns-width=auto]
        + & + & \sbullet & \sbullet & \sbullet & \sbullet \\
        \sbullet & + & \sbullet & + & \sbullet \\
        + & + & \sbullet & \sbullet \\
        + & + & \sbullet \\
        \sbullet & \sbullet \\
        \sbullet
    \end{NiceMatrix}
    \]
    
\[
    D_1\vee D_2 = \begin{NiceMatrix}[columns-width=auto]
        + & + & \sbullet & + & \sbullet & \sbullet \\
        \sbullet & + & + & + & \sbullet \\
        \sbullet & + & \sbullet & \sbullet \\
        \sbullet & + & \sbullet \\
        \sbullet & \sbullet \\
        \sbullet
    \end{NiceMatrix}
\]
    \caption{Pipe dreams $D_1,D_2\in\rpd(316542)$ above and their join $D_1\vee D_2$ below.}
    \label{fig:toprowjoin}
\end{figure}

\begin{remark}\label{rem:northerly}
    One could alternatively define \emph{northerly} pipe dream tableaux to have cells filled by the number of left bumps north of a vertical cross and analogous results hold. Counting left bumps above a crossing gives an encoding where $Dtop(w)$ is the all zeros filling. Hence, the Rubey lattice reverses the entrywise order on these northerly pipe dream tableaux. Similarly, one could transpose the pipe dream before encoding as well, so there are four natural ways to consider pipe dream tableaux. 

    More efficiently, one can list the values $T_D(x,y)$ for
    all $\Inv(w^{-1})$ in lexicographic order. This linear format
    suffices to encode each $D \in \rpd(w)$ by \Cref{lem:T(D)}, provided $w$ is known. 
\end{remark}

\begin{figure}[h]
    \centering
\includegraphics[scale=0.3]{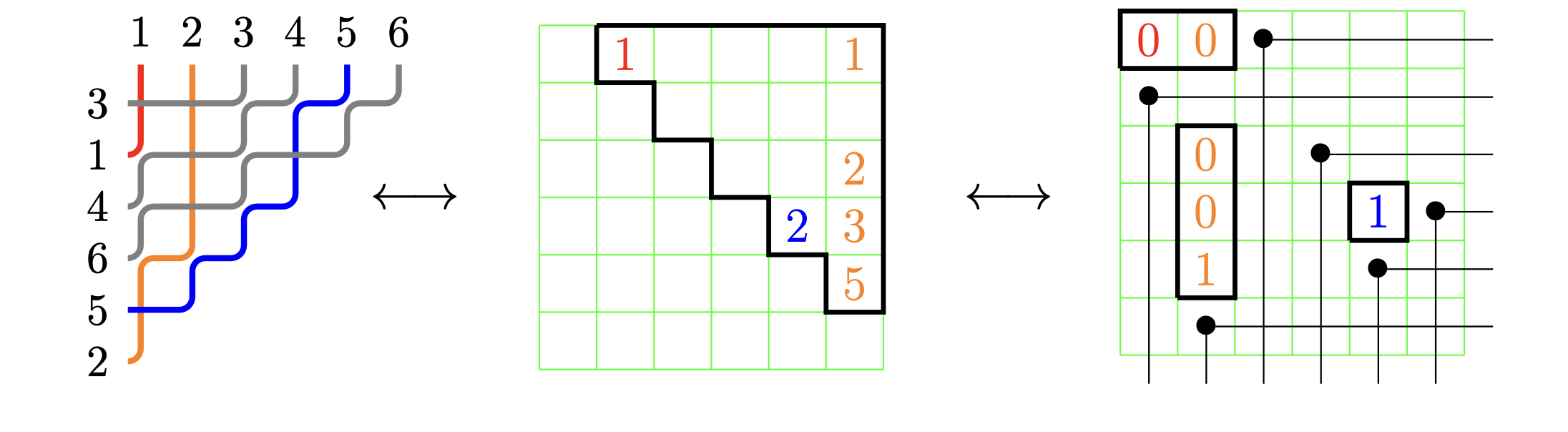}
    \caption{Let $w = 314652$. The pipe dream $D\in \rpd(w)$ on the left corresponds to the flagged inversion filling $F(D)$ in the middle and the northerly pipe dream tableau for $D$ on the right.}
    \label{fig:pipedreamtableauxbijs}
\end{figure}

\begin{remark}
    Pipe dream tableaux are closely related to several types of diagram fillings in the literature such as flagged balanced tableaux \cite{edelman1987balanced, fomin1997balanced}, flagged inversion fillings \cite[Ch.4]{Kelly.thesis}, $k$-flagged tableaux \cite{Serrano.Stump.2012}, perfect tableaux \cite{AdveRobichauxYong.2021}, and inversions and Lehmer tableaux \cite{axelrodfreed2025inversionstableaux,axelrodfreed2025chuteposetslattices}. 
    Any combinatorial objects that are naturally in bijection with reduced pipe dreams are examples of 
    \emph{Schubert objects}, as defined by Kelly in her 2007 thesis \cite{Kelly.thesis}.
    For example, the flagged inversion filling corresponding to $D\in \rpd(w)$, $F(D)$, records the row number
    of the crossing $(x,y) \in \Inv(w^{-1})$ in entry $(w^{-1}(y),w^{-1}(x))$ of an $n\times n$ array \cite[Rem. 4.3.10]{Kelly.thesis}. This definition is closely related to the one for inversions tableaux given by Axelrod-Freed 
    \cite[Def.3.12]{axelrodfreed2025inversionstableaux}.
    Similar to the map from inversions tableaux to Lehmer tableaux \cite[Def.6.13]{axelrodfreed2025inversionstableaux}, to go from $F(D)$ to a northerly pipe dream tableau for $D$, fill each $(w^{-1}(y),x)$ of
$\textbf{D}(w)$ with the number of positive integers that are smaller in
    value than $F(D)$ at $(w^{-1}(y),w^{-1}(x))$ and do not appear anywhere in column $w^{-1}(x)$ of $F(D)$. This is indeed the northerly pipe dream tableau for $D$ since there is a bump in row $i$ above the cross between $x<y$ if and only if there is no vertical cross of pipe $x$ in row $i$ by \Cref{obs:bumprows}. See \Cref{fig:pipedreamtableauxbijs} for an example.
\end{remark}

\section{Conclusion}\label{sec:future}

We conclude with a few open problems and remarks connecting this work
with related literature. We begin with the observation that the pipe
dream tableaux for $w \in S_{n}$ are in bijection with the reduced
pipe dreams for $w$, hence the Schubert polynomial is also
a generating function for these objects. Given a pipe dream tableau,
one can recover the pipe dream algorithmically as in the proof of
\Cref{lem:T(D)}. Can these objects be identified more explicitly?
Similarly, how can meets and joins in the Rubey lattice be computed
directly from pipe dream tableaux?

\begin{problem}
Characterize when a filling $T$ of the diagram of a permutation $w$ is
a pipe dream tableau for some $D\in \rpd(w)$, assuming $\T_{\Dbot(w)}(x,y)\leq
T(x,y) \leq \T_{\Dtop(w)}(x,y)$ for each $(x,y) \in \Inv (w^{-1})$.
\end{problem}

\begin{problem}
Give an explicit description of $\T(D_1\vee D_2)$ in terms of $\T(D_1)$ and $\T(D_2)$.
\end{problem}

In \cite{axelrodfreed2025chuteposetslattices}, the authors asked if there might be a description of the join and meet irreducibles in the Rubey lattice, where an element $D\in\rpd(w)$ is meet irreducible if and only if there exists a unique $(i,j)$ that is ladder movable in $D$. We give a more efficient characterization of the meet irreducibles in the Rubey lattice in terms of the statistics used in Section 3. 
\begin{lemma}
    Let $D\in\rpd(w)\setminus\{\Dtop(w)\}$. Then $D$ is meet irreducible if and only if there exist $H,K\in[n]$ such that for all movable cross tiles $(i,j)$ in $D$, $\upperbump_{ij}(D)=H$ and $\rightbump_{ij}(D)=K$. A similar statement holds for identifying join irreducibles in \rpd(w).
\end{lemma}
\begin{proof}
    Suppose that there exist $H,K\in[n]$ as described in the lemma. For any ladder movable $(i,j)$ in $D$, $H=\upperbump_{ij}(D)$ is the top row of $\rect_{ij}(D)$ and $K=\rightbump_{ij}(D)$ is the rightmost column of $\rect_{ij}(D)$ by hypothesis. Thus, any generalized ladder move will necessarily swap a cross tile $(i,j)$ southwest of $(H,K)$ with $(H,K)$. Since $D\neq\Dtop(w)$, there is at least one ladder-movable tile $(i,j)$ in $D$, so we must have $D(H,K)=\sbullet$. Furthermore, there can be at most one cross tile $(i,j)$ southwest of $(H,K)$ such that $(i,j)$ can be swapped with $(H,K)$ by a generalized ladder move, so there is at most one ladder-movable cross tile in $D$. Since the covering relations in $\rpd(w)$ are given by generalized ladder moves, $D$ has exactly one cover in $\rpd(w)$, so $D$ is meet irreducible.

    Conversely, suppose that there exist movable tiles $(i,j),(i',j')$ in $D$ with $(\upperbump_{ij}(D),\rightbump_{ij}(D))\neq(\upperbump_{i'j'}(D),\rightbump_{i'j'}(D))$. Let $\ladder_{st}$ be the first ladder move performed on $D$ in the $\move_{ij}$ recursion of \Cref{M_recursion}. It follows from the definition of the recursion that $\maxrow_{st}(D)=\upperbump_{ij}(D)$, so $\ladder_{st}$ swaps $(s,t)$ with a bump tile in row $\upperbump_{ij}(D)$. Similarly, if $\ladder_{s't'}$ is the first ladder move performed on $D$ in the $\move_{i'j'}$ recursion, then $\ladder_{s't'}$ swaps $(s',t')$ with a bump tile in row $\upperbump_{i'j'}(D)$. It follows that if $\upperbump_{ij}(D)\neq\upperbump_{i'j'}(D)$, then $D$ is covered by two distinct pipe dreams $\ladder_{st}(D)$ and $\ladder_{s't'}(D)$, so $D$ is not meet irreducible. On the other hand, if $\rightbump_{ij}(D)\neq\rightbump_{i'j'}(D)$, then one can use the $\move'$ recursion of \Cref{M'} to find generalized ladder moves $\ladder_{uv},\ladder_{u'v'}$ on $D$ such that $\ladder_{uv}$ swaps $(u,v)$ with a bump tile in column $\rightbump_{ij}(D)$ and $\ladder_{u'v'}$ swaps $(u',v')$ with a bump tile in column $\rightbump_{i'j'}(D)$. Then $\ladder_{uv}$ and $\ladder_{u'v'}$ must be distinct, so $D$ has two distinct covering relations, and therefore $D$ is not meet irreducible.
    
Transposing yields that $(j,i)$ is movable in $D^t$ if and only if there is a bump tile left of $(i,j)$ in $D$, i.e. $\{c<j:D(i,c)=\sbullet\}$ is nonempty. Therefore, $D$ is join irreducible in $\rpd(w)$ if and only if there exist $H,K\in[n]$ such that for all $(i,j)$ in $D$ with $\{c<j:D(i,c)=\sbullet\}$ nonempty, we have $\max\{c<j:D(i,c)=\sbullet\}=K$ and $\min\{r>i:D(r,j)=\sbullet\}=H$.
\end{proof}

In \cite{MPP.2019,MPPY.2025}, and recently \cite{anderson2026computationsamplingschubertspecializations},
many probabilistic aspects of pipe dreams were considered. They
address the problem of finding the most likely permutation given by
choosing a random pipe dream by independently inserting a cross or
bump tile in each of the $\binom{n }{2}$ positions of the triangle
for a given $n$, without regard to any fixed permutation. Such randomly generated pipe dreams can be
considered using both the usual product rule and the 
Demazure product for permutations related to Grothendieck polynomials.
Addressing a problem due to Stanley \cite{stanley2017schubertshenanigans}, it was shown in \cite{MPP.2019}
that
\begin{equation}
0.29 \leq \lim_{n \to \infty}\ \max_{w \in S_n}\ \frac{1}{n^{2}} \log_{2}
\mathfrak{S}_{w}(1,1,\ldots ,1) \leq 0.37,
\end{equation}
where $\mathfrak{S}_{w}(1,1,\ldots ,1)=\#\rpd(w)$ is the Schubert polynomial for $w$ evaluated at all $x_i=1$. Given the
exponential growth rate on the maximum number of reduced pipe dreams
for $w \in S_{n}$, it is challenging to choose one uniformly at random
even for $n=50$.

Toward the goal of choosing a uniformly random pipe dream for a fixed $w \in S_n$, we propose
using the following \textit{Markov process on the Rubey lattice} as a tool.
If $w$ has only one pipe dream, then the uniform process is trivial, so assume otherwise.  Start with $D=\Dbot(w)$. At each step, update $D$ as follows, assuming that $w$ is not the identity permutation. Uniformly at random, choose a cross tile of $D$. With probability $1/2$, attempt to apply a generalized ladder move to the chosen cross, and with probability $1/2$, attempt to apply a generalized chute move to the chosen cross. Such a move, when it exists, is uniquely determined by the cross. If the selected move exists, apply it and update $D$; otherwise, remain at $D$.  Observe that each generalized ladder move is the inverse of a generalized chute move, so the transition probabilities between any two adjacent reduced pipe dreams are equal in both directions.  The process is irreducible because the Hasse diagram of the Rubey lattice is connected, and it is aperiodic because every attempted chute move from $\Dbot(w)$ leaves the state unchanged.  Hence, this Markov process converges to the uniform distribution on $\rpd(w)$.   

\begin{problem} 
What is an upper bound, as a function of $n$, on the maximum mixing time of this Markov process over all permutations $w \in S_n$? 
\end{problem}

Alejandro Morales asked the following question. See \cite{MPPY.2025} for a review of the connection between nonreduced pipe dreams and Grothendieck polynomials. 

\begin{problem}
    Is there a natural lattice structure on the set of nonreduced 
pipe dreams associated with Grothendieck polynomials via the Demazure product?
\end{problem}

\section*{Acknowledgments} We thank Martin Rubey for introducing
the Rubey lattice. We would also like to thank Natasha Crepeau,
Elena Hafner, Cordelia Li, Ricky Liu, Bryan Lu, Alejandro Morales, Michael
Zeng,  and the anonymous referees for helpful conversations and suggestions on this paper. We thank
Yibo Gao for the pipe dream figure template used here and for
inspiration on this problem. 

\bibliographystyle{alphaurl}
\bibliography{bibliography}

\end{document}